\documentclass[draft]{amsart}

\usepackage{graphics}
\usepackage[all]{xy}

\usepackage{amssymb,amsmath}
\usepackage{psfrag}

\usepackage{stmaryrd}

\usepackage{amsfonts}
\usepackage{amscd}

\usepackage{graphicx}
\usepackage{epstopdf}
\usepackage{tikz}
\usetikzlibrary{arrows,automata}

\usepackage{graphicx,epsfig,float}

\usepackage{url}

\makeatletter
\def\url@leostyle{%
  \@ifundefined{selectfont}{\def\UrlFont{\sf}}{\def\UrlFont{\small\ttfamily}}}
\makeatother
\newtheorem{thm}{Theorem}[section]
\newtheorem{lem}[thm]{Lemma}
\newtheorem{prop}[thm]{Proposition}
\newtheorem{cor}[thm]{Corollary}
\newtheorem{fact}[thm]{Fact}

\newtheorem{clm}[thm]{Claim}

\newtheorem{defn}[thm]{Definition}

\newtheorem{nrmk}[thm]{Remark}
\newtheorem{expl}[thm]{Example}

\newcommand{\pf}{{\bf Proof. }}

\renewcommand{\tilde}{\widetilde}
\renewcommand{\bar}{\overline}

\newcommand{\ZZ}{\mathbb{Z}}

\newcommand{\tbA}{\tilde{{\mathbf A}}}

\newcommand{\bC}{{\mathbf C}}
\newcommand{\tbC}{\tilde{{\mathbf C}}}

\newcommand{\dex}{\textrm{\bf{\d{$\wr $}}}}

\renewcommand{\mod}{\mathrm{Mod}}

\newcommand{\cov}{\mathrm{Cov}}

\newcommand{\op}{\mathrm{Op}}
\newcommand{\opc}{\mathrm{Op}^{{\rm cons}}}

\newcommand{\df}{\mathrm{def}}
\newcommand{\Df}{\mathrm{Def}}
\newcommand{\Dfs}{\mathrm{Def}({\mathbb S})}
\newcommand{\tDf}{\tilde{\mathrm{Def}}}
\newcommand{\tDfs}{\tilde{\mathrm{Def}(\mathbb S)}}
\newcommand{\tDfw}{\bigwedge \textrm{}\tilde{\mathrm{Def}}}

\newcommand{\rh}{\mathit{R}\mathcal{H}\mathit{om}}
\newcommand{\ho}{\mathcal{H}\mathit{om}}

\newcommand{\Ho}{\mathrm{Hom}}
\newcommand{\Rh}{\mathrm{RHom}}

\newcommand{\tor}{\mathrm{Tor}}
\newcommand{\der}{\mathrm{D}}

\newcommand{\id}{\mathrm{id}}

\newcommand{\pt}{\mathrm{pt}}

\newcommand{\iso}{\stackrel{\sim}{\to}}

\newcommand{\supp}{\mathrm{supp}}

\newcommand{\imin}[1]{#1^{-1}}
\newcommand{\lind}[1]{\underset{#1}{\underrightarrow{\lim}}}
\newcommand{\Lind}{\underrightarrow{\lim}}  

\newcommand{\exs}[3]{0 \to {#1} \to {#2} \to {#3} \to 0}
\newcommand{\lexs}[3]{0 \to {#1} \to {#2} \to {#3}}

\begin{document}

\title {The six Grothendieck operations on o-minimal  sheaves}

\author {M\'{a}rio J. Edmundo}

\address{ Departamento de Matem\'atica\\
Faculdade de Ci\^encias da Universidade de Lisboa\\
Campo Grande, Edif\'icio C6, Piso 1\\
1749-016 Lisboa, Portugal}

\email{mjedmundo@fc.ul.pt}

\author{Luca Prelli}

\address{ Dipartimento di Matematica\\
Universit\`a degli Studi di Padova\\
via Trieste 63, 35121, Padova, Italia}

\email{lprelli@math.unipd.it}

\date{\today}
\thanks{The first author was supported by Funda\c{c}\~ao para a Ci\^encia e a Tecnologia, Financiamento Base 2008 - ISFL/1/209.  The second author  is a member of the Gruppo Nazionale per l'Analisi Matematica, la Probabilit\`a e le loro Applicazioni (GNAMPA) of the Istituto Nazionale di Alta Matematica (INdAM). This work is part of the FCT project PTDC/MAT/101740/2008.\newline
 }


\begin{abstract}
In this paper we develop  the  formalism of the Grothendieck six operations on o-minimal sheaves. The Grothendieck formalism allows us to obtain o-minimal versions  of: (i) derived projection formula; (ii) universal coefficient formula; (iii) derived base change formula;
(iv) K\"unneth formula; (v) local and global Verdier duality. 
\end{abstract}

\maketitle

\begin{section}{Introduction}\label{section intro}

The study of o-minimal structures (\cite{vdd}) is the analytic part of model theory which deals with theories of ordered, hence topological, structures satisfying certain tameness properties. It generalizes piecewise linear geometry (\cite[Chapter 1, \S 7]{vdd}), semi-algebraic geometry (\cite{BCR}) and globally sub-analytic geometry (\cite{kr}, also called finitely sub-analytic in \cite{vdd86}) and it is claimed to be the  formalization of  Grothendieck's notion of tame topology (topologie mod\'er\'ee). See \cite{vdd} and \cite{dmi}.

The most striking successes of this model-theoretic point of view of sub-analytic geometry include, on the one hand, an understanding of the behavior  at infinity of certain important classes of sub-analytic sets as in Wilkie's (\cite{w}) as pointed out by Bierstone and Milman  \cite{BM}, and on the other hand, the recent, somehow surprising,  first unconditional proof  of the Andr\'e-Oort conjecture for mixed Shimura varieties expressible as products of curves by Pila \cite{pl} following previous work also using o-minimality by Pila and Zannier (\cite{pz}), Pila and Wilkie (\cite{pw}) and  Peterzil and Starchenko (\cite{pst}).


The goal of this paper is to contribute further  to the claim that o-minimality does indeed realize Grothendieck's notion of topologie mod\'er\'ee by developing the formalism of the Grothendieck six operations on o-minimal sheaves, extending Delfs results for  semi-algebraic   sheaves (\cite{D3}) as well as Kashiwara and Schapira results for sub-analytic sheaves (defined and studied in \cite{ks2} and successively in \cite{lucap}) restricted to globally sub-analytic spaces   (as we deal, for now,  only with definable spaces and not locally definable spaces - the later more general case will be dealt with in a sequel to this paper). 

In the semi-algebraic case nearly all  of our results are completely new. Indeed, in \cite[Section 8]{D3}, Delfs constructs the semi-algebraic proper direct image functor and only proves two basic results about this functor (base change and commutativity with  small inductive limits) and then conjectures in \cite[Remark 8.11 ii)]{D3} that: ``It seems that the results of this section suffice to prove a semi-algebraic analogue of Verdier duality (by the same proof as in the theory of locally compact spaces c.f \cite{ver}).''

In the globally sub-analytic case we introduce   a new  globally sub-analytic proper direct image functor which, unlike  Kashiwara and Schapira (\cite{ks2}) sub-analytic proper direct image functor, generalizes to arbitrary o-minimal structures including: (i) arbitrary real closed fields; (ii)  the non-standard models 
 of sub-analytic geometry (\cite{dmm}); (iii) the non-standard o-minimal structure
which does not came from a standard one as in \cite{LpRob06,HrPet07}. See Remark \ref{nrmk comp ks1} for further details on this. Moreover, unlike \cite{ks2, lucap}, our definition is compatible with the restriction to open subsets (see Remark \ref{nrmk comp ks2}).\\

The Grothendieck formalism developed here allows us to obtain o-minimal versions  of: (i) derived projection formula;  (ii) derived base change formula; (iii) universal coefficient formula; (iv) K\"unneth formula; (v) local and the global Verdier duality.  It also sets up the framework for: (a) defining new o-minimal homology theories $H_*(X,{\mathbb Q}):=H^*(a_{X\dex }\circ a_X^{\dex}{\mathbb Q})$ extending o-minimal singular homology or the o-minimal Borel-Moore homology $H_*^{{\rm BM}}(X,{\mathbb Q}):=H^*(a_{X* }\circ a_X^{\dex}{\mathbb Q})$; (b)  the full  development of o-minimal geometry including the theory of characteristic classes, Hirzebruch-Riemann-Roch formula and Atiyah-Singer theorem in the non-standard o-minimal context.

The results of the paper were recently used: (i) To  settle  Pillay's conjecture for definably compact definable groups (\cite{p2} and \cite{hpp1})  
in arbitrary o-minimal structures. See \cite{oh5p}. 
This conjecture is an o-minimal analogue of Hilbert's fifth problem and it says roughly that a definably compact group has an infinitesimal normal subgroup such that the quotient equipped with a certain logic topology is a compact real Lie group of the same dimension. The solution to the conjecture is known to imply much information on the topology of definable groups as well as the structure of the definable subsets and, for example, explains how to obtain a compact real Lie group from an abelian variety over an arbitrary algebraically closed field of characteristic zero taking a quotient by a certain infinitesimal subgroup. 
(ii) In relation to the integral Hodge conjecture for real varieties by Benoist and Wittenberg (\cite{BWi1}, \cite{BWi2}).

We expect applications also in algebraic analysis. With our definition it is possible to treat the global sub-analytic sites in which are defined sheaves of functions with growth conditions up to infinity (e.g. tempered and Whitney $\mathcal{C}^{\infty }$ or holomorphic functions). This kind of objects are very important for applications as in \cite{Da13} and in \cite{DaK} (in which globally sub-analytic sheaves are hidden under the notion of ind-sheaves on a bordered space).\\

The structure of the paper is the following. In Section \ref{section prelim} we introduce the setting where we will work on, listing all the preliminary results that will be needed throughout  the rest of the paper. In Section \ref{section proper dir im} we define the proper direct image operation on o-minimal sheaves and prove its fundamental properties. In Section \ref{section pvd} we obtain the local and the global Verdier duality and its consequences.

Finally we note that, for the readers convenience, at the beginning of each section, after we introduce the set up, we include an extended summary in  which we  compare the ideas involved in the proofs of the section with those present in the references  cited in the bibliography. \\

The authors wish to thank the referee and also Olivier Wittenberg and Olivier Benoist for their very helpful comments on earlier versions of the paper. \\
\end{section}

\begin{section}{Preliminaries}\label{section  prelim}
In this paper we work in an  arbitrary o-minimal structure ${\mathbb M}=(M,<, (c)_{c\in {\mathcal C}}, $ $ (f)_{f\in  {\mathcal F}}, (R)_{R\in {\mathcal R}})$   with definable Skolem functions. We refer the reader to \cite{vdd} for basic o-minimality.

We let $\Df$ be the category whose objects are definable spaces  and whose morphisms are continuous definable maps between definable spaces - here and below ``definable'' always means definable in ${\mathbb M}$  with parameters. (See \cite[Chapter 10, \S 1]{vdd}). If ${\mathbb M}$ is a real closed field $(R,< , 0,1, +, \cdot )$, then $\Df$ is the category whose objects are semi-algebraic spaces over $R$  and whose morphisms are continuous semi-algebraic maps between such semi-algebraic space (\cite[Chapter 2]{vdd} and \cite[Chapter I, Example 1.1]{D3}); if ${\mathbb M}$ is ${\mathbb R}_{\rm an}=({\mathbb R},<, 0 ,1,+,\cdot, (f)_{f\in {\rm an}})$ - the field of real numbers expanded by restricted analytic functions, then $\Df$ is the category whose objects are globally sub-analytic spaces together with continuous  maps with globally sub-analytic graphs between such spaces (\cite{dd}); if ${\mathbb M}$ is an ordered vector space $(V, <, 0, +, (d)_{d\in D})$ over an ordered division ring $D$, then $\Df$ is the category whose objects are  the piecewise linear  spaces in this vector space together with continuous piecewise linear maps between such spaces (\cite[Chapter 1, \S 7]{vdd}).

Objects of $\Df$ are equipped with a topology determined by the order topology on $(M,<)$. However, if $(M,<)$ is non-archimedean  then infinite definable spaces are totally disconnected and not locally compact, so one studies definable spaces equipped with the o-minimal site and replaces topological notions (connected, normal, compact, proper) by their definable analogues (definably connected, definably normal, definably compact, definably proper). The o-minimal site (\cite{ejp}) generalizes  both the semi-algebraic site (\cite{D3}) and the sub-analytic site (\cite{ks2}). Given an object $X$ of $\Df$ the o-minimal site $X_{\df}$ on $X$ is the category $\op (X_{\df})$ whose objects are open (in the topology of $X$ mentioned above) definable subsets of $X$, the morphisms are the  inclusions and the admissible covers $\cov (U)$ of $U\in \op (X_{\df})$ are covers by open definable subsets of $X$ with  finite sub-covers.

As shown in \cite[Proposition 3.2]{ejp}, if $A$ is a commutative ring, then the category $\mod (A_{X_{\df}})$ of sheaves of $A$-modules on $X$ (relative to the o-minimal site) is isomorphic to the category $\mod (A_{\tilde{X}})$ of sheaves of $A$-modules on a certain spectral topological space $\tilde{X}$, the o-minimal spectrum of $X$, associated to $X$. The o-minimal spectrum $\tilde{X}$ of a definable space $X$ is the set of ultra-filters of definable subsets of $X$ (also know in model theory as types on $X$) equipped with the topology generated by the subsets $\tilde{U}$ with $U\in \op (X_{\df}).$ If $f:X\to Y$ is a morphism in $\Df$, then one has a corresponding continuous map $\tilde{f}: \tilde{X}\to \tilde{Y}: \alpha \mapsto \tilde{f}(\alpha )$ where $f(\alpha )$ is the ultrafilter in $\tilde{Y}$ determined by the collection $\{A: f^{-1}(A)\in \alpha \}.$ (See \cite[Definitions 2.2 and 2.18]{ejp} or \cite{c} and \cite{p} where these notions were first introduced).

If  ${\mathbb M}$ is a real closed field $(R,< , 0,1, +, \cdot )$,  and $V\subseteq R^n$ is an affine real algebraic variety over $R$, then $\tilde{V}$ is homeomorphic to ${\rm Sper}\,R[V]$, the real spectrum of the coordinate ring $R[V]$ of $V$ (\cite[Chapter 7, Section 7.2]{BCR} or \cite[Chapter I, Example 1.1]{D3}); and the isomorphism $\mod (A_{V_{\df}})\simeq \mod (A_{\tilde{V}})$ from \cite{ejp} corresponds in this case  to \cite[Chapter 1, Proposition 1.4]{D3}.

Below we denote by $\tDf$ the corresponding category of o-minimal spectra of definable spaces and continuous definable maps and
$$\Df \to \tDf$$
the functor just defined.  Due to the isomorphism $\mod (A_{X_{\df}})\simeq \mod (A_{\tilde{X}})$ for every object $X$ of $\Df$ we will work in this paper in $\tDf$. \\

\begin{subsection}{Extended summary}
In this section we will recall and prove some preliminary results that will be crucial for the rest of the paper. These preliminary results are necessary since the category $\tDf$ in which we will be working is rather different from the category ${\rm Top}$ of topological spaces.\\

An object of $\tDf$ is a $T_0$, quasi-compact and a spectral topological space, i.e.,   it has a basis of quasi-compact open subsets, closed under taking finite intersections and  each irreducible closed subset is the closure of a unique point. Such objects are not $T_1$ (unless they are finite), namely not every point is closed, so they are not  Hausdorff (unless they are finite). In particular, the usual (topological) definition of compact topological space cannot be used in $\tDf$. Similarly the usual (topological) definition of proper maps is of no use in $\tDf.$

Since $\tDf$ has cartesian squares (which one should point out are not cartesian squares in the category of topological spaces ${\rm Top}$), in this section we introduce a category theory definition of these notions in $\tDf$ just like in semi-algebraic geometry (\cite[Section 9]{dk1}) (and also  in algebraic geometry \cite[Chapter II, Section 4]{Har} or \cite[Chapter II, Section 5.4]{ega2}) and point out all the  properties needed later. \\

A fiber $\tilde{f}^{-1}(\alpha )$ of  a morphism $\tilde{f}:\tilde{X}\to \tilde{Y}$ in $\tDf$ is  not in general  an object of $\tDf$, but following \cite[Lemma 3.1]{bera} in the affine case, such fiber is homeomorphic to an object $\tilde{(f^{{\mathbb S}})^{-1}(a)}$ of $\tDfs$ where $a$ is a realization of the type $\alpha $ and ${\mathbb S}$ is the prime model of the first-order theory of ${\mathbb M}$ over $\{a\}\cup M.$ Here $\tDfs$ is the same as $\tDf$ but defined in the o-minimal structure ${\mathbb S}$. This phenomena is the model theoretic analogue of what happens in real-algebraic geometry (\cite[Proposition 2.4]{D3}, \cite[Chapter II, 3.2]{S2}) (and also in algebraic geometry \cite[Chapter II, Section 3]{Har}). Indeed, if $f:X\to Y$ is a morphism of real schemes over a real closed field $R$ and $\alpha \in Y$, then $f^{-1}(\alpha )$ with the underlying topology is homeomorphic  to the real scheme  $X\times _Y{\rm Sper} \, k(\alpha )$ over the residue ordered field  $k(\alpha )$ of $\alpha $ (recall that $\alpha $ is a prime cone). By \cite{p} the real closure of $k(\alpha )$ is isomorphic to the prime model over $R$ and a realization of $\alpha .$

Due to the homeomorphism $\tilde{f}^{-1}(\alpha )\simeq \tilde{(f^{{\mathbb S}})^{-1}(a)},$ when $\alpha \in \tilde{Y}$ is a closed point, we will be able to use the theory from \cite[Section 3]{ep1} of normal and constructible families of supports on the object $\tilde{(f^{{\mathbb S}})^{-1}(a)}$  also on the fiber $\tilde{f}^{-1}(\alpha )$ after we show that working in certain  full subcategories $\tbA$ of $\tDf$ the family $c$ of complete supports on $\tilde{(f^{{\mathbb S}})^{-1}(a)}$  is normal and constructible. Normal and constructible families of supports are the o-minimal analogue of the semi-algebraic and paracompactifying families of supports in semi-algebraic geometry (\cite[Chapter II, Section 1]{D3}). \\

\end{subsection}

\begin{subsection}{Morphisms proper in $\tDf$}\label{subsection proper in tdef}
Here we will introduce a category theory definition of the notions proper and complete in $\tDf$ just like in semi-algebraic geometry (\cite[Section 9]{dk1}) (and also  in algebraic geometry \cite[Chapter II, Section 4]{Har} or \cite[Chapter II, Section 5.4]{ega2}) and point out the  properties needed later. Note that these notions and properties are noting but the corresponding notions and properties in $\Df,$ introduced and studied in \cite{emp}, under the isomorphism $\Df\to \tDf.$\\

If $X$ is an object of $\tDf$, then a subset $Z\subseteq X$ is called {\it constructible} if $Z$ is also an object of $\tDf.$\\

Let $f:X\to Y$ be a morphism in $\tDf.$   We say that:
\begin{itemize}
\item
$f:X\to Y$ is  {\it closed in $\tDf$} if for every  closed  constructible subset  $A$ of $X$, its image $f(A)$ is a closed (constructible) subset of $Y.$
\item
$f:X\to Y$ is  a {\it closed (resp. open)  immersion} if $f:X\to f(X)$ is a  homeomorphism and $f(X)$ is  a closed (resp. open) subset of $Y$.\\
\end{itemize}

Since $\Df \to \tDf$ is an isomorphism of categories and cartesian squares exist in $\Df$ we have:\\

\begin{fact}
{\em
In the category $\tDf$  the {\it cartesian square} of any two morphisms $f:X\to Z$ and $g:Y\to Z$ in $\tDf$ exists and is given by a commutative diagram
$$
\xymatrix{X\times _ZY \ar[r]^{p_Y} \ar[d]^{p_X} & Y \ar[d]^g \\
X \ar[r]^f &Z}
$$
where the morphisms $p_X$ and $p_Y$ are known as projections. The Cartesian square satisfies the following universal property:   for any other object $Q$ of $\tDf$ and morphisms $q_X:Q\to X$ and $q_Y:Q\to Y$ of $\tDf$ for which the following diagram commutes,
$$
\xymatrix{
Q \ar@/_/[ddr]_{q_X} \ar@/^/[drr]^{q_Y}
\ar@{-->}[dr]^-{u} \\
& X\times _ZY \ar[d]^{p_X} \ar[r]^{p_Y}
& Y \ar[d]^g \\
& X \ar[r]^f & Z }
$$
there  exist a unique natural morphism $u:Q\to X\times _ZY$ (called mediating morphism) making the whole diagram commute. As with all universal constructions, the cartesian square  is unique up to a definable homeomorphism.\\
}
\end{fact}

The following is a very important remark that one should always have in mind:

\begin{nrmk}
{\em
The cartesian square in $\tDf$ of two morphisms $f:X\to Z$ and $g:Y\to Z$ in $\tDf$ is not the same as the cartesian square in ${\rm Top}$ of the two morphisms $f:X\to Z$ and $g:Y\to Z$ in ${\rm Top}$.  In particular, if ${\rm pt}$ denotes a fixed one point object of $\tDf$, then the {\it cartesian product}
$$
\xymatrix{X\times _{{\rm pt}}Y \ar[r]^{p_Y} \ar[d]^{p_X} & Y \ar[d]^{a_Y} \\
X \ar[r]^{a_X} &{\rm pt}}
$$
 in $\tDf$, also denoted by $X\times Y$, of two objects $X$ and $Y$ in $\tDf$ is not the same as usual the cartesian product  (in ${\rm Top}$) of the objects  $X$ and $Y$ in ${\rm Top}$.\\
}
\end{nrmk}

Given a morphism $f:X\to Y$ in $\tDf$, the corresponding diagonal morphism is the unique morphism $\Delta :X\to X\times _YX$ in $\tDf$ given by the universal property of cartesian squares:
$$
\xymatrix{
X \ar@/_/[ddr]_{{\rm id}_X} \ar@/^/[drr]^{{\rm id}_X}
\ar@{-->}[dr]^-{\Delta } \\
& X\times _YX \ar[d]^{p_X} \ar[r]^{p_X}
& X \ar[d]^f \\
& X \ar[r]^f & Y. }
$$
We say that:
\begin{itemize}
\item
 $f:X\to Y$ is {\it separated in $\tDf$} if the corresponding diagonal morphism  $\Delta :X\to X\times _YX$ is a closed  immersion. \\
 \end{itemize}

 \noindent
 We say that an object $Z$ in $\tDf$ is {\it separated in $\tDf$} if the morphism $Z\to {\rm pt}$ to a point is separated.\\



Let $f:X\to Y$ be a morphism in $\tDf.$   We say that:
\begin{itemize}
\item
$f:X\to Y$ is {\it universally  closed in $\tDf$} if for any morphism $g:Y'\to Y$  in $\tDf$ the morphism $f':X'\to Y'$ in $\tDf$ obtained from the cartesian square
$$
\xymatrix{X' \ar[r]^{f'} \ar[d]^{g'} & Y' \ar[d]^g \\
X \ar[r]^f &Y}
$$
in $\tDf$ is   closed in $\tDf$.  \\
\end{itemize}

\begin{defn}\label{defn proper in tdf}
{\em
We say that a morphism $f:X\to Y$ in $\tDf$ is {\it proper in $\tDf$} if $f:X\to Y$ is separated and universally closed  in $\tDf$.\\
}
\end{defn}

\begin{defn}\label{defn complete in tdf}
{\em
We say that an object $Z$ of $\tDf$ is {\it complete in $\tDf$} if the morphism $Z\to {\rm pt}$  is proper in $\tDf.$ \\
}
\end{defn}

Directly from the definitions (as in \cite[Chapter II, Proposition 5.4.2 and Corollary 5.4.3]{ega2}, see also \cite[Section 9]{dk1}) one has the  following. See  \cite[Proposition 3.7]{emp} for a detailed proof in $\Df$ which transfers to $\tDf$ due to the isomorphism $\Df\to \tDf$:\\

\begin{prop}\label{prop proper in tdef}
In the category $\tDf$ the following hold:
\begin{enumerate}
\item
Closed  immersions are proper in $\tDf$.
\item
A composition of two morphisms  proper in $\tDf$   is proper in $\tDf$.
\item
Let $
\xymatrix{ X \ar [rd]_p \ar[r]^f & Y \ar[d]^q \\
& Z}
$
be a morphism over $Z$ in $\tDf$ and $Z'\to Z$ a base extension in $\tDf.$
If $f:X\to Y$ is  proper  in $\tDf,$  then the corresponding base extension morphism $f':X\times _ZZ'\to Y\times _ZZ'$ is proper in $\tDf$.
\item
Let $
\xymatrix{ X \ar [rd]_p \ar[r]^f & Y \ar[d]^q \\
& Z}
$
and
$
\xymatrix{ X' \ar [rd]_{p'} \ar[r]^{f'} & Y' \ar[d]^{q'} \\
& Z}
$
be morphisms over $Z$ in $\tDf.$
If $f:X\to Y$ and $f':X'\to Y'$ are proper in $\tDf$, then the corresponding product morphism $f\times f':X\times _ZX' \to Y\times _ZY'$   is proper in $\tDf$.
\item
If $f:X\to Y$ and $g:Y\to Z$ are morphisms such that $g\circ f$ is proper in $\tDf$, then:
\begin{itemize}
\item[(i)]
$f$ is proper in $\tDf;$
 \item[(ii)]
 if $g$ is separated in $\tDf$ and $f$ is surjective, then $g$ proper in $\tDf.$
 \end{itemize}
\item
A morphism $f:X\to Y$  is proper in $\tDf$ if and only if $Y$ can be covered by finitely many open definable subsets $V_i$ such that $f_|:f^{-1}(V_i)\to V_i$ is proper in $\tDf$.\\
\end{enumerate}
\end{prop}

From Proposition \ref{prop proper in tdef} we easily have:

\begin{cor}\label{cor basic compl}
Let $f:X\to Y$ be a morphism in $\tDf$ and $Z\subseteq X$ a complete object of $\tDf$. Then the following hold:
\begin{enumerate}
\item
$Z$ is a closed (constructible) subset of $X.$
\item
$f_{|Z}:Z\to Y$ is proper in $\tDf.$
\item
$f(Z)\subseteq Y$ is  (constructible) complete in $\tDf.$
\item
If $f:X\to Y$ is proper in $\tDf$ and $C\subseteq Y$ is a complete object of $\tDf$, then $f^{-1}(C)\subseteq X$ is (constructible) complete in $\tDf.$ \\
\end{enumerate}
\end{cor}

Below, if $\bC$ is a subcategory of $\Df$ we denote by $\tbC$ its image under $\Df \to \tDf.$ The following is a standard consequence of  Proposition \ref{prop proper in tdef}. See \cite[Corollary 3.9]{emp} for a detailed proof in $\Df$ which transfers to $\tDf$ due to the isomorphism $\Df\to \tDf$:\\

\begin{cor}\label{cor completion in C}
Let $\bC$ be a full a subcategory  of the category of definable spaces $\Df$ whose set of objects is:
\begin{itemize}
\item
closed under taking locally closed definable subspaces of objects of $\bC$,
\item
closed under taking cartesian products of objects of $\bC$,
\end{itemize}
Then the following are equivalent:
\begin{enumerate}
\item
Every object $X$ of $\tbC$ is {\it completable in $\tbC$}  i.e., there exists an object $X'$ of $\tbC$ which is complete in $\tDf$ together with an open immersion $i: X\hookrightarrow X'$ in $\tbC$ with  $i(X)$ dense in $X'$. Such $i: X\hookrightarrow X'$ is called a {\it completion} of $X$ in $\tbC$.

\item
Every  morphism $f:X\to Y$ in $\tbC$ is {\it  completable in $\tbC$} i..e,  there exists a commutative diagram
$$
\xymatrix{ X \ar [d]_f \ar[r]^-{i} & X' \ar[d]^{f'} \\
Y \ar[r]^-{j} & Y'}
$$
of  morphisms in $\tbC$ such that:  (i) $i:X\to X'$ is a completion of $X$  in $\tbC$; (ii) $j$ is a  completion of $Y$ in $\tbC$.
\item
Every morphism $f:X\to Y$ in $\tbC$ has {\it a proper extension in $\tbC$} i.e.,  there exists a commutative diagram
$$
\xymatrix{ X \ar [rd]_f \ar[r]^\iota & P \ar[d]^{\overline{f}} \\
& Y}
$$
of morphisms in $\tbC$ such that $\iota$ is a  open immersion with $\iota (X)$ dense in $P$ and $\overline{f}$ a  proper  in $\tDf$.\\
\end{enumerate}
\end{cor}

\begin{nrmk}\label{nrmk P}
{\em
Given a morphism $f:X\to Y$ in $\tbC,$ in Corollary \ref{cor completion in C} (3), we have $P=f'^{-1}(j(Y))\subseteq X',$ an open subset of  $X',$  where $i:X\to X'$ is a completion of $X$ in $\tbC,$ $j:Y\to Y'$ is a  completion of $Y$ in $\tbC$ and $f'$ is such that 
$$
\xymatrix{ X \ar [d]_f \ar[r]^-{i} & X' \ar[d]^{f'} \\
Y \ar[r]^-{j} & Y'}
$$
is a commutative diagram of  morphisms in $\tbC.$ Moreover, $\bar{h}=j^{-1}\circ f'_{|P}:P\to Y$ where $j^{-1}:j(Y)\to Y$ is the inverse of $j:Y\to j(Y)$ and $\iota = i:X\to P\subseteq X'.$

In particular, if $\pi :X\times Y\to Y$ is a morphism of $\tbC,$ then in the commutative diagram 
$$
\xymatrix{ X\times Y \ar [rd]_{\pi } \ar[r]^\iota & P \ar[d]^{\overline{\pi}} \\
& Y}
$$
of morphisms in $\tbC$ we have $P=X'\times Y,$ $\iota =i\times {\rm id}$ (where $i:X\to X'$ is the completion of $X$ in $\tbC$) and $\bar{\pi }:P\to Y$ is the projection onto $Y.$\\
}
\end{nrmk}

Finally we observe the following which was proved for $\Df$  and $\Dfs$ in \cite[Theorem 4.3]{emp}, and transfers to $\tDf$ and $\tDfs$ due to the isomorphisms  $\Df\to \tDf$ and $\Dfs\to \tDfs$:\\

\begin{prop} \label{prop sep and proper in s}
Let ${\mathbb S}$ be an elementary extension of ${\mathbb M}.$ Let $f:X\to Y$ a morphism in $\tDf.$ Then the following are equivalent:
\begin{enumerate}
\item
$f$ is separated (resp. proper) in $\tDf .$
\item
$f^{{\mathbb S}}$ is separated (resp. proper) in $\tDfs.$ 
\end{enumerate}
In particular, $X$ is complete in $\tDf$ if and only if $X({\mathbb S})$ is complete in $\tDfs.$\\
\end{prop}

We end the subsection by  describing the closed subsets of objects of $\tDf$ and also showing that  a morphism in $\tDf$ is closed if and only if it closed in $\tDf.$\\

\begin{nrmk}\label{nrmk closed in qc}
{\em
Let $X$ be an object of $\tDf$ and $C=\bigcap _{i\in I}C_i\subseteq X$ with the $C_i$'s constructible subsets. The following hold:
\begin{itemize}
\item
$C$ is quasi-compact subset  of $X$ (with the induced topology).
\item
A closed subset $B$ of  $C$ is a quasi-compact subset of $X.$
Indeed, $B=D\cap C$ with $D$ a closed subset of $X$. So $D$ is quasi-compact
and $B$ is also quasi-compact.

\item
A quasi-compact subset $B$ of $C$ is a quasi-compact subset of $X.$\\
\end{itemize}
}
\end{nrmk}

\begin{lem}\label{lem spec f qc}
Let $f:X\to Y$ be a morphism in $\tDf$. Then the following hold:
\begin{enumerate}
\item
If $C=\bigcap _{i\in I}C_i\subseteq X$ with the $C_i$'s constructible subsets, then $f(C)=\bigcap _{i\in I}f(C_i).$
\item
If $D=\bigcap _{j\in J}D_j\subseteq Y$ with the $D_j$'s constructible subsets, then $f^{-1}(D)=\bigcap _{j\in J}f^{-1}(D_j)$.
\item
$C\subseteq X$ is closed if and only if $C=\bigcap _{i\in I} C_i$ with each $C_i\subseteq X$  constructible and closed.
\end{enumerate}
\end{lem}

\pf
(1)  Suppose that $C=\bigcap _{i\in I}C_i$ with the  $C_i$'s   constructible. Then $C$ is closed and hence compact in the Stone topology of $X$ (the topology generated by the constructible subsets). By \cite[Remark 2.19]{ejp}, $f:X\to Y$ is continuous with respect to the  Stone topologies on $X$ and $Y$, so $f(C)$ is compact in the Stone topology of $X$. Since $Y$ is Hausdorff in the Stone topology, $f(C)$ is closed (in the Stone topology). Therefore,
$$f(C)=\bigcap \{E: f(C)\subseteq E\,\,\textrm{and $E$ is  constructible}\}.$$
If $E$ is a constructible subset such that $f(C)\subseteq E$, then $C\subseteq f^{-1}(E)$ and by compactness, there is $C_i$ such that  $C\subseteq C_i\subseteq f^{-1}(E)$ and hence $f(C)\subseteq f(C_i)\subseteq E$. Therefore, $f(C)=\bigcap _{i\in I}f(C_i)$.

(2) Obvious.

(3) Suppose that $C$ is closed. Then its complement is open and so a union of open constructible subsets. Then   $C$ is the intersection of closed constructible subsets (\cite[Remark 2.3]{ejp}). The converse is clear.
\qed \\

\begin{prop}\label{prop closed on cons}
Let $f:X\to Y$ be a morphism in $\tDf$.  Then the following are equivalent:
\begin{enumerate}
\item
$f$ is closed.
\item
$f$ is closed in $\tDf.$ 
\end{enumerate}
\end{prop}

\pf
Assuming (1) clearly  (2) is immediate. Assume (2). If $D\subseteq X$ be  a closed constructible subset, then $f(D)$ is closed (and constructible by \cite[Remark 2.19]{ejp}); if $C\subseteq X$ is a general closed subset, we have $C=\bigcap _{i\in I}C_i$ with the  $C_i$'s  closed and constructible (Lemma \ref{lem spec f qc} (3)) and by Lemma \ref{lem spec f qc} (1),  $f(C)=\bigcap _{i\in I}f(C_i)$ with  each $f(C_i)$  closed and constructible and hence, $f(C)$ is closed as required. \qed \\


\end{subsection}

\begin{subsection}{Definably proper, definably compact and definably normal}
Here we recall the relation between the notions of proper, complete and normal in $\tDf$ and notions definably proper, definably compact and definably normal.\\

Recall from \cite{ps} that if  $X$ is an object of $\Df$ and  $C\subseteq X$ is a definable subset (i.e. it is also in $\Df$), then we say that $C$ is {\it definably compact} if  for every morphism  $\alpha :(a,b)\to C\subseteq X$ in $\Df,$ where $a<b$ are in $M\cup \{-\infty, +\infty \},$  the limits $\lim _{t\to a^+}\alpha (t)$ and $\lim _{t\to b^-}\alpha (t)$ exist in $C.$\\

For  definable subsets $X\subseteq M^n$ with their induced topology we have (\cite[Theorem 2.1]{ps}):

\begin{fact}\label{fact def comp affine}
A definable subset $X\subseteq M^n$ is definably compact if and only if it is closed and bounded in $M^n$  \\ 
\end{fact}

Recall also that a morphism $f:X\to Y$ in $\Df$  is called {\it definably proper} if for every definably compact definable subset $K$ of $Y$ its inverse image $f^{-1}(K)$ is a definably compact definable subset of $X$. (See  \cite[Chapter 6, Section 4]{vdd}).\\

From the definitions we see that:

\begin{nrmk}\label{nrmk proper compact}
{\em
An object $X$ of $\Df$  is definably compact if and only if the morphism $X\to \pt$ in $\Df$  is definably proper.\\
}
\end{nrmk}

Recall that we have separated, proper or complete  in $\Df$ if and only if under the isomorphim $\Df \to \tDf$ we have  separated, proper or complete  in $\tDf.$ \\

From  the way cartesian squares are defined in $\Df$ we easily obtain the following:

\begin{nrmk}\label{nrmk sep df2}
{\em
Let $f:X\to Y$ be a morphism in $\Df$. 
Then the  following are equivalent:
\begin{enumerate}
\item
$f:X\to Y$  is separated in $\Df$.
\item
The fibers $f^{-1}(y)$ of $f$ are Hausdorff (with the induced topology). 
\end{enumerate}
In particular, $X$ is separated in $\Df$ if and only if $X$ is Hausdorff.\\
}
\end{nrmk}

Under our  assumption that ${\mathbb M}$ has definable Skolem functions we have (see \cite[Theorem 3.15 and Corollary 3.17]{emp}):

\begin{fact}\label{fact def proper2}
 Let $X$ and $Y$ be objects of $\Df$ such that $X$ and $Y$ are  Hausdorff and $Y$ locally definably compact (i.e., every $y\in Y$ has a definably compact neighborhood). 
Let $f:X\to Y$ be a morphim in $\Df$. Then the following are equivalent:
\begin{enumerate}
 \item
$f$ is  proper in $\Df$. 
\item
$f$ is definably proper.
\end{enumerate}
In particular, $X$ is Hausdorff and definably compact if and only if $X$ is complete in $\Df.$\\
\end{fact}

Recall that an object $X$  of $\Df$ is \textit{definably normal} if one of the following equivalent conditions holds:
\begin{enumerate}
\item
for every disjoint closed definable subsets $Z_1$ and $Z_2$ of $X$ there are disjoint open definable subsets $U_1$ and $U_2$ of $X$ such that $Z_i\subseteq U_i$ for $i=1,2.$

\item
for every  $S\subseteq X$ closed definable  and $W\subseteq X$ open definable such that $S\subseteq W$, there is an open definable subsets $U$  of $X$ such that $S\subseteq U$ and $\bar{U}\subseteq W$.\\
\end{enumerate}

Under our  assumption that ${\mathbb M}$ has definable Skolem functions we have (see \cite[Theorem 2.11]{emp}):

\begin{fact}\label{fact def comp is normal}
If  $X$ is an object of $\Df$ which is  Hausdorff and definably compact (i.e complete in $\Df$), then $X$ is definably normal.\\
\end{fact}

The relationship between definably normal (in $\Df$) and normal in $\tDf$ (as spectral topological spaces) is given by the following. See \cite[Proposition 2.12 and Theorem 2.13]{ejp}.

\begin{fact}\label{fact ominspec closed point}
If $X$ is an object of $\Df,$ then the following are equivalent:
\begin{enumerate}
\item $\tilde{X}$ is normal. In fact, if $F$ and $G$ are two disjoint closed
subsets of $\tilde{X}$ then there
exist two disjoint constructible open subsets $U$ and $V$ of $\tilde{X}$ such
that $F\subseteq U$ and $G\subseteq V$.
\item
Every point $\alpha \in \tilde{X}$ has a unique closed specialization.
\item 
$X$ is definably normal.\\ 
\end{enumerate}
\end{fact}

Definable normality gives the shrinking lemma (see \cite[Chapter 6, (3.6)]{vdd} in $\Df$ or \cite[Proposition 2.17]{ejp} in $\tDf$):\\

\begin{fact}[The shrinking lemma]\label{fact shrinking lemma}
Let $X$ be an object of $\Df$ which is definably normal. If $\{U_i:i=1,\dots ,n\}$ is a covering of $X$ by open definable subsets, then there are definable open subsets $V_i$ and definable closed subsets $C_i$ of $X$ ($1\leq i\leq n$) such that $V_i\subseteq C_i\subseteq U_i$ and $\{V_i:i=1,\dots, n\}$ is a covering of $X.$ \\
\end{fact}

We end this subsection with the following result: \footnote{In topology, Borsuk asked in 1937 if the product of a normal space and the unit interval is normal and the question was only solved, negatively, without any set theoretic conditions beyond the axiom of choice,  in 1971 by M. E. Rudin (\cite{Rud}). We thank the referee for pointing us to this. So Theorem \ref{thm CxN is normal} is yet another manifestation of tameness of o-minimal structures.}

\begin{thm}\label{thm CxN is normal}
Let $Z$ and $K$ be objects of $\Df$ with $Z$ definably normal and $K$ Hausdorff and definably compact (i.e. complete in $\Df$). Then $Z\times K$ is definably normal.
\end{thm}

This theorem is a consequence of the following proposition, the affine case of the theorem, together with the shrinking lemma and Facts \ref{fact def comp is normal} and \ref{fact ominspec closed point}.\\

\begin{prop}\label{prop CxN affine is normal}
Let $Z\subseteq M^n$ be an affine definably normal, definable subset and let $K\subseteq M^k$ be a closed and bounded definable subset. Then $Z\times K\subseteq M^{n+k}$ is definably normal. In particular, every closed and bounded definable subset of $M^k$ is definably normal.\\
\end{prop}

\noindent
{\bf Proof of Theorem \ref{thm CxN is normal}:} By Fact \ref{fact def comp is normal} $K$ is also definably normal. So by the shrinking lemma there are $Z_1,\ldots, Z_l\subseteq Z$ closed affine definable subsets such that $Z=\bigcup _{i}Z_i$ and there are $K_1, \ldots , K_n\subseteq K$ closed affine definable subsets such that $K=\bigcup _{j}K_j.$ Now $Z\times K=\bigcup _{i,j}Z_i\times K_j$ with each $Z_i\times K_j$ a closed definable subset of $Z\times K$ and each $K_j$  definably compact. So each $Z_i\times K_j$ is definably normal by Fact \ref{fact def comp affine} and Proposition \ref{prop CxN affine is normal}. Now going to $\tDf$ but omitting the tilde, if $\alpha \in Z\times K,$ then $\alpha \in Z_i\times K_j$ for some $i$ and $j,$ and so, by Fact \ref{fact ominspec closed point}, $\alpha $ has a unique closed specialization $\rho $ in $Z_i\times K_j.$ Since $Z_i\times K_j$ is closed in $Z\times K,$ every specialization of $\alpha $ in $Z\times K$ is in $Z_i\times K_j$ (\cite[Proposition 2.7]{ejp}) and $\rho $ is a closed specialization of $\alpha $ in $Z\times K.$ So $\rho $ is the unique closed specialization of $\alpha $ in $Z\times K.$ By Fact \ref{fact ominspec closed point}, the tilde of $Z\times K$ in $\tDf$ is normal and $Z\times K$ (in $\Df$) is definably normal. \qed\\

Proposition \ref{prop CxN affine is normal} will be obtained after a couple of lemmas showing special cases of it.  \\

Below we will use the usual notation 
$$\Gamma (h):=\{(x,y)\in B\times M: h(x)=y\},$$  
$$(f,g)_B:=\{(x,y)\in B\times M: f(x)<y<g(x)\}$$
 and 
 $$[f,g]_B:=\{(x,y)\in B\times M: f(x)\leq y\leq g(x)\}$$ 
 where  $h, f, g :B\to M$ are functions.\\

The following is obtained from the definition of cells (\cite[Chapter 3, \S2]{vdd}):\\

\begin{nrmk}\label{nrmk cells perm}
{\em
Let $C\subseteq M^n$ be a $d$-dimensional cell. Then by the definition of cells,  $C$ is a $(i_1, \ldots , i_n)$-cell for some unique sequence $(i_1, \ldots , i_n)$ of $0$'s and $1$'s.  Moreover, if     $\lambda (1)< \dots <\lambda (d)$  are the indices $\lambda \in \{1, \ldots , n\}$ for which $i_{\lambda }=1$ and 
$$p_{\lambda (1), \dots , \lambda (d)}:M^n\to M^d:(x_1, \ldots , x_n)\mapsto (x_{\lambda (1)}, \ldots , x_{\lambda (d)})$$ 
is the projection, then $C':=p_{\lambda (1), \dots , \lambda (d)}(C)$ is an open $d$-dimensional cell in $M^d$ and the restriction $p_C:=p_{\lambda (1), \dots , \lambda (d)|C}:C\to C'$ is a definable homeomorphism  (\cite[Chapter 3, (2.7)]{vdd}). 

Let $\tau (1)<\dots <\tau (n-d)$ be the indices $\tau \in \{1, \ldots , n\}$ for which $i_{\tau }=0.$ For each such $\tau $, by the definition of cells, there is a definable continuous function $h_{\tau }: \pi _{\tau -1}(C)\subseteq M^{\tau -1}\to M$ where, for each $k=1, \ldots , n$,  $\pi _{k}:M^n\to M^{k}$ is the projection onto the first $k$-coordinates. Moreover we have $\pi _{\tau }(C)=\{(x,h_{\tau }(x)): x\in \pi _{\tau -1}(C)\}.$

Let $f=(f_1, \ldots , f_{n-d}):C'\to M^{n-d}$ be the definable continuous map where for each $l=1, \ldots , n-d$ we set $f_l=h_{\tau (l)}\circ \pi _{\tau (l) -1}\circ p_C^{-1}.$ Let  $\sigma :M^n\to M^n:(x_1,\ldots , x_n)\mapsto (x_{\lambda (1)}, \ldots , x_{\lambda (d)}, x_{\tau (1)}, \ldots , x_{\tau (n-d)}).$ Then we clearly have 
$$\sigma (C)=\left\{\left(x,f(x)\right) : x\in C'\right\}.$$

Now let $h:C\to M$ be any continuous definable map. Then the definable set $W= (p_{\lambda (1), \dots , \lambda (d)})^{-1}(C')$  is an open definable neighborhood of $C$ in $M^n$ and $g:W\to M$ given by $g(z)=h(\sigma ^{-1}(p_{\lambda (1), \dots , \lambda (d)}(z), f(p_{\lambda (1), \dots , \lambda (d)}(z))))$ is  a continuous definable map such that $g_{|C}=h.$ \\
}
\end{nrmk}

\begin{lem}\label{lem closed fiber in open}
Let $Z\subseteq M^n$ be a definably normal definable subset. Let $S\subseteq Z\times [a,b]$ be a closed definable subset and $W\subseteq Z\times [a,b]$ an open definable subset. Then for every closed definable subset $F\subseteq \pi (S)$ such that $S\cap \pi ^{-1}(F)\subseteq W$ there is an open definable neighborhood $O$ of $F$ in $Z$ such that $O\subseteq \bar{O}\cap Z\subset \pi (W)$ and  $S\cap \pi ^{-1}(\bar{O}\cap Z)\subseteq W.$
\end{lem}

\pf
Let $\pi :Z\times [a,b]\to Z$ be the projection and let $\pi ':Z\times [a,b]\to [a,b]$ be the other projection.  Since $[a,b]$ is definably compact, it is complete in $\Df$ (Fact \ref{fact def proper2}). In particular, the projection $\pi :Z\times [a,b]\to Z$ is closed in $\Df.$ 

Let  $W^c=(Z\times [a,b])\setminus W.$ If $S\subseteq W$ then since $\pi (S) \subseteq \pi (W)$ is closed in $Z$ and $Z$ is definably normal, there is an open definable neighborhood $O$ of $\pi (S)\supseteq F$ in $Z$ such that $O\subseteq \bar{O}\cap Z\subset \pi (W)$ and  so $S\cap \pi ^{-1}(\bar{O}\cap Z)=S\subseteq W.$ So we may  suppose that $S\cap W^c\neq \emptyset .$

For $z\in Z$ let 
$$D(z)=\{((d_1^-,d_1^+),\ldots ,(d^-_n,d^+_n)\in M^{2n}: z\in \Pi _{i=1}^n(d^-_i,d^+_i)\cap Z\}$$
and for $d\in D(z)$ let
$$U(z,d)=\Pi _{i=1}^n(d^-_i,d^+_i)\cap Z.$$
Then $\{U(z,d)\}_{d\in D(z)}$ is a uniformly definable system of fundamental open neighborhoods of $z$ in $Z.$ Moreover, since the relation $d\preceq d'$ on $D(z)$ given by $U(z,d)\subseteq U(z,d')$ is a definable downwards directed order on $D(z),$ by \cite[Lemma 4.2.18]{HrLo} (or \cite[Lemma 2.19]{Hr04}), there is a definable type $\beta $ on $D(z)$ such that for every $d\in D(z)$ we have $\{d' \in D(z): d'\preceq d\}\in \beta .$ \footnote{Recall that a  type $\beta $ on $X$ (i.e. an ultrafilter of definable subsets of $X$) is a {\it definable type on $X$} if  and only if for every uniformly definable family $\{Y_t\}_{t\in T}$ of definable subsets of $X,$ the set $\{t\in T: Y_t\in \beta \}$ is a definable set.}

Suppose that $z\in F$ and  for all $d\in D(z),$ we have $(S\cap \pi ^{-1}(U(z,d)))\cap W^c\neq \emptyset .$  Then by definable Skolem functions, there is a definable map $$h:D(z)\to S\cap W^c\subseteq Z\times [a,b]$$
 such that for every $d\in D(z)$ we have $h(d)\in (S\cap \pi ^{-1}(U(z,d)))\cap W^c.$ 
 
 Let $\alpha $ be the definable type on $S\cap W^c$ determined by the collection $\{A\subseteq S\cap W^c: h^{-1}(A)\in \beta \}.$ Let $\alpha _1$ be the  definable type on $\pi (S)$ determined by the collection $\{A\subseteq \pi (S): \pi ^{-1}(A)\in \alpha \}$ and let $\alpha _2$ be the  definable type on $[a,b]$ determined by the collection $\{A\subseteq [a,b]: \pi '^{-1}(A)\in \alpha \}.$

We have that $z$ is the limit of $\alpha _1$ i.e., for every open definable subset  $V$ of $Z$ such that $z\in V$ we have $V\in \alpha _1.$ Indeed, given any such $V$ there is $d'$ such that $U(z,d')\subseteq V$ and 
\begin{eqnarray*}
h^{-1}(\pi ^{-1}(V))&\supseteq &h^{-1}(\pi ^{-1}(\pi (S)\cap U(z,d')))\\
&\supseteq & h^{-1}((S\cap \pi ^{-1}(U(z,d')))\cap W^c)\\
&\supseteq & \{d''\in D(z):d''\preceq d'\}.
\end{eqnarray*}
On the other hand, since $\alpha _2$ is a definable type on the closed and bounded definable set $[a,b],$ it  has a limit, say  $c\in [a,b]$ (\cite[Fact 5.1]{emp}). It follows that  $(z,c)\in Z\times [a,b]$ is the limit of $\alpha .$ Since $S\cap W^c$ is closed and $\alpha $ is a definable type on $S\cap W^c,$ its limit  $(z,c)$ is in $S\cap W^c.$ But then $(z,c)\in S\cap \pi ^{-1}(z)\subseteq W^c$ which is a contradiction. 

So for each $z\in F$ there is $d\in D(z)$ such that $S\cap \pi ^{-1}(U(z,d))\subseteq W.$ By definable Skolem functions there is a definable map $\epsilon :F\to M^{2n}$ such that for each $z\in F$ we have $\epsilon (z)\in D(z)$ and $S\cap \pi ^{-1}(U(z,\epsilon (z)))\subseteq W.$  Then 
$$U(F,\epsilon )=\bigcup _{z\in F}U(z,\epsilon (z))$$
 is an open definable neighborhood of $F$ in $Z$ such that 
$$S\cap \pi ^{-1}(U(F,\epsilon )) =\bigcup _{z\in F}S\cap \pi ^{-1}(U(z,\epsilon (z)))\subseteq W.$$

Since $U(F,\epsilon )\cap \pi (W)$ is an open definable neighborhood of $F$ in $Z,$ $F$ is closed in $Z$ and $Z$ is definably normal, there is an open definable neighborhood $O$ of $F$ in $Z$ such that $O\subseteq \bar{O}\cap Z\subseteq U(F,\epsilon )\cap \pi (W)\subset \pi (W)$ and  $S\cap \pi ^{-1}(\bar{O}\cap Z)\subseteq W.$
\qed \\

\begin{lem}\label{lem CxN special}
Let $Z\subseteq M^n$ be a definably normal definable subset. Let $S$ and $T$ be disjoint, closed definable subsets of $Z\times [a,b]$ such that $\pi (S)=\pi (T).$ Then there are disjoint open definable subsets $U$ and  $V$ of $Z\times [a,b]$ such that $S\subseteq V$ and  $T\subseteq U.$  
\end{lem}

\pf
Since $S\cap T=\emptyset $ and by o-minimality (\cite[Chapter 3, (3.6)]{vdd}) there is  $N$ such that for each $z\in Z$ the definable sets  $S\cap \pi ^{-1}(z)$ and $T\cap \pi ^{-1}(z)$ are each a union of at most $N$ points and intervals (in $\pi ^{-1}(z)$),  by definable Skolem functions, there are  definable maps $f^-_l,f^+_l, g^-_l,g^+_l : Z\to M$ for  $l=1,\ldots , N,$ such that 
$$f^-_1<f^+_1<g^-_1<g^+_1<f^-_2<f^+_2<\ldots $$ 
and 
\begin{enumerate}
\item[-]
 $\bigcup _l(f^-_l,f^+_l)_Z$ contains $S;$ 
 \item[-]
 $\bigcup _l(g^-_l,g^+_l)_Z$ contains $T;$
\item[-]
 $\bigcup _l(f^-_l,f^+_l)_Z\cap T=\emptyset;$ 
 \item[-]
 $\bigcup _l(g^-_l,g^+_l)_Z\cap S=\emptyset .$
\end{enumerate}

Note that by construction, for each $z\in \pi (S)$  the definable set $(f^-_l,f^+_l)_Z\cap (S\cap \pi ^{-1}(z))$ is either empty, a point or a closed interval (in $\pi ^{-1}(z)$). Similarly,   for each $z\in \pi (T)$  the definable set $(g^-_l,g^+_l)_Z\cap (T\cap \pi ^{-1}(z))$ is either empty, a point or a closed interval (in $\pi ^{-1}(z)$).

Let $D=\pi (S)=\pi (T).$ We prove the result by induction on $\dim D.$ Suppose that $\dim D=0,$  so $D=\{d_0,\ldots , d_k\}\subseteq Z.$ For $i=0,\ldots , k$ and $l=1,2,\ldots ,N,$ let $s^-_{l,i}=f^-_l(d_i), s^+_{l,i}=f^+_l(d_i), t^-_{l,i}=g^-_l(d_i), t^+_{l,i}=g^+_l(d_i).$ Then 
$$s^-_{1,i}<s^+_{1,i}<t^-_{1,i}<t^+_{1,i}<s^-_{2,i}<s^+_{2,i}<\ldots $$ and 
\begin{enumerate}
\item[-]
$\bigcup _l(s^-_{l,i}, s^+_{l,i})_{\{d_i\}}$ contains $S\cap \pi ^{-1}(d_i);$ 
\item[-]
$\bigcup _l(t^-_{l,i},t^+_{l,i})_{\{d_i\}}$ contains $T\cap \pi ^{-1}(d_i);$
\item[-]
$\bigcup _l(s^-_{l,i}, s^+_{l,i})_{\{d_i\}}\cap T=\emptyset ;$ 
\item[-]
$\bigcup _l(t^-_{l,i},t^+_{l,i})_{\{d_i\}}\cap S=\emptyset .$
\end{enumerate}

Because $Z$ is definably normal, for each $i=0,\ldots , k,$ let  $B_{i}$  be  an open definable neighborhood of $d_i$ in $Z,$ such that  $B_i\cap B_{i'}=\emptyset $ for $i\neq i'.$ Let $\sigma ^-_{l,i}, \sigma ^+_{l,i}, \tau ^-_{l,i}, \tau ^+_{l,i}:B_{i}\to M$ be the constant definable maps with values $s^-_{l,i}, s^+_{l,i}, t^-_{l,i}$ and $t^+_{l,i}$ respectively.  Then 
$$\sigma^-_{1,i}<\sigma^+_{1,i}<\tau^-_{1,i}<\tau ^+_{1,i}<\sigma^-_{2,i}<\sigma^+_{2,i}<\ldots .$$
Moreover, if     
$$V_{i}=\bigcup _l(\sigma ^-_{l,i},\sigma ^+_{l,i})_{B_i}\cap Z\times [a,b]$$ 
and   
$$U_{i}=\bigcup _l (\tau ^-_{l,i}, \tau ^+_{l,i})_{B_i}\cap Z\times [a,b],$$
then $V_i$ and $U_i$ are  open definable subsets of $Z\times [a,b]$ such that
$$S\cap \pi ^{-1}(\{d_i\})\subseteq V_i,$$
$$T\cap \pi ^{-1}(\{d_i\})\subseteq U_i$$ 
and  $U_i\cap V_{i'}=\emptyset $ for every $i$ and $i'.$ 
Therefore, if  $V=\bigcup _iV_i$ and $U=\bigcup _iU_{i},$ then $V$ and $U$ are  open definable subsets of $Z\times [a,b]$ such that $S\subseteq V,$  $T\subseteq U$ and $V\cap U=\emptyset .$

Suppose that $\dim D=k$ and the result holds for every pair $(S',T')$ of disjoint, closed definable subsets of $Z\times [a,b]$ such that  $\pi (S')=\pi (T')$ and this pojection has dimension smaller than $k.$  

Let $C$ be the subset of $D$ on which the definable maps $f^-_l,f^+_l,$ $ g^-_l,g^+_l$ for $l=1,2, \ldots $ are all continuous and let $E=D\setminus C.$ Then $C$ and $E$ are definable, $C=C_0\sqcup \ldots \sqcup C_k$ with each $C_i$ a cell with $\dim C_i=\dim D$ and $\dim E<\dim D$ (\cite[Chapter 3, (2.11)]{vdd}). Since,  $\dim \bar{E}\cap Z=\dim E$ (\cite[Chapter 4, (1.8)]{vdd}),  by the induction hypothesis,  there are $V''$  and $U''$   open definable subsets of $Z\times [a,b]$ such that
$$S\cap \pi ^{-1}(\bar{E}\cap Z)\subseteq V'',$$
$$T\cap \pi ^{-1}(\bar{E}\cap Z)\subseteq U''$$
and  $V''\cap U''=\emptyset .$

Now  the definable maps $f^-_{l|C_i}, f^+_{l|C_i}, g^-_{l|C_i}, $ $g^+_{l|C_i}$ are all continuous and 
$$f^-_{1|C_i}< f^+_{1|C_i}< g^-_{l|C_i}< g^+_{l|C_i}<f^-_{2|C_i}< f^+_{2|C_i}<\ldots .$$ By the last paragraph of Remark \ref{nrmk cells perm} we may extend them to continuous definable maps  $\sigma ^-_{l,i}, \sigma ^+_{l,i}, \tau ^-_{l,i}, \tau ^+_{l,i}: B_i \to M$ for $l=1,2,\ldots $  such that 
$$\sigma ^-_{1,i}< \sigma ^+_{1,i}< \tau ^-_{1,i}< \tau ^+_{1,i}< \sigma ^-_{2,i}< \sigma ^+_{2,i}<\ldots ,$$  
where $B_i$ is an open definable neighborhood of $C_i$ in $Z.$ 
Moreover, if 
$$V'_{i}=\bigcup _l(\sigma ^-_{l,i},\sigma ^+_{l,i})_{B_i}\cap Z\times [a,b]$$ 
 and 
$$U'_{i}=\bigcup _l (\tau ^-_{l,i}, \tau ^+_{l,i})_{B_i}\cap Z\times [a,b],$$
then  $V'_i$ and $U'_i$ are  open definable subsets of $Z\times [a,b]$ such that 
$$S\cap \pi ^{-1}(C_i)\subseteq V'_i,$$
$$T\cap \pi ^{-1}(C_i)\subseteq U'_i$$ 
and  $U'_i\cap V'_i=\emptyset $ for each $i.$ However we might have $U'_{i}\cap V''\neq \emptyset $ or $V'_i\cap U''\neq \emptyset $ or $U'_i\cap V'_{i'}\neq \emptyset$ for $i\neq i'.$  We will first modify $U'_i, V'_i, U''$ and $V''$ so that $U'_{i}\cap V''=\emptyset $ and $V'_i\cap U''= \emptyset $ 

Since $\pi (V'')\cap \pi (U'')$ is an open definable neighborhood of $\bar{E}\cap Z$ in $Z$ and $Z$ is definably normal, by Lemma \ref{lem closed fiber in open},  there is an open definable neighborhood $O''$ of $\bar{E}\cap Z$ in $Z$ such that $O''\subseteq \bar{O''}\cap Z\subseteq \pi (V'')\cap \pi (U''),$  
$$S\cap \pi ^{-1}(\bar{E}\cap Z)\subseteq S\cap \pi ^{-1}(\bar{O''}\cap Z)\subseteq  V''$$
 and  
 $$T\cap \pi ^{-1}(\bar{E}\cap Z)\subseteq T\cap \pi ^{-1}(\bar{O''}\cap Z)\subseteq U''.$$ 
 In particular, we also have $S\cap \pi ^{-1}(\bar{O''}\cap Z)\cap \bar{U''}=\emptyset $ and   $T\cap \pi ^{-1}(\bar{O''}\cap Z)\cap \bar{V''}=\emptyset .$

Let 
$$U'=\pi ^{-1}(O'')\cap U''$$ 
and 
$$V'=\pi ^{-1}(O'')\cap V''.$$ 
Then $U'$ and $V'$ are  open definable subsets of $Z\times [a,b]$  such that $V'\cap U'=\emptyset ,$  
$$S\cap \pi ^{-1}(\bar{E}\cap Z)\subseteq S\cap \pi ^{-1}(O'')\subseteq V'$$ 
and 
$$T\cap \pi ^{-1}(\bar{E}\cap Z)\subseteq T\cap \pi ^{-1}(O'')\subseteq  U'.$$

For each $i=0,\ldots , k,$ let 
$$V_i=V'_i\setminus (\pi ^{-1}(\bar{O''}\cap Z)\cap \bar{U''})$$ 
and 
$$U_i=U'_i\setminus (\pi ^{-1}(\bar{O''}\cap Z)\cap \bar{V''}).$$
 Then $V_i$  and $U_i$ are  open definable subsets of $Z\times [a,b]$ such that 
$$S\cap \pi ^{-1}(C_i)\subseteq V_i$$ and 
$$T\cap \pi ^{-1}(C_i)\subseteq U_i.$$ 
Moreover,  $U_i\cap V_i=\emptyset ,$ $U_i\cap V'=\emptyset $  and $V_i\cap U'=\emptyset $ for each $i.$ We would also like to have $U_i\cap V_{i'}=\emptyset $ for all $i,i'.$ For that we will need to shrink the $U_i$'s and $V_i$'s without destroying what we already achieved.

For each $i=0,\ldots ,k,$ let $D_i=C_i\setminus O''.$ Then $D_i$ is a closed definable subset of $Z.$ In fact, let $z\in Z$ such that every open definable neighborhood of $z$ in $Z$ intersects $D_i.$ Then $z\in D$ and so either $z\in E$ or $z\in C_{i'}$ for some $i'.$ It cannot be the case that $z\in E$ for then $O''$ is an open definable neighborhood of $z$ in $Z$ which does not intersect $D_i.$ Now since $C_{i'}$ is relatively open in $D$ we have $i'=i$ and $z\in D_i.$ 

Since the $D_i$'s are closed definable subsets of $Z,$  two by two disjoint and $B_i\setminus (\bar{E}\cap Z)$ is an open definable neighborhood of $D_i$ in $Z$, because $Z$ is definably normal, for each $i$ there is an open definable neighborhood $O_i$ of $D_i$ in $Z$ such that  $O_i\subseteq \bar{O_i}\cap Z\subseteq B_i\setminus (\bar{E}\cap Z)$ and  $O_i\cap O_{i'}=\emptyset $ for $i\neq i'.$

For each $i$ let 
 
$${\mathcal V}_{i}=V_i\cap \pi ^{-1}(O_i)$$ 
and 
$${\mathcal U}_i=U_i\cap \pi ^{-1}(O_i).$$
Then ${\mathcal V}_i$ and ${\mathcal U}_i$ are  open definable subsets of $Z\times [a,b]$ such that 
$$S\cap \pi ^{-1}(D_i)\subseteq {\mathcal V}_i$$
and
$$T\cap \pi ^{-1}(D_i)\subseteq {\mathcal U}_i$$ 
Moreover,   ${\mathcal U}_i\cap V'=\emptyset $  and ${\mathcal V}_i\cap U'=\emptyset $ and   ${\mathcal U}_i\cap {\mathcal V}_{i'}=\emptyset $ for all $i, i'.$ 

Finally, let 
$$V=(\bigcup _{i=1}^k{\mathcal V}_i)\cup V'$$
and 
$$U=(\bigcup _{i=1}^k{\mathcal U}_i)\cup U'.$$
Then  $V$ and $U$ are  open definable subsets of $Z\times [a,b]$ and since $C_i=D_i\sqcup (C_i\cap O''),$ we have that $S\subseteq V,$  $T\subseteq U$ and $U\cap V=\emptyset .$ 
\qed \\

\noindent
{\bf Proof of Proposition \ref{prop CxN affine is normal}:} We have $K\subseteq [a, b]^k.$  So if $Z \times [a,b]^k$ is definably normal, then $Z \times K$ is definably normal as well being a closed definable subset. By induction we can easily reduce to $k=1,$ so we just need to show that $Z\times [a,b]$ is definably normal, i.e., given disjoint, closed definable subsets $S$ and $T$ in $Z\times [a,b]$, we need to find disjoint definable open neighborhoods $V \supseteq S$ and $U\supseteq T$  in $Z\times [a,b].$

If  $\pi (S)\cap \pi (T)=\emptyset ,$ since $Z$ is definably normal and $\pi :Z\times [a,b]\to Z$ is closed in $\Df,$ there are disjoint, definable open neighborhoods $V'\supseteq \pi (S)$ and $U'\supseteq \pi (T)$ in  $Z.$ But then $V=\pi ^{-1}(V')\supseteq S$ and $U=\pi ^{-1}(U')\supseteq T$ are disjoint, definable open neighborhoods in  $Z\times [a,b].$ 

So suppose that $D=\pi (S)\cap \pi (T)\neq \emptyset.$ Then by Lemma \ref{lem CxN special} there are disjoint definable open neighborhoods $V''\supseteq S\cap \pi ^{-1}(D)$ and $U''\supseteq T\cap \pi ^{-1}(D)$ in $Z\times [a,b].$ Since $\pi (V'')\cap \pi (U'')$ is an open definable neighborhood of $D$ in $Z$ and $Z$ is definably normal, by Lemma \ref{lem closed fiber in open},  there is an open definable neighborhood $O''$ of $D$ in $Z$ such that $O''\subseteq \bar{O''}\cap Z\subseteq \pi (V'')\cap \pi (U''),$  
$$S\cap \pi ^{-1}(D)\subseteq S\cap \pi ^{-1}(\bar{O''}\cap Z)\subseteq  V''$$
 and  
 $$T\cap \pi ^{-1}(D)\subseteq T\cap \pi ^{-1}(\bar{O''}\cap Z)\subseteq U''.$$ 
 In particular, we also have $S\cap \pi ^{-1}(\bar{O''}\cap Z)\cap \bar{U''}=\emptyset $ and   $T\cap \pi ^{-1}(\bar{O''}\cap Z)\cap \bar{V''}=\emptyset .$

Let 
$$U'=\pi ^{-1}(O'')\cap U''$$ 
and 
$$V'=\pi ^{-1}(O'')\cap V''.$$ 
Then $U'$ and $V'$ are  open definable subsets of $Z\times [a,b]$  such that $V'\cap U'=\emptyset ,$  
$$S\cap \pi ^{-1}(D)\subseteq S\cap \pi ^{-1}(O'')\subseteq V'$$ 
and 
$$T\cap \pi ^{-1}(D)\subseteq T\cap \pi ^{-1}(O'')\subseteq  U'.$$

The definable sets  $\pi (S)\setminus O''$ and $\pi (T)\setminus O''$  are disjoint closed definable subsets of $Z.$ Therefore, since  $Z\setminus D$ is an open definable neighborhood of them  in $Z$ and  $Z$ is definably normal,  there are disjoint open definable neighborhoods $W_S\supseteq \pi (S)\setminus O''$  and $W_T\supseteq \pi (T)\setminus O''$  in $Z$ such that  $W_S\subseteq \bar{W_S}\cap Z\subseteq Z\setminus D$ and  $W_T\subseteq \bar{W_T}\cap Z\subseteq Z\setminus D.$

Let 
 
$${\mathcal V}=\pi ^{-1}(W_S)\setminus (\pi ^{-1}(\bar{O''}\cap Z)\cap \bar{U''})$$ 
and 
$${\mathcal U}=\pi ^{-1}(W_T)\setminus (\pi ^{-1}(\bar{O''}\cap Z)\cap \bar{V''}).$$
Then ${\mathcal V}$ and ${\mathcal U}$ are  open definable subsets of $Z\times [a,b]$ such that 
$$S\cap \pi ^{-1}(\pi (S)\setminus O'')\subseteq {\mathcal V}$$
and
$$T\cap \pi ^{-1}(\pi (T)\setminus O'')\subseteq {\mathcal U}$$ 
Moreover,   ${\mathcal U}\cap V'=\emptyset ,$   ${\mathcal V}\cap U'=\emptyset $ and   ${\mathcal U}\cap {\mathcal V}=\emptyset .$ 

Finally, let 
$$V={\mathcal V}\cup V'$$
and 
$$U={\mathcal U}\cup U'.$$
Then  $V$ and $U$ are  open definable subsets of $Z\times [a,b]$ and since $\pi (S)=(\pi (S)\setminus O'')\sqcup (\pi (S)\cap O'')$ and $\pi (T)=(\pi (T)\setminus O'')\sqcup (\pi (T)\cap O'')$ we have that $S\subseteq V,$  $T\subseteq U$ and $U\cap V=\emptyset .$ 
\qed \\

\end{subsection}

\begin{subsection}{Fibers in $\tDf$}\label{subsection  spec fibers}
A fiber $f^{-1}(\alpha )$ of a morphism $f:X\to Y$ in $\tDf$ is not in general an object of $\tDf$, but by a model theoretic trick such fibers can each  be seen as objects of $\tDfs$ for some elementary extension ${\mathbb S}$   of ${\mathbb M}$. This follows from our assumption that ${\mathbb M}$ has definable Skolem functions. \\


Before we proceed we recall the basic model theoretic consequences of our assumption on ${\mathbb M}$ that will be required later.\\

By \cite[Theorem 5.1]{pis}:

\begin{fact}\label{fact def skolem4}
Since ${\mathbb M}$ has definable Skolem functions, then for any  parameter $v$ in some model of the  first-order  theory of ${\mathbb M}$, there is a prime model ${\mathbb S}$ of the  first-order  theory of ${\mathbb M}$ over $M\cup \{v\}$ such that, for every $s\in S$, there is a definable set $D$ and a definable map $f:D\to M$ (all definable in ${\mathbb M}$), such that $v\in D({\mathbb S})$ and $s=f^{{\mathbb S}}(v).$
\end{fact}

By Fact \ref{fact def skolem4} we have:

\begin{fact}\label{fact prime exchange}
Since ${\mathbb M}$ has definable Skolem functions, if  ${\mathbb S}$ is a prime model of the  first-order  theory of ${\mathbb M}$ over $M\cup \{v\}$ where $v$ is a  parameter  in some model of the  first-order theory of ${\mathbb M},$ then if $u\in S$ and ${\mathbb K}$ is a prime model of the theory of ${\mathbb M}$ over $M\cup \{u\}$, then  we have either ${\mathbb K}={\mathbb M}$ or ${\mathbb K}={\mathbb S}.$
\end{fact}

Indeed, consider the model theoretic algebraic closure operator ${\rm acl}_{{\mathbb M}}(\bullet )$ in ${\mathbb S}$ given by $a\in {\rm acl}_{{\mathbb M}}(C)$ with $C\subseteq S$ if and only if there is $c\in C$ and a definable set $D$ and a definable map $f:D \to M$ (all definable in ${\mathbb M}$) such that $c\in D({\mathbb S})$ and $f^{{\mathbb S}}(c)=a$.
By \cite[Theorem 4.1]{pis}, this model theoretic algebraic closure operator satisfies the exchange property: if $a\in {\rm acl}_{{\mathbb M}}(Cb)\setminus {\rm acl}_{{\mathbb M}}(C)$, then $b\in {\rm acl}_{{\mathbb M}}(Ca)$.
Thus, by Fact \ref{fact def skolem4}, if $u\in S$ and ${\mathbb K}$ is a prime model of the theory of ${\mathbb M}$ over $M\cup \{u\}$, then  we have either ${\mathbb K}={\mathbb M}$ or ${\mathbb K}={\mathbb S}.$\\

Let ${\mathbb S}$ be an elementary extension   of ${\mathbb M}$. For each object $\tilde{X}$ of $\tDf$ we have a restriction map $r: \tilde{X({\mathbb S})} \to \tilde{X}$ such that given an ultrafilter $\alpha \in \tilde{X({\mathbb S})}$, $r(\alpha )$ is the  ultrafilter in $\tilde{X}$ determined by the collection $\{A: A({\mathbb S}) \in \alpha \}$. This is a continuous surjective map which is neither open nor closed.\\

The following result was proved in \cite[Lemma 3.1]{bera}   in the affine case - for continuous definable maps between definable sets,  under the assumption that ${\mathbb M}$ is an o-minimal expansion of an ordered group. However, the only thing needed from ${\mathbb M}$ is that it has definable Skolem functions. In fact all that is required is Fact \ref{fact def skolem4}. A special case of this result, when the map is a projection $X\times [a,b]\to X$, was proved before  in \cite[Claim 4.5]{ejp}. \\

\begin{lem}\label{lem bera}
Let $\tilde{f}:\tilde{X}\to \tilde{Y}$ be a morphism in $\tDf$, $\alpha \in \tilde{Y}$, $a\models \alpha $ a realization of $\alpha $ and ${\mathbb S}$ a prime model of  the first-order  theory of ${\mathbb  M}$ over $\{a\}\cup M.$ Then there is a homeomorphism
$$r_{|}: \tilde{(f^{\mathbb S})^{-1}(a)}\to \tilde{f}^{-1}(\alpha )$$
induced by the restriction $r:\tilde{X({\mathbb S})}\to \tilde{X}.$
\end{lem}

\pf
Since $(f^{-1}(B))({\mathbb S})=(f^{{\mathbb S}})^{-1}(B({\mathbb S}))$ for every definable subset $B\subseteq Y$, we have a commutative diagram
$$
\xymatrix{
\tilde{(f^{{\mathbb S}})^{-1}(a)} \ar@{^{(}->}[r] \ar[d]^{r_|} & \tilde{X({\mathbb S})} \ar[d]^{r}\ar[r]^{\tilde{f^{{\mathbb S}}}}  & \tilde{Y({\mathbb S})} \ar[d]^{r} \\
\tilde{f}^{-1}(\alpha ) \ar@{^{(}->}[r] &\tilde{X} \ar[r]^{\tilde{f}} &\tilde{Y}
}
$$
and so $r_{|}: \tilde{(f^{\mathbb S})^{-1}(a)}\to \tilde{f}^{-1}(\alpha )$ is well defined.

We now have to show that $r_{|}: \tilde{(f^{\mathbb S})^{-1}(a)}\to \tilde{f}^{-1}(\alpha )$ is a continuous and open bijection.  Let $Y_i$ be a definable chart of $Y$ such that $\alpha \in \tilde{Y_i}$ and let $Z_i=f^{-1}(Y_i).$  Since $\tilde{((f_{|Z_i})^{{\mathbb S}})^{-1}(a)}=  \tilde{(f^{{\mathbb S}})^{-1}(a)}$ and  $\tilde{(f_{|Z_i})}^{-1}(\alpha )= \tilde{f}^{-1}(\alpha )$, the map we are interested on is the same as the restriction in the following commutative diagram
$$
\xymatrix{
\tilde{((f_{|Z_i})^{{\mathbb S}})^{-1}(a)} \ar@{^{(}->}[r] \ar[d]^{r_|} & \tilde{Z_i({\mathbb S})} \ar[d]^{r}\ar[r]^{\tilde{f^{{\mathbb S}}}}  & \tilde{Y_i({\mathbb S})} \ar[d]^{r} \\
\tilde{f_{|Z_i}}^{-1}(\alpha ) \ar@{^{(}->}[r] &\tilde{Z_i} \ar[r]^{\tilde{f}} &\tilde{Y_i}.
}
$$
Let $X_j$ be a definable chart of $X.$   Since $\tilde{((f_{|Z_i\cap X_j})^{{\mathbb S}})^{-1}(a)}=  \tilde{((f_{|Z_i})^{{\mathbb S}})^{-1}(a)}\cap \tilde{X_j({\mathbb S})}$ and   $\tilde{(f_{|Z_i\cap X_j})}^{-1}(\alpha )= \tilde{f_{|Z_i}}^{-1}(\alpha )\cap \tilde{X_j}$, then  the restriction in the following commutative diagram
$$
\xymatrix{
\tilde{((f_{|Z_i\cap X_j})^{{\mathbb S}})^{-1}(a)} \ar@{^{(}->}[r] \ar[d]^{r_|} & \tilde{(Z_i\cap X_j)({\mathbb S})} \ar[d]^{r}\ar[r]^{\tilde{f^{{\mathbb S}}}}  & \tilde{Y_i({\mathbb S})} \ar[d]^{r} \\
\tilde{f_{|Z_i\cap X_j}}^{-1}(\alpha ) \ar@{^{(}->}[r] &\tilde{Z_i\cap X_j} \ar[r]^{\tilde{f}} &\tilde{Y_i}
}
$$
is a continuous and open bijection by  \cite[Lemma 3.1]{bera}. Therefore, since the $\tilde{X_j}$'s (resp. $\tilde{X_j({\mathbb S})}$) cover $\tilde{f_{|Z_i}}^{-1}(\alpha )$ (resp. $\tilde{((f_{|Z_i})^{{\mathbb S}})^{-1}(a)}$), it follows that $r_{|}: \tilde{(f^{\mathbb S})^{-1}(a)}\to \tilde{f}^{-1}(\alpha )$ is a continuous and open surjection.

Let $\beta _1, \beta _2\in \tilde{(f^{{\mathbb S}})^{-1}(a)}=\tilde{((f_{|Z_i})^{{\mathbb S}})^{-1}(a)}$ be such that $r_{|}(\beta _1)=r_{|}(\beta _2)\in \tilde{f}^{-1}(\alpha )=\tilde{(f_{|Z_i})}^{-1}(\alpha )$. Let $X_j$ be a  definable chart   of $X$ such that $r_|(\beta _1)=r_|( \beta _2)\in \tilde{X_j}.$ Then
$r_|(\beta _1)=r_|( \beta _2)\in \tilde{f_{|Z_i}}^{-1}(\alpha )\cap \tilde{X_j} = \tilde{(f_{|Z_i\cap X_j})}^{-1}(\alpha ).$ On the other hand, we also have $\beta _1, \beta _2\in \tilde{X_j({\mathbb S})}$ and so $\beta _1, \beta _2\in \tilde{((f_{|Z_i})^{{\mathbb S}})^{-1}(a)}\cap \tilde{X_j({\mathbb S})} = \tilde{((f_{|Z_i\cap X_j})^{{\mathbb S}})^{-1}(a)}.$  Hence,   by the above, $\beta _1=\beta _2.$ This proves that $r_{|}: \tilde{(f^{\mathbb S})^{-1}(a)}\to \tilde{f}^{-1}(\alpha )$ is also an injection as required.
\qed \\


\end{subsection}

\begin{subsection}{On  families of supports on fibers in $\tDf$}\label{section supp on spec fibers}
In the paper \cite{ep1} normal and constructible families of supports on objects of $\tDf$ played a fundamental role.  Here we introduce the notion of a normal and constructible family of supports $\Phi $ on a fiber $f^{-1}(\alpha )$ of a morphism $f:X\to Y$ in $\tDf.$ \\

First recall the following definitions from \cite[page 1267]{ep1}.

\begin{defn}\label{def const supp}
{\em
Let $X$ be an object in $\tDf$. A family  $\Phi $ of closed subsets of  $X$ is a {\em family of supports on $X$} if:
\begin{itemize}
\item
every closed subset of a member of $\Phi $ is in $\Phi$;
\item
$\Phi $ is closed under finite unions.
\end{itemize}
A family $\Phi $ of supports on $X$ is said to be {\em constructible} if:
\begin{itemize}
\item
every member of $\Phi $ is contained in a  member of $\Phi $ which is constructible.
\end{itemize}
A family $\Phi $ of supports on $X$  is said to be {\em normal} if:
\begin{itemize}
\item
every member of $\Phi $ is normal;

\item
for each member $S$ of $\Phi $, if $U$ is an open neighborhood of $S$ in $X$, then there exists a closed constructible neighborhood of $S$ in $U$ which is a member of $\Phi.$\\
\end{itemize}
}
\end{defn}

\begin{nrmk}\label{nrmk cons sup}
{\em
Let $X$ be an object in $\tDf$ and let $\Phi $ be a constructible family of supports on $X.$ If $C\in \Phi $ and $V$ is an open neighborhood of $C$ in $X,$ then there is a constructible closed subset $B$ of $X$ and a constructible open subset $U$ of $X$ such that $C\subseteq B\subseteq U\subseteq V$ and $B\in \Phi.$ When $\Phi $ is also normal, we may take $B$ to closed and constructible neighborhood of $C$ in $V.$

In fact, since $C$ is quasi-compact we find  $U\subseteq V.$ Since, by Lemma \ref{lem spec f qc}, $C=\bigcap _{i\in I}C_i\subseteq U$  with $C_i$'s closed constructible subsets of $X,$ we have $X\setminus U\subseteq \bigcup _{i\in I}X\setminus C_i.$ Since $X\setminus U$ is quasi-compact, being closed, there is $I_0\subseteq I$ finite such that $\bigcap _{i\in I_0}C_i\subseteq U.$ On the other hand, since $\Phi$ is constructible, there is $D\in \Phi$ constructible such that $C\subseteq D.$ Take $B=D\cap \bigcap _{i\in I_0}C_i.$
}
\end{nrmk}

Working with the previous definition in $\tDfs$ we can  use the homeomorphism $r_{|}: (f^{\mathbb S})^{-1}(a)\to f^{-1}(\alpha )$ of Lemma \ref{lem bera} to define the notion of a normal and constructible  family of supports on the fibers $f^{-1}(\alpha )$:\\

\begin{defn}\label{defn c supp on fiber}
{\em
Let $f:X\to Y$ be a morphism in $\tDf$, $\alpha \in Y$, $a\models \alpha $ a realization of $\alpha $ and ${\mathbb S}$ a prime model of  the first-order  theory of ${\mathbb  M}$ over $\{a\}\cup M.$
\begin{itemize}
\item
A family of supports $\Psi $ on the quasi-compact space $f^{-1}(\alpha )$ is {\em constructible} if its inverse image $(r_{|})^{-1}\Psi $  is a constructible family of supports  on the object $(f^{\mathbb S})^{-1}(a)$ of $\tDfs$.

\item
A family of supports $\Psi $ on the quasi-compact space $f^{-1}(\alpha )$ is {\em normal} if its inverse image $(r_{|})^{-1}\Psi $  is a normal family of supports  on the object $(f^{\mathbb S})^{-1}(a)$ of $\tDfs$.

\end{itemize}
If there is no risk of confusion, and since  $r_{|}$ is a homeomorphism, we also use $\Psi $ to denote the inverse image of $\Psi $ by $r_{|}$.
}
\end{defn}

We also say that a subset $Z$ of $f^{-1}(\alpha )$ is {\em constructible}  if its inverse image $(r_{|})^{-1}(Z) $  is a constructible  subset   of the object $(f^{\mathbb S})^{-1}(a)$ of $\tDfs$. Of course, a subset $Z$ of $f^{-1}(\alpha )$ is quasi-compact  if its inverse image $(r_{|})^{-1}(Z) $  is a quasi-compact subset   of the object $(f^{\mathbb S})^{-1}(a)$ of $\tDfs$.\\

Let $X$ be an object  of $\tDf.$ Then:
\begin{itemize}
\item[$\bullet $]
The   {\em family of complete supports on  $X$}, denoted $c$, is the constructible family of  supports on  $X$ of all closed subsets $A$ of $X$ with $A\subseteq Z$ for some  constructible complete in $\tDf$ subset $Z$ of $X.$
\item[$\bullet $]
The   {\em family of complete supports on  $X({\mathbb S})$}, also denoted $c$, is the constructible family of  supports on  $X({\mathbb S})$ of all closed subsets $A$ of $X({\mathbb S})$ with $A\subseteq Z$ for some  constructible complete in $\tDfs$ subset $Z$ of  $X({\mathbb S}).$\\
\end{itemize}

We can now introduce one of the main definitions of the paper:

\begin{defn}\label{defn compact supp}
{\em
Let $f:X\to Y$ be a morphism in $\tDf$, $\alpha \in Y$, $a\models \alpha $ a realization of $\alpha $ and ${\mathbb S}$ a prime model of the first-order theory of ${\mathbb M}$ over $\{a\}\cup M.$
\begin{itemize}
\item[$\bullet $]
The  {\em family of complete supports on  $(f^{\mathbb S})^{-1}(a)$} (an object of $\tDfs$), denoted $c$, is the constructible family of  supports on  $(f^{\mathbb S})^{-1}(a)$ of all closed subsets $A$ of $(f^{\mathbb S})^{-1}(a)$ with $A\subseteq Z$ for some  constructible complete in $\tDfs$ subset $Z$ of $(f^{\mathbb S})^{-1}(a).$

\item[$\bullet $]
The  {\em family of complete supports on $f^{-1}(\alpha )$}, denoted $c$, is the (constructible) family of supports on  $f^{-1}(\alpha )$ whose  inverse image by $r_{|}$ is the constructible family of complete supports on $(f^{\mathbb S})^{-1}(a).$ \\
\end{itemize}
}
\end{defn}

\end{subsection}

\begin{subsection}{The full subcategories $\tbA$ of $\tDf$}\label{subsection tbA} Here we introduce the full subcategories $\tbA$ of $\tDf$ on which we develop the Grothendieck formalism of the six operations on o-minimal sheaves. We  explain exactly which properties of objects and morphisms of $\tbA$  are used later and give the main examples of such subcategories.\\

The full subcategories $\tbA$ of $\tDf$ on which we must work  are
such that: \label{cat tA}
\begin{itemize}
\item[(A0)]
cartesian products (in $\tDf$) of objects of $\tbA$ are objects of $\tbA$ and locally closed constructible subsets of objects of $\tbA$ are objects of $\tbA;$
\item[(A1)]
in every object of $\tbA$  every open constructible
subset is a finite union of open and
normal constructible
subsets;
\item[(A2)]
every object of $\tbA$ has a 
completion in $\tbA.$
\end{itemize}
For some of the results about the proper direct image we will required that the morphisms $\tilde{f}:\tilde{X}\to \tilde{Y}$ in $\tbA$ involved satisfy the following: 

\begin{itemize}
\item[(A3)]
if $u\in Y$,  then for every elementary extension  ${\mathbb S}$ of ${\mathbb M}$ and every $F\in \mod (A_{X_{\df}})$ we have an isomorphism
$$H^*_c(f^{-1}(u); F_{|f^{-1}(u)})\simeq H^*_c((f^{{\mathbb S}})^{-1}(u); F({\mathbb S})_{|(f^{{\mathbb S}})^{-1}(u)})$$
where $\tilde{F}({\mathbb S})=r^{-1}\tilde{F},$ $r:\tilde{X({\mathbb S})}\to \tilde{X}$ is the restriction and $H_c^*$ is the cohomology with definably compact supports (\cite[Example 2.10 and Definition 2.12]{ep1}).\\
\end{itemize}

Examples of categories $\tbA$ and morphisms in such categories satisfying (A3) 
will be given below. Before let us compare our conditions (A0) to (A3) with those used by Delfs in the semi-algebraic case (\cite{D3}):\\

\begin{nrmk}\label{nrml a1}
{\em
In the semi-algebraic setting instead of (A1) the following stronger version is used, e.g. in the proof of \cite[Chapter II, Section 8, Proposition 8.2 and Corollaries 8.3 and 8.4]{D3}:\\
\begin{itemize}
\item[(A1)$^*$]
in every object of $\tbA$  every open constructible
subset is normal.\\
\end{itemize}
The reason is that, in that setting (\cite[Chapter II, Section 8]{D3}), one works with the tilde of  locally complete semi-algebraic spaces which are in particular regular, and regular semi-algebraic spaces are affine and hence semi-algebraically normal. If we were only interested in working with the tilde of regular definable spaces in o-minimal expansions of real closed fields, then for the same reason (\cite[Chapter 10, (1.8) and Chapter 6 (3.5)]{vdd}) the stronger assumption (A1)$^*$ would also hold.

In our application of the results of this paper to the solution of Pillay's conjecture for definably compact groups in arbitrary o-minimal structures (\cite{oh5p}) we have to work with the tilde of locally closed definable subsets of definable manifolds defined in cartesian products of definable group-intervals. In that setting (A1) always holds (\cite[Corollary 3.15]{oh5p}) but if these definable group-intervals are orthogonal (in a model theoretic sense) to each other then (A1)$^*$ fails (\cite[Example 3.16]{oh5p}).\\
}
\end{nrmk}

The only  consequence of  (A1) which will be  used  throughout the paper 
 is the following result. Thus this result could be a replacement for (A1).

 Before we prove the result let us define an object $X$ of $\tDf$ to be  {\it affine} if it is the tilde of an affine object of $\Df.$ Recall that an affine definable space is a definable space which is definably homeomorphic to a definable set equipped with its induced topology. \\

\begin{prop}\label{prop obj tA}
Let $X$ be an object of $\tbA.$ 
If $\alpha \in X$, then there is an open, affine, normal  constructible subset $U$ of $X$ such that $\alpha \in U$ and $\alpha $ is closed in $U$.
\end{prop}



\pf
The specializations of a point $\alpha $ in an object $Y$ of $\tDf$ form  finite chains  by \cite[Lemma 2.11]{ejp} and if the object is  normal, $\alpha $ has a  unique closed specialization $\rho $ (Fact \ref{fact ominspec closed point}) which must be the minimum of any chain of specialization of $\alpha .$

Since $X$ has a finite cover by affine, open, constructible subsets, we may assume that $X$ is affine.  By (A1) let $U$ be an open, normal constructible open subset of $X$ such that $\alpha \in U$. Let $\rho $ be the unique closed specialization of $\alpha $ in $U$. If $\alpha $ is closed in $U$ we are done, otherwise,  $U\setminus \{\rho \}$ is open (in $U$ and so also in $X$), and so, by (A1), there is an open, normal constructible open subset $V$ of $U\setminus \{\rho \}$ (and so of $X$) such that $\alpha \in V.$ Since every specialization of $\alpha $ in $V$ is  also a specialization of $\alpha $ in $U$ and the maximum length of chains of specializations of $\alpha $ in $V$ is smaller than the maximum length of chains of specializations of $\alpha $ in $U,$ repeating the process finitely many times we get the result.



\qed \\

\begin{nrmk}\label{nrml a2}
{\em
In the semi-algebraic setting condition (A2) holds for the tildes of locally complete semi-algebraic spaces, and similarly it holds for the tildes of  locally definably compact definable spaces in o-minimal expansions of real closed fields. 

In arbitrary o-minimal structures  locally definably compact does not imply the existence of  completions in $\Df.$ For example, in ${\mathbb M}=({\mathbb R},<, 0,-,+, (q)_{q\in {\mathbb Q}}),$ if $X=(0,+\infty ),$ then $X$ is locally definably compact but $X$ has no completion in $\Df.$ Indeed, let  $i:X\to Z$ be a completion of $X$ in $\Df.$ Then because $Z$ is complete in $\Df,$ it is Hausdorff and definably compact (Fact \ref{fact def proper2}, ${\mathbb M}$ has definable Skolem functions by \cite[Chapter 6 (1.2)]{vdd}). Therefore the morphism $i:(0,+\infty )\to Z$ has the limit $\lim _{t\to +\infty }i(t),$ say $a,$ in $Z.$  Since $Z$ is definably normal (Fact \ref{fact def comp is normal}), by shrinking lemma, $a$ has an affine definably compact neighborhood. So we obtain a morphism $\iota : [d,+\infty )\to  [-b,b]$ in $\Df$ which is injective and  piecewise linear  by \cite[Chapter 1 (7.8)]{vdd}. This is absurd.

An analogue of  (A2) is required even in algebraic geometry  since when one has to define the proper direct image functor of a separated morphism of finite type of schemes one has to use Nagata's theorem on the existence of proper extensions of such morphisms (i.e.,  existence of completions as in (A2)).\\
}
\end{nrmk}

\begin{nrmk}\label{nrmk a2 x' norm}
{\em
If $X$ is an object of $\tbA$ and $i:X\to X'$ is a completion of $X$ in $\tbA,$ then $X'$ is normal. Indeed, under the inverse of the isomorphism $\Df \to \tDf,$ we have that $X'$ is complete in $\Df,$ therefore, by Fact \ref{fact def proper2}, $X'$ in $\Df$ is Hausdorff and definably compact, by \ref{fact def comp is normal}, $X'$ in $\Df$ is definably normal and by Fact \ref{fact ominspec closed point} $X'$ in $\tDf$ is normal.\\
}
\end{nrmk}

The first  consequence of  (A2) that we require is:\\

\begin{prop}\label{prop a2 c is normal}
Let $X$ be an object of $\tDf$ with a  completion in $\tDf.$ Then the family $c$ of  complete supports on $X$ is a normal and constructible family of supports on $X.$

In particular, if $f:X\to Y$ is a morphism in $\tbA$ and $\alpha \in Y,$ then the family $c$ of  complete supports on the fiber $f^{-1}(\alpha )$ is a normal and constructible family of supports on $f^{-1}(\alpha )$.
\end{prop}

\pf
Let $i:X\to X'$ be a completion of $X$ in $\tDf.$ Let $C\in c$ and  let $V\subseteq X$ be an open  neighborhood of $C$ in $X.$ By Remark \ref{nrmk cons sup} we may assume that $C$ and $V$ are constructible.  Then $i_{|C}:C\to X'$ is proper in $\tDf$ (Corollary \ref{cor basic compl} (2)) and so $i(C)$ is closed constructible in $X'.$ Since $i(V)$ is an open constructible neighborhood of $i(C)$ in $X'$ and $X'$ is normal (Remark \ref{nrmk a2 x' norm}),  by shrinking lemma  and Corollary \ref{cor basic compl} (1) and (3) applied to $i^{-1}:i(X)\to X,$ there is a complete constructible neighborhood $D$ of $C$ in $X$ such that $D\subseteq V.$

 Now let $f:X\to Y$ be a morphism in $\tbA$ and $\alpha \in Y.$ Let $a\models \alpha $ be a realization of $\alpha  $ and ${\mathbb S}$ a prime model of the first-order theory of ${\mathbb M}$ over $\{a\}\cup M$. By (A2) let $i:X\to X'$ be a completion of $X$ in $\tbA.$ By Proposition \ref{prop sep and proper in s}  $i^{{\mathbb S}}:X({\mathbb S})\to X'({\mathbb S})$ is a completion of $X({\mathbb S})$ in $\tDfs.$ 

 Since ${\mathbb S}$ has definable Skolem functions, by the first paragraph in $\tDfs,$ we see that the family $c$ of  complete supports on $X({\mathbb S})$ is a normal and constructible family of supports on $X({\mathbb S}).$ Since $(f^{{\mathbb S}})^{-1}(a)$ is closed constructible subset of $X({\mathbb S})$ it follows that the family $c$ of  complete supports on $(f^{{\mathbb S}})^{-1}(a)$ is a normal and constructible family of supports on $(f^{{\mathbb S}})^{-1}(a).$ Therefore, by Definition \ref{defn compact supp},  the family $c$ of  complete supports on the fiber $f^{-1}(\alpha )$ is a normal and constructible family of supports on $f^{-1}(\alpha ).$
\qed \\

In the semi-algebraic setting (A1)$^*$ is used  in combination with (A2) to verify the results \cite[Chapter II, Section 8,  Corollaries 8.3 and 8.4]{D3}.
Here, with slightly more work, namely Theorem \ref{thm CxN is normal}, from the weaker (A1) we are able to obtain the following suitable replacement of those results:



\begin{prop}\label{prop a2 c in s and proper}
 Let $f:X\to Y$ be a morphism in $\tbA$ and let  $\alpha  \in Y$ be closed.   If $K$ is a constructible  element of the family $c$ of  complete supports on the fiber $f^{-1}(\alpha ),$ then  there exists a  constructible and complete 
 neighborhood  $B$ of $K$ in  $X.$
\end{prop}

\pf
First we observe the we may assume that $Y$ is affine and normal. Indeed, by Proposition \ref{prop obj tA}, there is an open, affine, normal  constructible subset $U$ of $Y$ such that $\alpha \in U$ and $\alpha $ is closed in $U$. Let $V=f^{-1}(U)$ so that $f^{-1}(\alpha )=(f_{|V})^{-1}(\alpha ).$ If $B$ is a constructible and complete neighborhood of $K\subseteq (f_{|V})^{-1}(\alpha )$ in $V$ 
then $B$ is still a constructible and complete neighborhood of $K\subseteq f^{-1}(\alpha )$ in $X.$ 

To proceed note that  in $\Df$ the morphism $f:X\to Y$  is the composition of the graph embedding $X\to X\times Y: x\mapsto (x,f(x))$ and the projection $X\times Y \to Y.$  Since the graph embedding is a closed immersion ($f$ is continuous and both $X$ and $Y$ (in $\Df$) are Hausdorff, for example by (A1)), it follows that it is proper in $\Df$ (Proposition \ref{prop proper in tdef} (1)). Thus it is enough to show the  result for a projection $\pi :X\times Y \to Y$ in $\tbA.$ In this case, by (A2) and Corollary \ref{cor completion in C} (3), we have a commutative diagram
$$
\xymatrix{
X \times Y\ar [rd]_{\pi} \ar[r]^{\iota} & P \ar[d]^{\bar{\pi}} \\
& Y}
$$
of morphisms in $\tbA$ such that $\iota $ is an open immersion and $\bar{\pi }$ is  proper in $\tDf.$ Moreover, by Remark \ref{nrmk P},  $P=X'\times Y,$ $\iota =i\times {\rm id}$ with $i:X\to X'$ a completion of $X$ in $\tbA$ and  $\bar{\pi }:X'\times Y\to Y$ is the projection onto $Y.$  

Since $X'$ is complete in $\tDf,$ under the inverse of the isomorphism $\Df \to \tDf,$ we have that $X'$ is complete in $\Df,$ and so,  by Fact \ref{fact def proper2}, $X'$ in $\Df$ is Hausdorff and definably compact.  By Fact \ref{fact ominspec closed point} $Y$ in $\Df$ is definably normal.  Therefore, by Theorem \ref{thm CxN is normal}, $P=X'\times Y$ in $\Df$ is definably normal, and so by Fact \ref{fact ominspec closed point}, $P=X'\times Y$ in $\tDf$ is normal.

Let $a\models \alpha $ be a realization of $\alpha  $ and ${\mathbb S}$ a prime model of the first-order theory of ${\mathbb M}$ over $\{a\}\cup M$. Let $L$ be a constructible and  complete in $\tDfs$ subset  of $(\pi ^{{\mathbb S}})^{-1}(a)$ such that $K=r_{|}(L).$

By Corollary \ref{cor basic compl} (3) in $\tDfs,$ we have that   $\iota ^{{\mathbb S}}(L)$ is constructible and complete in $\tDfs$ subset of $(\bar{\pi }^{{\mathbb S}})^{-1}(a)\subseteq P({\mathbb S})$. In particular, by Corollary \ref{cor basic compl} (1) in $\tDfs,$  $\iota ^{{\mathbb S}}(L)$ is constructible and closed subset of $(\bar{\pi }^{{\mathbb S}})^{-1}(a)$.

From the commutative diagram,
$$
\xymatrix{
(\pi ^{{\mathbb S}})^{-1}(a) \ar [d]_{\iota ^{{\mathbb S}}} \ar[r]^{r_{|}} & \pi ^{-1}(\alpha ) \ar[d]^{\iota _{|}} \\
(\bar{\pi }^{{\mathbb S}})^{-1}(a) \ar[r]^{r_{|}} & \bar{\pi }^{-1}(\alpha )
}
$$
and the fact that the $r_{|}$ are homeomorphisms (Lemma \ref{lem bera}) it follows  that $\iota _{|}(K)=r_{|}(\iota ^{{\mathbb S}} (L))$ is a constructible  element of the family $c$ of  complete supports on the fiber $\bar{\pi }^{-1}(\alpha ).$ In particular, by the above  $\iota _{|}(K)$ is a closed constructible subset of  $\bar{\pi }^{-1}(\alpha)$ which is closed in $P$ ($\alpha$ is closed in $Y$).
Hence $\iota _{|}(K)$ is closed in $P$ and moreover $\iota _{|} (K) \cap (P \setminus \iota(X\times Y)) = \emptyset$.

Since $P$ is normal, by the shrinking lemma, we can find disjoint constructible open neighborhoods  $U$ and  $V$ of $\iota _{|} (K)$ and $P \setminus \iota(X\times Y)$ respectively in $P$. Set $C=P \setminus V$. Then  $\iota _{|} (K) \subset C\subset \iota (X\times Y)\subseteq P$ with $C$   constructible and complete in  $\tDf.$ By Proposition \ref{prop a2 c is normal} we can assume that $C$ is a constructible complete in $\tDf$  neighborhood of $\iota _{|} (K)$ in $\iota (X\times Y).$

Since $\iota :X\times Y\rightarrow \iota (X\times Y)$ is a homeomorphism
$B=\iota ^{-1}(C)$ is a constructible and complete 
neighborhood of $K$ in  $X\times Y.$
\qed \\




\begin{nrmk}\label{nrml a3}
{\em
In the semi-algebraic setting (A3) is proved for all tildes of semi-algebraic morphisms between the tildes of locally complete semi-algebraic spaces (\cite[Chapter II, Theorem 6.10]{D3} ). A similar proof works in o-minimal expansions of real closed fields as it is based on the definable/semi-algebraic triangulation theorems.

Condition  (A3) is required below for the base change formula (Proposition \ref{prop base change formula}), derived base change formula (Theorem \ref{thm der base change formula}), the K\"unneth formula (Theorem \ref{thm kunneth formula}) and dual base change formula (Proposition \ref{prop der base change formula dual}).

In each case where  (A3) is used that fact will be explicitly mentioned. We also note that although  (A3) looks rather technical  and in some specific cases  may follow from  simpler assumptions, we preferred to state exactly what is used in the proofs since in some situations we only need to apply the consequences of  (A3) to very specific morphisms and not all morphisms. \\
}
\end{nrmk}

We have:

\begin{expl}\label{expl A}
{\em
Categories $\tbA$ satisfying also (A3)  include the tildes of:
\begin{itemize}
\item[(i)] Regular, locally definably compact definable spaces in o-minimal expansions of real closed fields. ((A1) is by \cite[Chapter 10, Theorem (1.8)]{vdd} and  \cite[Chapter 6, Lemma (3.5)]{vdd};  (A2) and (A3) by \cite[Fact 4.7 and Corollary 4.8]{ep2}.)


\item[(ii)] Hausdorff locally definably compact definable spaces  in o-minimal expansions of ordered groups with  completions in $\Df.$ ((A1) by \cite[Chapter 6, Lemma (3.5)]{vdd};  (A2) by assumption and  (A3) by \cite[Theorem 4.5]{ep2}.)


\item[(iii)] Locally closed definable subsets of cartesian products of a given definably compact definable group  in an arbitrary o-minimal structure. (A1) by   \cite[Corollary 3.15]{oh5p}; (A2) follows since: (i) definably compact groups are definably normal (\cite[Corollary 2.3]{et} or \cite[Theorem 2.11]{emp}) and (ii) a locally closed definable subset of a cartesian product of a given definably compact definable group has a definable completion, namely its closure; (A3) by \cite[Theorem 1.1]{ep2} 
\end{itemize}
}
\end{expl}

In these examples, (A3) is obtained in \cite{ep2} after extending an invariance result for  closed and bounded definable sets in o-minimal expansions of ordered groups from \cite{bf}.\\

Finally we note that in topology we do not have the problem of fibers described above and what we need from (A0), (A1) and (A2)  holds on Hausdorff locally compact topological spaces. \\

\end{subsection}

\end{section}

\begin{section}{Proper direct image}\label{section proper dir im}
In this section  we will develop the theory of proper direct image in  full subcategories $\tbA$ of $\tDf$ introduced in Subsection \ref{subsection tbA}.\\

\noindent
{\bf Notation:} For the rest of the paper we let  $A$ be a commutative ring with unit.  If  $X$ is a topological space we denote by $\mod(A_X)$ the category of sheaves of $A$-modules on $X$ and we call its objects $A$-sheaves on $X$.

Since $\tDf $ is a subcategory of ${\rm Top}$, if $X$ is an object of $\tDf$ then we have the classical operations
$$\ho_{A_X}(\bullet ,\bullet ), \,\, \bullet \otimes _{A_X}\bullet , \,\, f_*,\,\, f^{-1}, \,\, (\bullet )_Z, \,\,\Gamma _Z(X;\bullet ), \,\, \Gamma (X; \bullet )$$
on $A$-sheaves on $X,$ where $Z\subseteq X$ is a locally closed subset. Below we may use freely these operations and refer the reader to \cite[Chapter II, Sections 2.1 - 2.4]{ks1} for the details on  sheaves on topological spaces and on the properties of  these basic operations.

Later in this section we may also  use the derived versions of many of the properties relating the above operations  and we refer to reader to \cite[Chapter II, Section 2.6]{ks1} for details. The reader can also see the classical references \cite{b}, \cite{g} and \cite{i} for similar  details on sheaves on topological spaces.\\

\noindent
{\bf Convention:}
For the rest of the paper we will work in the category $\tDf$ (resp. $\tbA$ or $\tDfs$) and omit the tilde on objects and morphisms.  We will say ``proper'' (resp. ``separated'' or ``complete'') instead of ``proper in $\tDf$'' (resp. ``separated in $\tDf$'' or ``complete in $\tDf$''). 

We will assume that all morphism and objects that we consider in $\tDf$ are separated. This is the case for morphisms and objects in $\tbA$ (by (A2) all objects of $\tbA$ are tildes of open definable subsets of objects complete in $\Df$ (i.e. Hausdorff and definably compact) and so are Hausdorff). 

If $X$  is an object of $\tDf$ (resp. $\tbA$), then $\op (X)$ denotes the category of open subsets of $X$ with inclusions and $\opc (X)$  is the full sub-category of constructible open subsets of $X.$\\


\begin{subsection}{Extended summary}
Since our definition of the o-minimal proper direct image functor $f_{\dex }$ associated to a morphism $f:X\to Y$ in $\tDf$ is similar to the definition of the topological proper direct image functor $f_!$, we only replace proper in ${\rm Top}$ by proper in $\tDf$,  we follow closely the proofs in topology (\cite[Chapter II, Sections 2.5 and 2.6]{ks1}) but we have to deal with: (i)  the fact that the fibers $f^{-1}(\alpha )$ are often not objects of $\tDf$; (ii) the fact that cartesian squares in $\tDf$ are not cartesian squares in ${\rm Top}$; (iii) the fact that objects of $\tDf$ are not locally compact topological spaces.\\

The most technical result is the Fiber formula (Corollary \ref{cor fibers of f!}) which is a consequence of  the Relative fiber formula (Proposition \ref{prop f! and c}) in whose proof we use heavily the consequences of our assumptions (A0), (A1) and (A2) discussed in the Subsection \ref{subsection tbA}. Our methods are different from those in semi-algebraic case: the fiber formula there (\cite[Chapter II, Corollary 8.8]{D3}) is obtained as a consequence of a Base change formula (\cite[Chapter II, Theorem 8.7]{D3}) and no semi-algebraic analogue of the Relative fiber formula is proved by Delfs.\\

We then proceed with the theory of $f$-soft sheaves which gives us the $f_{\dex}$-injective objects (Proposition \ref{prop f-soft are f!-inj}) and later the bound on the cohomological dimension of $f_{\dex }$ (Theorem \ref{thm cohomo  dim f!}). The $f$-soft sheaves are those whose restriction to fibers $f^{-1}(\alpha )$ are $c$-soft. We show our results about $f$-soft sheaves after we remark (Remark \ref{nrmk met}) that by our assumptions (A0), (A1) and (A2) we can always reduce to the case where $\alpha $ is a closed point in which case $c$ is a normal and constructible family of supports on $f^{-1}(\alpha )$ (Proposition \ref{prop a2 c is normal}) and  the results follow, via $\tilde{(f^{{\mathbb S}})^{-1}(a)}\simeq f^{-1}(\alpha )$,  from similar results already proved in \cite[Section 3]{ep1}. This theory has not been considered in the semi-algebraic case. \\

The Projection formula (Proposition \ref{prop proj formula}) is obtained in a standard way but the Base change formula (Proposition \ref{prop base change formula}) requires assumption  (A3) and, after we prove a technical lemma  (Lemma \ref{lem homeo fibers final cart sq}), we use (A3) more or less as the corresponding fact (\cite[Chapter I, Theorem 6.10]{D3}) is used in  the prove of the semi-algebraic Base change formula (\cite[Chapter II, Theorem 8.7]{D3}).  Using the $f$-soft resolutions we obtain in a standard way the Derived projection formula (Theorem \ref{thm der proj formula}) and the Derived base change formula (Theorem  \ref{thm der base change formula}) after we prove a technical lemma  (Lemma \ref{lem f-acy to f'-acy}) requiring assumption  (A3). This part was also not considered in the semi-algebraic case.

The complication with the Base change formula and the Derived base change formula and the need to use (A3) arises from the fact that in a cartesian square
$$
\xymatrix{X' \ar[r]^{f'} \ar[d]^{g'} & Y' \ar[d]^g \\
X \ar[r]^f &Y.}
$$
in $\tDf,$ if $\gamma \in Y'$, then we may have $\tilde{(f'^{{\mathbb S}})^{-1}(a)}\simeq f'^{-1}(\gamma)$ and $\tilde{(f^{{\mathbb K}})^{-1}(b)}\simeq f^{-1}(g(\gamma ))$ with ${\mathbb K}={\mathbb M}$ and ${\mathbb S}\neq {\mathbb M}.$ So after we identify $(f'^{{\mathbb S}})^{-1}(a)$ with $(f^{{\mathbb S}})^{-1}(b)$ via $g'^{{\mathbb S}}$ we need to know that $H^*_c(f^{-1}(b); F_{|f^{-1}(b)})\simeq H^*_c((f^{{\mathbb S}})^{-1}(b); F({\mathbb S})_{|(f^{{\mathbb S}})^{-1}(b)}).$ \\

\end{subsection}

\begin{subsection}{Proper direct image}\label{subsection proper dir im}
Here we define the proper direct image functor and prove some of its basic properties.\\

\begin{defn}\label{defn  dir image proper }
{\em
Let $f:X \to Y$ be a morphism in $\tDf$ and let $F \in \mod(A_{X})$. The {\em proper direct image} is the subsheaf of $f_*F$ defined by setting for $U \in \opc(Y),$ an open constructible subset of $Y,$
$$
\Gamma(U;f_{\dex }F) = \lind Z \Gamma_Z(\imin f(U);F),
$$
where $Z$ ranges through the family of closed constructible subsets of $\imin f(U)$ such that $f_{|Z}:Z \to U$ is  proper.
}
\end{defn}

From the definition we have:

\begin{nrmk}\label{nrmk basic fdex}
{\rm
The functor $f_{\dex }$ is clearly left exact; if $f:X\to Y$ is proper, then $f_{\dex }=f_*$; if $i:Z\to Y$ is the inclusion of a locally closed subset, then $i_{\dex }$ is the extension by zero functor (i.e. $(i_{\dex}F)_{\alpha }\simeq F_{\alpha }$ or $0$ according to $\alpha \in Z$ or $\alpha \not \in Z$) and $(\bullet )_Z=i_{\dex }\circ i^{-1}(\bullet )$. Compare with \cite[Chapter II, Proposition 2.5.4]{ks1}.
}
\end{nrmk}

\begin{nrmk}\label{nrmk ax!}
{\em
If we consider  the morphism $a_X:X \to \{\pt \}$ to a point in $\tDf$,  then we have
$$(a_{X\dex }F)_{\pt}\simeq \Gamma (\pt; a_{X\dex }F)=\Gamma_c(X;F).$$
}
\end{nrmk}

\begin{prop}\label{prop compos of !}
Let $f:X \to Y$ and  $g:Y \to Z$ be  morphisms in $\tDf$. Then  $(g \circ f)_{\dex } \simeq g_{\dex } \circ f_{\dex }.$
\end{prop}

\pf
Indeed, if $V \in \opc(Z)$, then a section of $\Gamma(V,g_{\dex }\circ f_{\dex }F)$ is represented by $s \in \Gamma(\imin f(\imin g(V));F)$ such that $\supp(s) \subset \imin f(S) \cap Z$, where $Z,S$ are closed constructible subsets of $\imin f(\imin g(V))$ and $\imin g(V)$ respectively and such that $f_{|Z}:Z \to \imin g(V)$, $g_{|S}:S \to V$ are  proper. The set $Z \cap \imin f(S)$ is closed constructible and the restriction $(g \circ f):Z \cap \imin f(S) \to V$ is  proper (Proposition \ref{prop proper in tdef} (1) and (2)). Conversely if $V \in \opc(Z)$, then a section of $\Gamma(V, (g\circ f)_{\dex }F)$ is represented by  $s \in \Gamma_Z(\imin f(\imin g(V));F)$ where $Z$ is closed constructible subset of  $\imin f(\imin g(V))$ such that $(g \circ f)_{|Z}:Z \to V$ is  proper. But  then $f_{|Z}:Z \to \imin g(V)$ is proper and $g_{|f(Z)}:f(Z) \to V$ is  proper  (Proposition \ref{prop proper in tdef} (5)).
 \qed \\

In particular, we have:

\begin{nrmk}\label{nrmk f! and aX!}
{\em
For $f:X \to Y$ a morphism in $\tDf$, since $a_X = a_Y \circ f$, we have $\Gamma _c(Y;f_{\dex }F)\simeq \Gamma_c(X;F).$
}
\end{nrmk}

\begin{prop}\label{prop f! and lim} Let $f:X \to Y$ be a morphism in $\tDf$. The functor $f_{\dex }$ commutes with filtrant inductive limits.
\end{prop}

\pf
Let $(F_i)_{i \in I}$ be a filtrant inductive system in $\mod(A_X)$. Since open constructible sets are quasi-compact, then for any  $V \in \opc(X)$ and any closed constructible subset $S$ of $V$ the functors $\Gamma(V;\bullet)$ and $\Gamma(V \setminus S;\bullet)$ commute with filtrant $\Lind$ (\cite[Remark 2.7]{ep1}). Hence $\Gamma_S(V;\bullet)$ commutes with filtrant $\Lind$. Therefore we have
\begin{eqnarray*}
\Gamma(U;\lind i f_{\dex }F_i) & \simeq & \lind i \Gamma(U;f_{\dex }F_i) \\
& \simeq & \lind {i,Z} \Gamma_Z(f^{-1}(U);F_i) \\
& \simeq & \lind Z \Gamma_Z(f^{-1}(U);\lind i F_i) \\
& \simeq & \Gamma(U;f_{\dex }\lind i F_i),
\end{eqnarray*}
where $Z$ ranges through the family of closed constructible subsets of $\imin f(U)$ such that $f_{|Z}:Z \to U$ is  proper. \qed \\

The following lemma follows immediately from the definition of proper direct image and we leave the details to the reader

\begin{lem}\label{lem f!  reduc open}
Let $f:X \to Y$ be a  morphism in $\tDf$ and $W$ an open subset  of $Y$. Consider the commutative diagram
$$
\xymatrix{
f^{-1}(W)  \ar@{^{(}->}[r]^{j}  \ar[d]^{f_|} & X \ar[d]^{f} \\
W \ar@{^{(}->}[r]^{i} & Y.
}
$$
If $F\in \mod (A_X),$ then $(f_{\dex }F)_{|W}\simeq (f_|)_{\dex }(F_{|f^{-1}(W)}).$\\
\end{lem}

\begin{prop}[Relative fiber formula]\label{prop f! and c}
Let $f:X \to Y$ and  $g:Y \to Z$ be  morphisms in $\tbA$, $\beta \in Z$ and let $F$ be a sheaf  in $\mod(A_{X}).$ If $K \in c$ in $g^{-1}(\beta ),$  then $\Gamma(K; f_{\dex }F) \simeq \Gamma_{c}(f^{-1}(K);F)$.
\end{prop}

\pf
First we show that we may assume without loss of generality that $\beta $ is a closed point in $Z.$

By Proposition \ref{prop obj tA},  there exists $Z'  \in \opc(Z)$ with $\beta \in Z'$ such that $\beta $ is closed in $Z'$. Let $Y'=g^{-1}(Z')$ and $X'=f^{-1}(Y').$ Let also $f'=f_{|X'}$ and $g'=g_{|Y'}.$  By Lemma \ref{lem f!  reduc open} we have $(f_{\dex }F)_{|K}=((f_{\dex }F)_{|Y'})_{|K})\simeq (f'_{\dex }F_{|X'})_{|K}),$ and  since $f'^{-1}(K)=f^{-1}(K)$, we have also $F_{|f^{-1}(K)}=F_{|f'^{-1}(K)}).$ Thus after replacing $Z$ (resp., $Y$ and $X$) by $Z'$ (resp. $Y'$ and $X'$) and $f$ (resp., $g$) by $f'$ (resp., $g'$), we may assume without loss of generality that $\beta $ is a closed point of $Z.$ \\

Let $b\models \beta $ be a realization of $\beta $ and ${\mathbb S}$ a prime model of the first-order theory of ${\mathbb M}$ over $\{b\}\cup M.$ Then we have the following commutative diagram
\begin{equation*}\label{+++}
\xymatrix{
(f ^{{\mathbb S}})^{-1}((g^{{\mathbb S}})^{-1}(b))  \ar[r]^{r_|}  \ar[d]_{(f_|)^{{\mathbb S}}} & \imin f(\imin g (\beta )) \ar[d]^{f_|} \\
(g^{{\mathbb S}})^{-1}(b) \ar[r]^{r_|} & \imin g(\beta )
}
\end{equation*}
given by Lemma \ref{lem bera} (where we have omitted the tildes and used the fact that $(g\circ f)^{{\mathbb S}}=(g^{{\mathbb S}})\circ (f^{{\mathbb S}})$).

By (A0), (A2)  and Corollary \ref{cor completion in C} (2), we  a commutative diagram
$$
\xymatrix{ X \ar [d]_{f} \ar[r]^{i _X}& X' \ar[d]^{f'} \\
Y \ar [d]_{g} \ar[r]^{i _Y}& Y' \ar[d]^{g'} \\
Z  \ar[r]^{i_Z} & Z'}
$$
of completions in $\tbA.$  \\ 

Let $\Phi_{f}=\{A: A\subseteq X$ is closed, $A\subseteq B$ for some closed constructible subset $B$ of $X$ such that $f_{|B}:B\to Y$ is proper in $\tDf$$\}.$ Then:

\begin{clm}\label{clm fiber1}
$$\Phi _f\cap f^{-1}(K)=c_{|f^{-1}(K)}.$$ 
In particular, if $Y=Z,$  $g={\rm id}$ and $\alpha \in Y$ is closed, then for $c$ the family of complete supports of $f^{-1}(\alpha ),$ we have
$$c=\Phi_{f}\cap f^{-1}(\alpha ).$$
\end{clm}

Let $Z \in \Phi_f.$  We may assume that $Z$ is a closed constructible subset of $X$ such that $f_{|Z}:Z\to Y$ is proper in $\tDf.$ Then, $f^{{\mathbb S}}_{|Z({\mathbb S})}:Z({\mathbb S}) \to Y({\mathbb S})$ is proper in $\tDf ({\mathbb S})$ (Proposition \ref{prop sep and proper in s}). Since $K \in c$ in $g^{-1}(\beta )$, by definition, there is a $C\subseteq (g^{{\mathbb S}})^{-1}(b)$ constructible and complete in $(g^{{\mathbb S}})^{-1}(b)$ such that $K\subseteq r_|(C).$ By Corollary \ref{cor basic compl} (4), $(f^{{\mathbb S}})^{-1}(C)\cap Z({\mathbb S})$ is a constructible and complete in $(f^{{\mathbb S}})^{-1}((g^{{\mathbb S}})^{-1}(b)).$ Therefore, $r_|((f^{{\mathbb S}})^{-1}(C)\cap Z({\mathbb S}))=Z \cap f^{-1}(r_|(C))$ is a constructible and complete in  $f^{-1}(g^{-1}(\beta )).$ Since $Z\cap f^{-1}(K)$ is a closed subset of $Z \cap f^{-1}(r_|(C))$ we have that $Z\cap f^{-1}(K) \in c$ in $f^{-1}(g^{-1}(\beta )).$

Conversely, let $Z \in c_{|f^{-1}(K)}$, then there is $C$ a constructible complete subset of  $\imin f(\imin g(\beta ))$ such that $Z \subseteq C\subseteq \imin f(K).$ Since $\beta $ is closed in $Z,$ it follows from Proposition  \ref{prop a2 c in s and proper} applied to $g\circ f$ that there is a $B\subseteq X$  constructible and complete in $X$ such that $C\subseteq B.$ Hence $f_{|B}$ is proper (Corollary \ref{cor basic compl} (2)). So $C\subseteq B\cap \imin f(K)\in \Phi _f\cap \imin f(K).$
\qed \\

We proceed with the proof of the Proposition. We have
$$\Gamma(K;f_{\dex }F) \simeq \lind U \Gamma(U;f_{\dex }F) \simeq \lind {U,Z}\Gamma_Z(\imin f(U);F),$$
 where $U$ ranges through the family of open constructible neighborhoods of $K$ and $Z$ is closed constructible in $\imin f(U)$ and such that $f_{|Z}:Z \to U$ is proper.

 By Claim \ref{clm fiber1} we can conclude the proof  after the following  two steps.

\begin{clm}
We have an isomorphism
$$\lind {U,Z}\Gamma_Z(\imin f(U);F) \simeq \lind U \Gamma_{\Phi_f \cap \imin f(U)}(\imin f(U);F),$$ where $U$ ranges through the family of open constructible neighborhoods of $K$ and $Z$ is a closed constructible subset of $\imin f(U)$ and such that $f_{|Z}:Z \to U$ is  proper.
\end{clm}

Let $s \in \Gamma(\imin f(U);F)$ whose support is contained in  $Z$ a closed constructible subset of  $\imin f(U)$ such that $f_{|Z}:Z \to U$ is  proper. Since $K \in c$ in $g^{-1}(\beta )$ and $\beta $ is closed in $Z,$ it follows from Proposition  \ref{prop a2 c in s and proper} applied to $g$ that there is  $B\subseteq Y$ a  constructible and complete in $Y$ neighborhood of  $K$ in $Y.$ Since $K$ is closed in $Y$ ($g^{-1}(\beta )$ is closed in $Y$), it is closed in $B.$ Since $i_{Y|B}:B\to Y'$ is proper (Corollary \ref{cor basic compl} (2)), $i_Y(K)=i_{Y|B}(K)$ is a closed subset of $i_Y(B)$ (Proposition \ref{prop closed on cons}). Since $i_Y(B)$  is a closed constructible subset of $Y'$ (Corollary \ref{cor basic compl} (1)), it follows that $i_Y(K)$ is a closed subset of $Y'.$ Since $Y'$ is normal (Remark \ref{nrmk a2 x' norm}), there exists an open  constructible neighborhood $V$ of $K$ such that $\bar{V} \subset U$. Then $f_{|Z \cap \imin f(\overline{V})}:Z \cap \imin f(\overline{V}) \to \overline{V}$ is proper and composing with the inclusion $\overline{V} \hookrightarrow Y$, the morphism $f_{|Z \cap \imin f(\overline{V})}:Z \cap \imin f(\overline{V}) \to Y$ is proper (Proposition \ref{prop proper in tdef} (1) and (2)). Hence $Z \cap \imin f(\overline{V}) \in \Phi_f$. The support of the restriction of $s$ to $\imin f(V)$ is contained in $Z \cap \imin f(V) = Z  \cap \imin f(\overline{V}) \cap \imin f(V) \in \Phi_f \cap \imin f(V)$ and the result
follows. \qed\\

\begin{clm}
We have an isomorphism
$$\lind U \Gamma_{\Phi_f \cap \imin f(U)}(\imin f(U);F) \simeq \Gamma_{\Phi_f \cap \imin f(K)}(\imin f(K);F),$$
where $U$ ranges through the family of open constructible neighborhoods of $K$.
\end{clm}

Let $s \in \Gamma(\imin f(K);F)$ and let $Z \in \Phi_f$ closed constructible subset of  $\imin f(U)$  containing its support. The section $s$ is represented by a section $t \in \Gamma(W;F)$ for some open constructible neighborhood $W$ of $\imin f(K)$. Since $Z \in \Phi_f$ and $K \in c$ in $g^{-1}(\beta ),$ we have  $Z\cap f^{-1}(K) \in c$ in $f^{-1}(g^{-1}(\beta ))$ (Claim \ref{clm fiber1}). Since $\beta $ is closed in $Z,$  it follows from Proposition \ref{prop a2 c in s and proper} applied to $g\circ f$ that there is  $B\subseteq X$ a  constructible and complete in $X$ neighborhood of  $Z\cap f^{-1}(K)$ in $X.$ Since $Z\cap f^{-1}(K)$ is closed in $X$ ($f^{-1}(g^{-1}(\beta ))$ is closed in $X$), it is closed in $B.$ Since $i_{X|B}:B\to X'$ is proper (Corollary \ref{cor basic compl} (2)), $i_X(Z\cap f^{-1}(K))=i_{X|B}(Z\cap f^{-1}(K))$ is a closed subset of $i_X(B)$ (Proposition \ref{prop closed on cons}). Since $i_X(B)$  is a closed constructible subset of $X'$ (Corollary \ref{cor basic compl} (1)), it follows that $i_X(Z\cap f^{-1}(K))$ is a closed subset of $X'$ .

Since $X'$ is normal (Remark \ref{nrmk a2 x' norm}) there exists an open constructible neighborhood $V$ of $Z \cap \imin f(K)$ with $\overline{V} \subset W$.  We have that $f'_{|\overline{V}}$ is proper, so $f_{|\overline{V}}$ is proper (Proposition \ref{prop proper in tdef} (5) (i)) and $\overline{V} \in \Phi_f$. Since $f'{}^{-1}(K)=f^{-1}(K)$, we have $f'{}^{-1}(K) \cap ((\overline{V}\cap {\rm supp}(t))\setminus V)=\emptyset $ and so $K\cap f'((\overline{V} \cap \supp(t) )\setminus V)=\emptyset .$ Since $K$ and $f'((\overline{V} \cap \supp(t) )\setminus V)$ are closed subsets of $Y'$ (recall that $f'_{|\overline{V}}$ is closed) and $Y'$ is normal, by  shrinking lemma,    there exists an open constructible subset $U$ of $Y$ with $K\subseteq U$ such that $\imin f(U) \cap \bar{V} \cap \supp (t) \subset V$. Let us define $\widetilde{t} \in \Gamma(\imin f(U);F)$ by
\begin{eqnarray*}
\widetilde{t}_{|\imin f(U) \setminus (\overline{V}\cap \supp (t))} & = & 0, \\
\widetilde{t}_{|\imin f(U) \cap V}
& = & t_{|\imin f(U) \cap V}.
\end{eqnarray*}
By construction $\widetilde{t}_{|\imin f(K)}=s$. Note that $\supp(\widetilde{t})$ is contained in $\imin f(U) \cap \overline{V}$ and $\overline{V} \in \Phi_f$.
\qed \\

Setting $Y=Z$ and $g=\id$ we obtain

\begin{cor}[Fiber formula] \label{cor fibers of f!}
Let $f:X \to Y$ be a morphism in $\tbA$ and let $F$ be a sheaf  in $\mod(A_{X})$. Let $\alpha \in Y$.  Then $(f_{\dex }F)_\alpha \simeq \Gamma_c(f^{-1}(\alpha);F)$.\\
\end{cor}

We end the subsection  comparing the o-minimal proper direct image functor $f_{\dex}$ with Verdier's (\cite{ver}) topological proper direct image functor $f_!$ and Kashiwara and Schapira (\cite{ks2}) sub-analytic proper direct image functor $f_{!!}.$

To make the comparison more clear we will use the isomorphism $\mod (A_{\tilde{X}})\simeq \mod (A_{X_{\df}})$ and the fact that $f:X\to Y$ is proper in $\Df$ if and only if $\tilde{f}:\tilde{X}\to \tilde{Y}$ is proper in $\tDf$ to introduce $f_{\dex}: \mod (A_{X_{\df}})\to \mod (A_{Y_{\df}})$ as $\tilde{f}_{\dex}: \mod (A_{\tilde{X}})\to \mod (A_{\tilde{Y}}).$ Recall that given an object $X$ of $\Df$ the o-minimal site $X_{\df}$ on $X$ is the category $\op (X_{\df})$ whose objects are open (in the topology of $X$ mentioned above) definable subsets of $X$, the morphisms are the  inclusions and the admissible covers $\cov (U)$ of $U\in \op (X_{\df})$ are covers by open definable subsets of $X$ with  finite sub-covers.\\

We now explain the relationship between the o-minimal  proper direct image functor and Verdier's (\cite{ver}) proper direct image functor.  

\begin{nrmk}\label{nrmk comp top}
{\em
Consider the category $\Df$ associated to  an o-minimal expansion   ${\mathbb M}=({\mathbb R}, <, (c)_{c\in {\mathcal C}},  (f)_{f\in {\mathcal F}}, (R)_{R\in {\mathcal R}})$  of the ordered set of real numbers.

For a continuous map $f:X\to Y$ between locally compact topological spaces  Verdier's proper direct image functor $f_!$ is given by: for $G \in \mod(A_{X})$ and $U \in \op(Y)$ we have
$$
\Gamma(U;f_{!}G) = \lind Z \Gamma_Z(\imin f(U);G),
$$
where $Z$ ranges through the family of closed  subsets of $\imin f(U)$ such that $f_{|Z}:Z \to U$ is   proper (in ${\rm Top}$). 

If $X$ is an object of $\Df$ then $X$ is also a topological space  (with the usual topology generated by open definable subsets) and we have  the natural morphism of sites  $\rho : X\to X_{\df}$  induced by the inclusion $\op(X_{\rm def}) \subset \op(X)$.

Since if $f_{|Z}:Z\to Y$ is proper in $\Df$ then  $f_{|Z}:Z\to Y$ is proper in ${\rm Top}$ (\cite[Theorem 4.11]{emp}), we have $\Gamma(U;f_{\dex }\circ \rho _*F)\subseteq \Gamma(U;\rho_*\circ f_{!}F)$ for $U \in \op(X_{\rm def}).$
However, in general, we have $\rho_*\circ f_{!}\neq f_{\dex }\circ \rho_*.$

For example, suppose that $X={\mathbb R}^2$, $Y={\mathbb R},$ $f:{\mathbb R}^2\to {\mathbb R}$ is the projection onto the first coordinate  and  $U=(0,+\infty )$. Let $\phi :(0, +\infty )\to (0, +\infty )$ be given by $\phi (t)=\frac{1}{t}$ (resp. $\phi (t)=\frac{1}{\ln (t)}$) if ${\mathbb M}$ is semi-bounded (\cite{e1}) (resp. polynomially bounded (\cite{dmi})). Set $S=\{(x,y)\in (0,+\infty )\times {\mathbb R}: y=\phi (x))\}$ and consider the sheaf $A_S$. Let $s \in \Gamma (\imin f(U);A_S)$. Then $f:\supp (s)\to U$ is proper in ${\rm Top}$ (it is a homeomorphism into its image), so $s\in \Gamma(U;f_{!}A_S)=\Gamma(U;\rho_*\circ f_!A_S),$ but $s\notin \Gamma(U; f_{\dex }\circ \rho _*A_S)$ since there is no closed definable   subset $Z$  of $\imin f(U)$ such that $f_{|Z}:Z \to U$ is  proper in $\Df$ and $S\subseteq Z$.\\

}
\end{nrmk}

Now we  explain the relationship between the o-minimal proper direct image functor and Kashiwara-Schapira's (\cite{ks2}) sub-analytic proper direct image functor.

\begin{nrmk}\label{nrmk comp ks1}
{\em
Consider the category $\Df$ associated to the o-minimal structure ${\mathbb R}_{\rm an}=({\mathbb R},<, 0 ,1,+,\cdot, (f)_{f\in {\rm an}})$ - the field of real numbers expanded by restricted globally analytic functions (\cite{dmi}). As explained in \cite{dmi}, in this case, $\Df$ is the category of globally sub-analytic spaces with continuous maps with globally sub-analytic graphs.

If   $X$  is a  real analytic manifold, then we can equip $X$ with its sub-analytic site, denoted $X_{sa},$ where the objects are open sub-analytic subsets $\op (X_{sa})$ and $S\subset \op (X_{sa})\cap \op (U)$ is covering of $U\in \op (X_{sa})$ if for any compact subset $K$ of $X$ there is a finite subset $S_0\subseteq S$ such that $K\cap \bigcup _{V\in S_0}V=K\cap U,$ see \cite{ks2}.
In this context we have a sub-analytic proper direct image functor $f_{!!}$ given by: for   $F \in \mod(A_{X_{sa}})$ and $U \in \op(Y_{sa})$ we have
$$
\Gamma(U;f_{!!}F) = \lind K \Gamma_{K\cap \imin f(U)}(\imin f(U);F),
$$
where $K$ ranges through the family of compact sub-analytic subsets of $X.$ We recall that this is a direct  construction of Kashiwara and Schapira (\cite{ks2}) sub-analytic proper direct image functor originally constructed as a special case of a more general construction within  the theory of ind-sheaves (the direct construction for subanalytic sheaves was performed in \cite{lucap}).

If $X$ is also globally sub-analytic (i.e. a definable subset in  the o-minimal structure ${\mathbb R}_{\rm an}$) we have  the natural morphism of sites $\nu :X_{sa} \to X_{\rm def}$ induced by the inclusion $\op(X_{\rm def}) \subset \op(X_{sa})$.

Since a compact sub-analytic subset is a compact globally sub-analytic subset (hence a definably compact definable subset in  the o-minimal structure ${\mathbb R}_{\rm an}$), we have $\Gamma(U;\nu_*\circ f_{!!}F)\subseteq \Gamma(U;f_{\dex }\circ \nu_*F)$ for $U \in \op(X_{\rm def}).$
However, in general, we have $\nu_*\circ f_{!!}\neq f_{\dex }\circ \nu_*.$

For example, suppose that $X={\mathbb R}^2$, $Y={\mathbb R},$ $f:{\mathbb R}^2\to {\mathbb R}$ is the projection onto the first coordinate  and  $U=(0, 1)$.  Let $\phi :(0, 1)\to (0, +\infty )$ be given by $\phi (t)=\frac{1}{x(1-x)}.$ Set $S=\{(x,y)\in (0, 1)\times {\mathbb R}: y=\phi (x))\}$ and consider the sheaf $A_S$. Let $s \in \Gamma (\imin f(U);A_S).$ 
Then $f_|:\supp (s)\to U$ is definably proper (and so proper in $\Df$), so $s\in \Gamma(U;f_{\dex}\circ \nu_*A_S),$ but $s\notin \Gamma(U;\nu_*\circ f_{!!}A_S)=\Gamma(U;f_{!!}A_S) $ since there is no compact (globally) sub-analytic subset $K$  of $X$ such that  $\supp (s)\subseteq K$. \\
}
\end{nrmk}

\begin{nrmk}\label{nrmk comp ks2}
{\em
The o-minimal proper direct image functor $f_{\dex}$ in comparison with the sub-analytic proper direct image functor $f_{!!}$ seems to be more natural since it commutes with restrictions in the sense of Lemma \ref{lem f!  reduc open}, as the classical Verdier proper direct image functor $f_!$ does in topological spaces. The property of that lemma is not satisfied by the functor $f_{!!}$ since there are less compact subsets on $f^{-1}(W)$  than there are  intersections of $f^{-1}(W)$ with compact subsets of $X.$  \\ 
}
\end{nrmk}

\end{subsection}

\begin{subsection}{$f$-soft sheaves}\label{subsection f-soft}

Here we introduce the  $f$-soft sheaves  where $f:X\to Y$ is a morphism in $\tDf$. We show that the full additive subcategory of $f$-soft sheaves is $f_{\dex }$-injective  and    is stable under $\Lind $  and its flat objects are also stable under $\bullet \otimes F$ for all $F\in \mod (A_X)$. \\ 


Consider a fiber $f^{-1}(\alpha)$ of a morphism $f:X \to Y$ in $\tDf$ and let $c$ be the family of complete supports on $f^{-1}(\alpha)$. Recall that a  sheaf $F$ on $f^{-1}(\alpha)$ is {\it $c$-soft} if and only if  the restriction $\Gamma (f^{-1}(\alpha); F)\to \Gamma (K; F)$ is surjective for every   $K\in c .$ \\

\begin{defn}
{\em
Let $f:X \to Y$ be a morphism in $\tDf$ and let $F$ be a sheaf  in $\mod(A_{X})$. We say that $F$ is {\em $f$-soft} if for any $\alpha \in Y$ the sheaf $F_{|f^{-1}(\alpha)}$ in $\mod(A_{f^{-1}(\alpha)})$ is $c$-soft.
}
\end{defn}

By  Remark \ref{nrmk  ax!} we have:

\begin{nrmk}\label{nrmk c-soft ax!}
{\em
Let $a_X:X \to {\rm pt}$ be  the morphism to a point in $\tDf$ and let $F$ be a sheaf  in $\mod(A_{X})$. Then $F$ is $a_X$-soft if and only if it is $c$-soft.
}
\end{nrmk}

\begin{nrmk}\label{nrmk met}
{\em
Let $f:X\to Y$ be a morphism in $\tbA.$ To  prove a property of $f$-soft sheaves  in $\mod(A_{X})$ we have to  take an arbitrary $\alpha \in Y$ and prove the corresponding property for  $c$-soft sheaves  in $\mod(A_{f^{-1}(\alpha )}).$    We will use Proposition \ref{prop obj tA} to be able to replace $f:X\to Y$ by a suitable morphism $f':X'\to Y'$ in $\tbA$ such that $\alpha \in Y'$ and is closed in $Y'.$ After that we take  $a\models \alpha $ a realization of $\alpha $, ${\mathbb S}$ a prime model of the  first-order theory of ${\mathbb  M}$ over $\{a\}\cup M$ and  $r_{|}: (f'^{\mathbb S})^{-1}(a)\to f'^{-1}(\alpha )$  the homeomorphism  of Lemma \ref{lem bera}. It follows that  $(f'^{{\mathbb S}})^{-1}(a )$ is an object of $\tDfs$ and $c$ is a normal and constructible family of supports on $(f'^{{\mathbb S}})^{-1}(a )$ because $f':X'\to Y'$ is a morphism in $\tbA$ and  by Proposition \ref{prop a2 c is normal}, $c$ is a normal and constructible family of supports on $f'^{-1}(\alpha ).$ Therefore, we will be able to transfer results for $c$-soft sheaves on o-minimal spectral spaces (\cite[Section 3]{ep1}) to $c$-soft sheaves on the fiber $f^{-1}(\alpha )$ since all we need is that  $c$ is a family of normal and constructible  supports. \\
}
\end{nrmk}

\begin{nrmk}\label{nrmk k field}
{\em
In the paper \cite{ep1} we assumed that $A$ is a field, but that was only used to ensure that $\bullet \otimes G\simeq G\otimes \bullet $, for $G\in \mod(A_{X})$,  is exact.
For this reason, our results here about $\bullet \otimes G\simeq G\otimes \bullet $ will come with the flatness assumption.\\
}
\end{nrmk}

The first application of the above method is:

\begin{prop} \label{prop lim f-soft}
Let $f:X \to Y$ be a morphism in $\tbA$. Then filtrant inductive limits of $f$-soft sheaves in ${\rm Mod}(A_{X})$  are $f$-soft.
\end{prop}

\pf
Let $(F_i)_{i\in I}$ be a filtrant inductive family of  $f$-soft sheaves in ${\rm Mod}(A_{X}).$ We have to show that $\varinjlim F_i$ is an $f$-soft sheaf in $\mod (A_X),$  i.e., for every $\alpha \in Y$ the sheaf $(\varinjlim F_i)_{|f^{-1}(\alpha)}=\varinjlim (F_{i|f^{-1}(\alpha )})$ in $\mod(A_{f^{-1}(\alpha)})$ is $c$-soft. 

So fix $\alpha \in Y$ and by Proposition \ref{prop obj tA}, let $Y'\in \opc (Y)$ be such that $\alpha \in Y'$ and $\alpha $ is closed in $Y'.$ Let $X'=f^{-1}(Y')$ and $f'=f_{|Y'}.$ Then $f'^{-1}(\alpha )=f^{-1}(\alpha )$ and so $F_{i|f^{-1}(\alpha )}=F_{i|f'^{-1}(\alpha )}.$ Hence the sheaf $\varinjlim (F_{i|f^{-1}(\alpha )})$ in $\mod(A_{f^{-1}(\alpha)})$ is $c$-soft if and only if the sheaf $\varinjlim (F_{i|f'^{-1}(\alpha )})$ in $\mod(A_{f'^{-1}(\alpha)})$ is $c$-soft.

Let $a\models \alpha $ a realization of $\alpha $, ${\mathbb S}$ a prime model of the  first-order theory of ${\mathbb  M}$ over $\{a\}\cup M$ and  $r_{|}: (f'^{\mathbb S})^{-1}(a)\to f'^{-1}(\alpha )$  the homeomorphism  of Lemma \ref{lem bera}. Then the sheaf $\varinjlim (F_{i|f'^{-1}(\alpha )})$ in $\mod(A_{f'^{-1}(\alpha)})$ is $c$-soft if and only if the sheaf $\varinjlim (r^{-1}F_{i|(f'^{{\mathbb S}})^{-1}(a )})$ in $\mod(A_{(f'^{{\mathbb S}})^{-1}(a)})$ is $c$-soft.

But $(f'^{{\mathbb S}})^{-1}(a )$ is an object of $\tDfs$ and $c$ is a normal and constructible family of supports on $(f'^{{\mathbb S}})^{-1}(a )$ because $f':X'\to Y'$ is a morphism in $\tbA$ and by Proposition  \ref{prop a2 c is normal}, $c$ is a normal and constructible family of supports on $f'^{-1}(\alpha ).$ Since the sheaves $F_{i|f^{-1}(\alpha )}$ in $\mod(A_{f^{-1}(\alpha)})$ are $c$-soft if and only if the sheaves $F_{i|f'^{-1}(\alpha )}$ in $\mod(A_{f'^{-1}(\alpha)})$ are  $c$-soft if and only if the sheaves $r^{-1}F_{i|(f'^{{\mathbb S}})^{-1}(a )}$ in $\mod(A_{(f'^{{\mathbb S}})^{-1}(a)})$ are $c$-soft, by \cite[Corollary 3.5]{ep1} in $\tDfs$, we conclude that the sheaf $\varinjlim (r^{-1}F_{i|(f'^{{\mathbb S}})^{-1}(a )})$ in $\mod(A_{(f'^{{\mathbb S}})^{-1}(a)})$ is $c$-soft as required.
\qed \\

Since restriction is exact and commutes with $\otimes $, by the method of Remark \ref{nrmk met} in combination with  \cite[Proposition 3.13]{ep1} and Remark \ref{nrmk k field} we have:

\begin{prop}\label{prop f-soft and tensor}
Let $f:X\to Y$ be a morphism in $\tbA$. If  $G\in \mod(A_{X})$ is  $f$-soft and flat, then  for every   $ F\in \mod(A_{X})$ we have that  $ G\otimes  F$ is $f$-soft.
\end{prop}

By Relative fiber formula (Proposition \ref{prop f! and c}) we have:

\begin{prop}\label{prop gf soft to g soft}
Let $f:X \to Y$ and  $g:Y \to Z$ be  morphisms in $\tbA$ and let $F$ be a sheaf  in $\mod(A_{X})$. If  $F$ is $(g \circ f)$-soft, then $f_{\dex }F$ is $g$-soft.
\end{prop}

\pf
Let $\beta \in Z$, we shall prove that the restriction  $(f_{\dex }F)_{|\imin g(\beta)}$ is $c$-soft. Let $K \subset K'$ be elements of the family of complete supports $c$ on $g^{-1}(\beta ).$ By Relative fiber formula (Proposition \ref{prop f! and c}),
there is a commutative diagram
$$
\xymatrix{
\Gamma (K'; f_{\dex }F) \ar [d]_{} \ar[r]^-{\psi_{K'}} & \Gamma _c(f^{-1}(K'); F) \ar[d]^{} \\
\Gamma (K; f_{\dex }F)  \ar[r]^-{\psi _{K}} & \Gamma _c(f^{-1}(K); F)
}
$$
induced by the restrictions, with $\psi_{K'}$ and $\psi _{K}$ isomorphisms. Since $F$ is $(g \circ f)$-soft, the arrow on the right is surjective by \cite[Proposition 3.4 (3)]{ep1}. Therefore, the arrow on the left is also surjective.
This implies that, if $K'$ is a constructible element of the family of complete supports $c$ on $g^{-1}(\beta ),$ then $(f_{\dex}F)_{|K'}$ is soft, which implies $(f_{\dex}F)_{|\imin g(\beta)}$ $c$-soft by \cite[Proposition 3.4 (4)]{ep1}.
\qed \\





A special case which follows by  Remark \ref{nrmk f! and aX!} and \cite[Proposition 3.4]{ep1} is the following. Compare also with \cite[Chapter II, Proposition 2.5.7 (ii)]{ks1}.

\begin{nrmk}\label{nrmk  soft f! and aX!}
{\em
Let $f:X \to Y$ be a morphism in $\tDf$ and let $F$ be a sheaf  in $\mod(A_{X})$. Suppose that the family $c$ of supports on $X$ (resp. on $Y$) is such that every $C\in c$ has a neighborhood $D$ in $X$ (resp. in $Y$) such that $D\in c$. If  $F$ is $c$-soft, then $f_{\dex }F$ is $c$-soft.
}
\end{nrmk}


\begin{lem}\label{def fl f soft}
Let $f:X \to Y$ be a morphism in $\tDf$ with $X$ an open subspace of a normal space in $\tDf$ and let $F$ be a sheaf  in $\mod(A_{X})$. If $F$ is  flabby, then $F$ is $f$-soft. In particular, the full additive subcategory of $\mod (A_X)$ of   $f$-soft sheaves is cogenerating, i.e. for every $ F'\in \mod (A_X)$ there exists an  $f$-soft
$ F\in \mod (A_X)$ and an exact sequence
$0\rightarrow  F'\rightarrow  F$.

\end{lem}

\pf
Let $\alpha \in Y$ and let  $K$ be an element of the family $c$ of complete supports on $f^{-1}(\alpha )$. Then $K$ is quasi-compact in $X$ (Remark \ref{nrmk closed in qc}). Since $X$ is an open subspace of a normal space in $\tDf,$  by \cite[Lemma 3.2]{ep1}, the canonical morphism
$$
\lind {K\subseteq U}\Gamma(U; F) \to \Gamma(K; (F_{|f^{-1}(\alpha)})_{|K})
$$
where $U$ ranges through the family of open constructible  subsets of $X$, is an isomorphism. Since $F$ is flabby $\Gamma(X;F) \to \Gamma(U;F)$ is surjective. The exactness of filtrant $\Lind$ implies that $\Gamma(X;F) \to \Gamma(K;F_{|f^{-1}(\alpha)})_{|K})$ is surjective. This morphism factors through $\Gamma(f^{-1}(\alpha);F_{|f^{-1}(\alpha)})$ and the result follows. \qed \\


\begin{prop}\label{prop f-soft are f!-inj}
Let $f:X \to Y$ be a morphism in $\tbA$. Then the full additive subcategory of $\mod (A_X)$ of   $f$-soft sheaves
is $f_{\dex }$-injective, i.e.:
\begin{enumerate}
\item
For every $ F\in \mod (A_X)$ there exists an  $f$-soft
$ F'\in \mod (A_X)$ and an exact sequence
$0\rightarrow  F\rightarrow  F'$.
\item
If $0\rightarrow F'\rightarrow  F\rightarrow  F''
\rightarrow 0$ is an exact sequence in $\mod (A_X)$ and $ F'$ is $f$-soft, then
$0\to f_{\dex }F'\to f_{\dex }F \to  f_{\dex }F''\to 0$ is an exact sequence.
\item
If $0\rightarrow  F'\rightarrow  F\rightarrow  F'' \rightarrow 0$ is an exact sequence in $\mod(A_X)$ and $ F'$ and $ F$  are $f$-soft, then $ F''$ is  $f$-soft.
\end{enumerate}
\end{prop}

\pf
(1) By (A2) $X$ is homeomorphic to an open subspace of a normal space in $\tDf,$ therefore, by Lemma \ref{def fl f soft} the full additive subcategory of $\mod (A_X)$ of   $f$-soft sheaves is cogenerating.

(2) Let $\exs{F'}{F}{F''}$ be an exact sequence in $\mod (A_X)$ with $F'$ is $f$-soft. Then for every $\alpha \in Y$, the  sequence $\exs{F'_{|f^{-1}(\alpha)}}{F_{|f^{-1}(\alpha)}}{F''_{|f^{-1}(\alpha)}}$ is exact and $F'_{|f^{-1}(\alpha)}$ is a $c$-soft sheaf in $\mod (A_{f^{-1}(\alpha)})$. By the method of Remark \ref{nrmk met} in combination with  \cite[Proposition 3.7 (2)]{ep1}, $0\to \Gamma  _{c}(f^{-1}(\alpha); F'_{|f^{-1}(\alpha)})\to \Gamma  _{c}(f^{-1}(\alpha); F_{|f^{-1}(\alpha)})\to \Gamma  _{c}(f^{-1}(\alpha); F''_{|f^{-1}(\alpha)})\to 0$ is an exact sequence for every $\alpha \in Y$. By the Fiber formula (Corollary \ref{cor fibers of f!}),  
the sequence $\exs{f_{\dex }F'}{f_{\dex }F}{f_{\dex }F''}$ is exact.

(3) Let $\exs{F'}{F}{F''}$ be an exact sequence in $\mod (A_X)$ with $ F'$ and $F$ both  $f$-soft. Then for every $\alpha \in Y$, the  sequence $\exs{F'_{|f^{-1}(\alpha)}}{F_{|f^{-1}(\alpha)}}{F''_{|f^{-1}(\alpha)}}$ is exact and $F'_{|f^{-1}(\alpha)}$ and $F_{|f^{-1}(\alpha)}$ are $c$-soft sheaves in $\mod (A_{f^{-1}(\alpha)})$. By the method of Remark \ref{nrmk met} in combination with  \cite[Proposition 3.7 (3)]{ep1},  $F''_{|f^{-1}(\alpha)}$ is $c$-soft. Since $\alpha $ was arbitrary, $F''$ is $f$-soft.
\qed \\


\end{subsection}

\begin{subsection}{Projection and base change formulas}\label{subsection proj base change}
Here we prove the projection and base change formulas.\\

For our next proposition  we need the following  topological lemma which is an adaptation of
\cite[Lemma 2.5.12]{ks1}:

\begin{lem}\label{lem qc and tensor ks1}
Let $X$ be a quasi-compact topological space with a basis of open quasi-compact neighborhoods and let  $\Phi $  be a family of supports on $X.$ Let  $F\in \mod (A_X)$ and let  $M$  be a flat $A$-module . Then there is a natural isomorphism.
$$\Gamma _{\Phi }(X;F)\otimes M\simeq \Gamma _{\Phi }(X;F\otimes M_X).$$
In particular, if $F$ is $\Phi $-soft, then $F\otimes M_X$ is $\Phi $-soft.
\end{lem}

\pf
We may assume that $X\in \Phi.$ If $X=\cup _jU_j$ is a finite cover of $X$ by open quasi-compact neighborhoods, then
$$0\to \Gamma (X;F)\stackrel{\lambda }\to \oplus _j\Gamma (U_j;F)\stackrel{\mu }\to \oplus _{j, i}\Gamma (U_j\cap U_i;F)$$
is exact. Applying the exact functor $\bullet \otimes M,$ recall $M$ is flat, we obtain the commutative diagram with exact rows:
{\tiny
$$
\xymatrix{
0\ar[r]  & \Gamma (X;F)\otimes M \ar[r]^{\lambda } \ar[d]^{\varphi } & \oplus _j\Gamma (U_j;F)\otimes M \ar[r]^{\mu } \ar[d]^{\psi } &  \oplus _{j, i}\Gamma (U_j\cap U_i;F)\otimes M \ar[d]^{\vartheta }\\
0\ar[r] & \Gamma (X;F\otimes M_X) \ar[r]^{\lambda '} & \oplus _j \Gamma (U_j;F\otimes M_X) \ar[r]^{\mu '} & \oplus _{j, i}  \Gamma (U_j\cap U_i;F\otimes M_X)
}
$$
}
Observe also that for every $x\in X$ we have
$$ \lind {U}(\Gamma (U;F)\otimes M)\simeq \lind {U}\Gamma (U;F\otimes M_X),$$
where $U$ ranges through the family of  open quasi-compact neighborhoods of $x.$ In fact, both sides of this auxiliary isomorphism are isomorphic to $F_x\otimes M.$ Thus, if $s\in \Gamma (X;F)\otimes  M$ is such that $\varphi (s)=0$, then we can find, by the auxiliary isomorphism and quasi-compactness of $X$,  a finite covering  $X=\cup _jU_j$ by open quasi-compact neighborhoods such that $\lambda (s)=0.$ Therefore, $s=0$ and $\varphi $ is injective. If we apply the same argument to $U_j$ and $U_j\cap U_i$ instead of $X$ we see that $\psi $ and $\vartheta $ are also injective.

To show that $\varphi $ is surjective, take $t\in \Gamma (X;F\otimes M_X).$ By the auxiliary isomorphism above, there exists a finite covering $X=\cup _jU_j$ by open quasi-compact neighborhoods such that $\lambda '(t)$ is in the image of $\psi .$ But by injectivity of $\vartheta $ it follows that $t$ is in the image of $\varphi .$
\qed \\

By the Fiber formula (Corollary \ref{cor fibers of f!}) and  Lemma \ref{lem qc and tensor ks1}  we have:

\begin{prop}[Projection formula]\label{prop proj formula}
Let $f:X \to Y$ be a morphism in $\tbA$, $F \in \mod(A_{X})$ and  $G \in \mod(A_{Y})$. If $G$ is flat, then the natural morphism
$$f_{\dex }F \otimes G \to f_{\dex }(F \otimes \imin fG)$$
is an isomorphism.
\end{prop}

\pf
The morphism $\imin f\circ f_{\dex } \to \imin f\circ f_* \to \id$ induces the morphism $\imin f(f_{\dex }F \otimes G) \simeq \imin f\circ f_{\dex }F \otimes \imin f G \to F \otimes \imin f G$ and we obtain the morphism $f_{\dex }F \otimes G \to f_{\dex }(F \otimes \imin fG)$ by adjunction.

To prove that it is an isomorphism, let $\alpha \in Y$. Then
\begin{eqnarray*}
(f_{\dex }(F \otimes \imin fG))_\alpha & \simeq & \Gamma_c(\imin f(\alpha); (F \otimes \imin fG)_{|\imin f(\alpha)}) \\
& \simeq & \Gamma_c(\imin f(\alpha);F _{|\imin f(\alpha)}\otimes (\imin f G)_{|\imin f(\alpha)}) \\
& \simeq & \Gamma_c(\imin f(\alpha);F _{|\imin f(\alpha)}\otimes G_\alpha ) \\
& \simeq & \Gamma_c(\imin f(\alpha);F_{|\imin f(\alpha)}) \otimes G_\alpha \\
& \simeq & (f_{\dex }F)_\alpha \otimes G_\alpha \\
& \simeq & (f_{\dex }F \otimes G)_\alpha,
\end{eqnarray*}
by the Fiber formula (Corollary \ref{cor fibers of f!}),  Lemma \ref{lem qc and tensor ks1} and using also the  fact that $(\imin fG)_{|\imin f(\alpha)} \simeq G_\alpha $. \qed \\



We now proceed to the proof of the base change formula. By Lemma \ref{lem bera}, we have:

\begin{nrmk}\label{nrmk rest on fiber}
Let $f:X\to Y$ be a morphism in  $\tDf$ and  $\alpha \in Y$. Let $a\models \alpha $ a realization of $\alpha $ and ${\mathbb S}$ a prime model of the theory of ${\mathbb M}$ over $\{a\}\cup M.$ Then we have a  commutative diagram
$$
\xymatrix{
(f^{{\mathbb S}})^{-1}(a) \ar@{^{(}->}[r] \ar[d]^{r_|} & X({\mathbb S}) \ar[d]^{r}\ar[r]^{f^{{\mathbb S}}}  & Y({\mathbb S}) \ar[d]^{r} \\
f^{-1}(\alpha ) \ar@{^{(}->}[r] &X \ar[r]^{f} &Y
}
$$
with the restriction $r_|$ a homeomorphism  (Lemma \ref{lem bera}).

If  $F \in \mod(A_{X})$ then
$$F({\mathbb S})_{|(f^{{\mathbb S}})^{-1}(a)}= (r_{|})^{-1} F_{|f^{-1}(\alpha )}$$
where here and below we use the notation $F({\mathbb S})=r^{-1}F.$

\end{nrmk}




\begin{lem}\label{lem homeo fibers s cart sq}
Consider a cartesian square in $\tDf$
$$
\xymatrix{X' \ar[r]^{f'} \ar[d]^{g'} & Y' \ar[d]^g \\
X \ar[r]^f &Y.}
$$
Let $\gamma \in Y'$, $v\models \gamma $  a realization of $\gamma $ and ${\mathbb S}$ a prime model of the first-order theory of ${\mathbb M}$ over $\{v\}\cup M.$ Let  $u=g^{{\mathbb S}}(v)$ (so $u\models g(\gamma )$ is a realization of $g(\gamma )$). Then we have a commutative diagram
$$
\xymatrix{
(f'^{{\mathbb S}})^{-1}(v)  \ar@{^{(}->}[r]  \ar[d]^{g'^{{\mathbb S}}_{|}} & X'({\mathbb S}) \ar[r]^{f'^{{\mathbb S}}} \ar[d]^{g'^{{\mathbb S}}} & Y'({\mathbb S}) \ar[d]^{g^{{\mathbb S}}} \\
(f^{{\mathbb S}})^{-1}(u)   \ar@{^{(}->}[r] & X({\mathbb S}) \ar[r]^{f^{{\mathbb S}}} &Y({\mathbb S})}
$$
in $\tDfs$ with  $g'^{{\mathbb S}}: (f'^{{\mathbb S}})^{-1}(v) \to (f^{{\mathbb S}})^{-1}(u)$ an homeomorphism.

\end{lem}

\pf
Indeed, after removing the tilde, we have a similar commutative diagram of continuous  ${\mathbb S}$-definable maps between ${\mathbb S}$-definable spaces with the restriction to the ${\mathbb S}$-definable fibers an ${\mathbb S}$-definable homeomorphism. \qed \\

By Lemmas \ref{lem bera} and \ref{lem homeo fibers s cart sq} we have:

\begin{lem}\label{lem homeo fibers not real cart sq}
Consider a cartesian square in $\tDf$
$$
\xymatrix{X' \ar[r]^{f'} \ar[d]^{g'} & Y' \ar[d]^g \\
X \ar[r]^f &Y.}
$$
Let $\gamma \in Y'$, $v\models \gamma $  a realization of $\gamma $ and ${\mathbb S}$ a prime model of the first-order theory of ${\mathbb M}$ over $\{v\}\cup M.$ Let  $u=g^{{\mathbb S}}(v)$ (so $u\models g(\gamma )$ is a realization of $g(\gamma )$) and let ${\mathbb K}$ is a prime model of the first-order theory of ${\mathbb M}$ over $\{u\}\cup M.$  If ${\mathbb K}={\mathbb S}$, then we have a homeomorphism
$$g'_{|}:f'^{-1}(\gamma ) \to f^{-1}(g(\gamma ))$$
which induces an isomorphism
$$\Gamma_c(f'^{-1}(\gamma ); (g'^{-1}F)_{|f'^{-1}(\gamma )})\simeq \Gamma_c(f^{-1}(g(\gamma )) ;F_{|f^{-1}(g(\gamma )) })$$
for every $F \in \mod(A_{X})$.
\end{lem}

\pf
We have the following commutative diagram
$$
\xymatrix{
(f'^{{\mathbb S}})^{-1}(v)  \ar@{^{(}->}[r]  \ar[d]^{g'^{{\mathbb S}}_{|}} & X'({\mathbb S}) \ar[r]^{r} \ar[d]^{g'^{{\mathbb S}}} & X' \ar[d]^{g'}  & \ar@{_{(}->}[l]  \ar[d]^{g'_{|}} f'^{-1}(\gamma )\\
(f^{{\mathbb S}})^{-1}(u)   \ar@{^{(}->}[r] & X({\mathbb S}) \ar[r]^{r} & X &  \ar@{_{(}->}[l] f^{-1}(g(\gamma )).
}
$$
Since  $g'^{{\mathbb S}}_{|}: (f'^{{\mathbb S}})^{-1}(v) \to (f^{{\mathbb S}})^{-1}(u)$ is a homeomorphism by Lemma \ref{lem homeo fibers s cart sq} and the restrictions $(r_{|})^{-1}:(f'^{{\mathbb S}})^{-1}(v) \to f'^{-1}(\gamma )$ and $(r_{|})^{-1}:(f^{{\mathbb S}})^{-1}(u) \to f^{-1}(g(\gamma ))$ are also homeomorphisms  (Lemma \ref{lem bera} and ${\mathbb K}={\mathbb S}$) the result follows. \qed \\

From now on until the end of the subsection we need to assume also (A3).\\

\begin{lem}\label{lem homeo fibers final cart sq}
Consider a cartesian square in $\tDf$
$$
\xymatrix{X' \ar[r]^{f'} \ar[d]^{g'} & Y' \ar[d]^g \\
X \ar[r]^f &Y.}
$$
Suppose that $f:X\to Y$ satisfies (A3).

If $\gamma \in Y'$, then there exists an isomorphism
$$\Gamma_c(f'^{-1}(\gamma ); (g'^{-1}F)_{|f'^{-1}(\gamma )})\simeq \Gamma_c(f^{-1}(g(\gamma )) ;F_{|f^{-1}(g(\gamma )) })$$
for every $F \in \mod(A_{X})$.
\end{lem}

\pf
Let $v\models \gamma $ be a realization of $\gamma $ and ${\mathbb S}$ a prime model of the first-order theory of ${\mathbb M}$ over $\{v\}\cup M.$ Set $u=g^{{\mathbb S}}(v)$ and note that $u\models g(\gamma )$ is a realization of $g(\gamma )$. Let ${\mathbb K}$ is a prime model of the first-order theory of ${\mathbb M}$ over $\{u\}\cup M.$  Thus since $u\in S$, by Fact \ref{fact prime exchange}, we have either ${\mathbb K}={\mathbb M}$ and ${\mathbb M}\neq {\mathbb S}$ or ${\mathbb K}={\mathbb S}.$
 So we proceed with the proof by considering the two cases.\\

\noindent
Case ${\mathbb K}={\mathbb M}$ and ${\mathbb M}\neq {\mathbb S}$: We have $u=g(\gamma )$. Then we have
\begin{eqnarray*}
\Gamma_c(f^{-1}(g(\gamma )); F_{|f^{-1}(g(\gamma ))}) &= &\Gamma_c(f^{-1}(u );  F_{|f^{-1}(u)})\\
&\simeq & \Gamma _{c} ((f^{{\mathbb S}})^{-1}(u) ;F({\mathbb S})_{|(f^{{\mathbb S}})^{-1}(u) })\\
&\simeq & \Gamma_{c}((f'^{{\mathbb S}})^{-1}(v) ; (g'^{{\mathbb S}})^{-1}F({\mathbb S})_{|(f'^{{\mathbb S}})^{-1}(v) })\\
&\simeq & \Gamma_c(f'^{-1}(\gamma ); (g'^{-1}F)_{|f'^{-1}(\gamma )})
\end{eqnarray*}
where the first isomorphism follows  by (A3), the second follows from Lemma \ref{lem homeo fibers s cart sq} and the third follows from Lemma \ref{lem bera} together with  Remark \ref{nrmk rest on fiber}.

\noindent
Case ${\mathbb K}={\mathbb S}$: Then by Lemma \ref{lem homeo fibers not real cart sq} we have
\begin{eqnarray*}
\Gamma_c(f^{-1}(g(\gamma )); F_{|f^{-1}(g(\gamma ))})& \simeq &  \Gamma_c(f'^{-1}(\gamma ); (g'^{-1}F)_{|f'^{-1}(\gamma )}).
\end{eqnarray*}
\qed \\

We are now ready to prove the base change formula:

\begin{prop}[Base change formula]\label{prop base change formula}
Consider a cartesian square in $\tbA$
$$
\xymatrix{X' \ar[r]^{f'} \ar[d]^{g'} & Y' \ar[d]^g \\
X \ar[r]^f &Y}
$$
and let $F \in \mod(A_{X})$. Suppose that $f:X\to Y$ satisfies (A3). Then
$$\imin g \circ f_{\dex }F \iso f'_{\dex }\circ g'{}^{-1}F.$$
\end{prop}

\pf
Let us construct the morphism $\imin g \circ f_{\dex } \to f'_{\dex } \circ g'{}^{-1}.$ We shall construct first the morphism
$$f_{\dex }\circ g'_* \to g_*\circ f'_{\dex }.$$
Let $U \in \opc(Y)$ and $G \in \mod(A_{X'})$. A section $t \in \Gamma(U;f_{\dex }\circ g'_*G)$ is defined by a section $s \in \Gamma((f \circ g')^{-1}(U); G)$ such that $\supp(s) \subset g'^{-1}(Z)$ for a closed constructible  subset $Z$ of $\imin f(U)$ such that $f_{|Z}:Z\to U$ is proper.  Since
$$
\xymatrix{
g'^{-1}(Z) \ar[r]^{f'_{|}} \ar[d]^{g'_|} & g^{-1}(U) \ar[d]^{g_|} \\
Z \ar[r]^{f_{|}} &U
}
$$
is a cartesian square in $\tDf$, $(f\circ g')^{-1}(U)=(g\circ f')^{-1}(U)$ and by Proposition \ref{prop sep and proper in s} (5), the restriction $f'_{|g'^{-1}(Z)}:g'^{-1}(Z) \to \imin g(U)$ is proper. Therefore, $s \in \Gamma((g \circ f')^{-1}(U); G)=\Gamma (f'^{-1}(g^{-1}(U)); G)$ such that $\supp(s) \subset g'^{-1}(Z)$ for a closed constructible  subset $g'^{-1}(Z)$ of $f'^{-1}(g^{-1}(U))$ such that $f'_{|g'^{-1}(Z)}:g'^{-1}(Z)\to g'^{-1}(U)$ is proper. Then $s$ defines a section of $\Gamma (U; g_*\circ f'_{\dex }G)$ and we obtain $f_{\dex }\circ g'_* \to g_*\circ f'_{\dex }.$

To construct $\imin g \circ f_{\dex } \to f'_{\dex }\circ g'{}^{-1}$ consider the morphism $f_{\dex } \to f_{\dex }\circ g'_*\circ g'{}^{-1} \to g_*\circ f'_{\dex }\circ g'{}^{-1}$, where the second arrow is induced by $f_{\dex }\circ g'_* \to g_*\circ f'_{\dex }.$ We define $\imin g \circ f_{\dex } \to f'_{\dex }\circ g'{}^{-1}$ by adjunction.

To prove that $\imin g \circ f_{\dex }F \to f'_{\dex }\circ g'{}^{-1}F$ is an isomorphism, let us take $\gamma  \in Y'$. Then by the Fiber formula (Corollary \ref{cor fibers of f!})
\begin{eqnarray*}
(\imin g\circ f_{\dex }F)_\gamma & \simeq & (f_{\dex }F)_{g(\gamma )} \\
& \simeq & \Gamma_c(\imin f(g(\gamma ));F_{|\imin f(g(\gamma ))})
\end{eqnarray*}
and
\begin{eqnarray*}
(f'_{\dex }\circ g'^{-1}F)_\gamma & \simeq & \Gamma_c(f'^{-1}(\gamma ); (g'^{-1}F)_{|f'^{-1}(\gamma )}).
\end{eqnarray*}
Therefore, we have to show that
$$\Gamma_c(\imin f(g(\gamma ));F_{|\imin f(g(\gamma ))}) \simeq  \Gamma_c(f'^{-1}(\gamma ); (g'^{-1}F)_{|f'^{-1}(\gamma )}).$$
But this is proved in Lemma \ref{lem homeo fibers final cart sq}. \qed \\

\end{subsection}

\begin{subsection}{Derived proper direct image}\label{subsection derived proper dir im}
Here we derive the proper direct image and prove the derived projection and base change formulas. As corollaries we obtain the universal coefficients formula and the K\"unneth formula.\\

Recall that $A$ is a commutative ring with unit and if  $X$ is a topological space, in particular an object of $\tDf,$ we denote by $\mod(A_X)$ the category of sheaves of $A$-modules on $X$ (called also $A$-sheaves on $X$). 

Since $\mod (A_X)$ is an abelian category we may consider its derived category 
$$\der (A_X):=\der (\mod (A_X))$$
and its full triangulated subcategories
$$\der ^*(A_X):=\der ^*(\mod (A_X))$$
where $*=-,+,b.$

Since $\mod (A_X)$ has enough injectives we may  right derive  (resp. derive) the classical left exact (resp. exact) functors 
$$R\ho_{A_X}(\bullet ,\bullet ), \,\, Rf_*,\,\, f^{-1}, \,\, (\bullet )_Z, \,\,R\Gamma _Z(X;\bullet ), \,\, R\Gamma (X; \bullet )$$
on $A$-sheaves on $X,$ where $Z\subseteq X$ is a locally closed subset. \\

\begin{nrmk}\label{nrmk  wgld of k}
{\em
In order to left derive the functor 
$$\,\, \bullet \otimes _{A_X}\bullet $$
we will need to assume that the ring $A$  has finite weak global dimension, ${\rm wgld}(A)<\infty .$  The weak global dimension of $A$ is the smallest $n$ such that every $A$-module has a flat resolution of length $n,$ equivalently, it is the smallest $n$ such that $\tor ^A_j(M,N)=0$ for any $j>n$ and any $A$-modules $M$ and $N.$ Alternatively we can assume  that $A$ has finite global homological dimension, ${\rm gld}(A)<\infty $, since ${\rm wgld}(A)\leq {\rm gld}(A).$ The global dimension of $A$ is the smallest $n$ such that every $A$-module has a projective resolution of length $n.$ See \cite[Exercises I.28 and  I.29]{ks1}.

If  ${\rm wgld}(A)<\infty ,$ then by the observations on page 110 in \cite{ks1},  if $F\in \der ^b(A_X)$ (resp. $F\in \der ^+(A_X)$), then $F$ is quasi-isomorphic to a bounded complex (resp. a complex bounded from below) of flat $A$-sheaves. Therefore, we may define the left derived functor
$$\bullet \otimes  \bullet :\der ^*(A_X)\times \der ^*(A_X)\to \der ^*(A_X)$$
with $*=-, +, b.$ \\
}
\end{nrmk}

Below we will use freely  the properties relating the above derived operations  and we refer to reader to \cite[Chapter II, Section 2.6]{ks1} for details. \\

Let $f:X\to Y$ be a morphism in $\tbA$. We are going to consider the right derived functor of proper direct image
$$Rf_{\dex }: \der ^+(A_{X})\to \der ^+ (A_Y).$$
If $F\in \der ^{+}(A_X)$  then since the $f$-soft sheaves are $f_{\dex}$-injective (Proposition \ref{prop f-soft are f!-inj}), there is a complex $F'$  of $f$-soft sheaves quasi-isomorphic to $F$ and 
$$Rf_{\dex }F\simeq f_{\dex }F'.$$
Furthermore, if $g:Y\to Z$ is another morphism in $\tbA$, then by Proposition \ref{prop gf soft to g soft},
$$R(g\circ f)_{\dex }\simeq Rg_{\dex }\circ Rf_{\dex }.$$
Note also that by Theorem \ref{thm cohomo  dim f!}, the functor $Rf_{\dex }$ induces a functor:
$$Rf_{\dex }:\der ^b(A_X) \to \der ^b(A_Y).$$\\

Deriving the projection formula (Proposition \ref{prop proj formula}) we have:

\begin{thm}[Derived projection formula]\label{thm der proj formula}
Let $f:X\to Y$ be a morphism in $\tbA$. Let $F \in \der ^+(A_{X})$ and $G \in \der^+(A_{Y}).$ Suppose that ${\rm wgld}(A)<\infty .$ Then there is a natural isomorphism
$$Rf_{\dex }F \otimes G \simeq Rf_{\dex }(F \otimes \imin f G).$$
\end{thm}

\pf
First note that, if $G\in \mod (A_Y)$ is flat, then   by Lemma \ref{lem qc and tensor ks1}, $\bullet \otimes f^{-1}G$  sends $f$-soft sheaves to $f$-soft sheaves. (Indeed, for every $\alpha \in Y$, the restriction
$(\bullet )_{|f^{-1}(\alpha )}$ commutes with $\otimes $ and $(f^{-1}G)_{|f^{-1}(\alpha )}\simeq G_{\alpha }$ which is also flat.)

Now let $F \in \der ^+(A_{X})$ and $G \in \der ^+(A_Y)$. Let $F'$ be a complex of $f$-soft sheaves quasi-isomorphic to $F$ (Lemma \ref{def fl f soft}). By Remark \ref{nrmk  wgld of k}, there exists a complex $G'$ bounded from below of flat sheaves quasi-isomorphic to $G.$ Then $F'\otimes f^{-1}G'$ is a complex of
$f$-soft sheaves quasi-isomorphic to $F\otimes f^{-1}G$. Therefore, by Proposition  \ref{prop f-soft are f!-inj},
$$Rf_{\dex }F \otimes G \simeq f_{\dex }F' \otimes G' \simeq f_{\dex }(F' \otimes \imin f G')
\simeq Rf_{\dex }(F \otimes \imin f G),$$
where the second isomorphism follows from Proposition \ref{prop proj formula}.
\qed \\

Recall that $F$ is {\it $\varphi $-acyclic} where $\varphi $ is a left exact functor if $R^k\varphi F=0$ for all $k \neq 0$. In such a situation, the full subcategory of $\varphi$-acyclic objects is $\varphi$-injective (\cite[Exercises I.19]{ks1}).\\

In order to prove the derived base change formula, we need the following lemma in which we assume (A3):

\begin{lem}\label{lem f-acy to f'-acy}
Consider a cartesian square
$$
\xymatrix{X' \ar[r]^{f'} \ar[d]^{g'} & Y' \ar[d]^g \\
X \ar[r]^f &Y}
$$
 in $\tDf$. Suppose that $f:X\to Y$  satisfies (A3). Then $g'^{-1}(\bullet )$ sends $f$-soft sheaves to $f'_{\dex } $-acyclic sheaves. \footnote{In topology it is trivial that $g'^{-1}(\bullet )$ sends $f$-soft sheaves to $f'$-soft sheaves. Here  this is not evident.}
\end{lem}

\pf
Let $F\in \mod (A_X)$ be $f$-soft. Then  $F_{|f^{-1}(\alpha )}$ is   $\Gamma _c(f^{-1}(\alpha ); \bullet )$-acyclic for every $\alpha \in Y$. We must show that  the restriction $(g'^{-1}F)_{|f'^{-1}(\gamma )}$ is $\Gamma _c(f'^{-1}(\gamma ); \bullet )$-acyclic for every $\gamma \in Y'$. 

So take $\gamma \in Y'$. Let $v\models \gamma $ be a realization of $\gamma $ and ${\mathbb S}$ a prime model of the first-order theory of ${\mathbb M}$ over $\{v\}\cup M.$ Set $u=g^{{\mathbb S}}(v)$ and note that $u\models g(\gamma )$ is a realization of $g(\gamma )$. Let ${\mathbb K}$ is a prime model of the first-order theory of ${\mathbb M}$ over $\{u\}\cup M.$  Thus since $u\in S$, by Fact \ref{fact prime exchange}, we have either ${\mathbb K}={\mathbb M}$ and ${\mathbb M}\neq {\mathbb S}$ or ${\mathbb K}={\mathbb S}.$ So we proceed with the proof by considering the two cases.\\

\noindent
Case ${\mathbb K}={\mathbb M}$ and ${\mathbb M}\neq {\mathbb S}$: We have $u=g(\gamma )$ and $F_{|f^{-1}(g(\gamma ))} = F_{|f^{-1}(u)}$. Since $F_{|f^{-1}(u)}$ is  $\Gamma _c(f^{-1}(u); \bullet )$-acyclic, by (A3), $F({\mathbb S})_{|(f^{{\mathbb S}})^{-1}(u) }$ is $\Gamma _c((f^{\mathbb S})^{-1}(u); \bullet )$-acyclic.  Since $g'^{{\mathbb S}}:(f'^{{\mathbb S}})^{-1}(v ) \to (f^{{\mathbb S}})^{-1}(u)$ is a homeomorphism (Lemma \ref{lem homeo fibers s cart sq}) and  on the other hand, $((g'^{{\mathbb S}})^{-1}F({\mathbb S}))_{|(f'^{{\mathbb S}})^{-1}(v) }=(g'^{{\mathbb S}}_|)^{-1}(F({\mathbb S})_{|(f^{{\mathbb S}})^{-1}(u)})$, we conclude that $((g'^{{\mathbb S}})^{-1}F({\mathbb S}))_{|(f'^{{\mathbb S}})^{-1}(v) }$ is $\Gamma _c((f'^{\mathbb S})^{-1}(v); \bullet )$-acyclic. As $(r_|)^{-1}: f'^{-1}(\gamma ) \to (f'^{{\mathbb S}})^{-1}(v )$ is a homeomorphism (Lemma \ref{lem bera}), by Remark \ref{nrmk rest on fiber}, we have that $(g'^{-1}F)_{|f'^{-1}(\gamma )}$  is $\Gamma _c(f'^{-1}(\gamma ); \bullet )$-acyclic. \\

\noindent
Case ${\mathbb K}={\mathbb S}$: Since $g'_|:f'^{-1}(\gamma ) \to f^{-1}(g(\gamma ))$ is  a homeomorphism  (Lemma \ref{lem homeo fibers not real cart sq}),   $(g'^{-1}F)_{|f'^{-1}(\gamma )}=(g'_|)^{-1}(F_{|f^{-1}(g(\gamma ))})$ and $F_{|f^{-1}(g(\gamma ))}$ is $\Gamma _c(f^{-1}(g(\gamma )); \bullet )$-acyclic, we have that $(g'^{-1}F)_{|f'^{-1}(\gamma )}$ is $\Gamma _c(f'^{-1}(\gamma ); \bullet )$-acyclic.
\qed \\

Let $f:X\to Y$ be a morphism in $\tbA$. 
The full additive subcategory of $\mod (A_X)$ of $f_{\dex } $-acyclic  sheaves is $f_{\dex } $-injective. Therefore, if $F\in \der ^{+}(A_X)$ and $F'$ is a complex of $f_{\dex } $-acyclic  sheaves quasi-isomorphic to $F$, then
$$Rf_{\dex }F\simeq f_{\dex }F'.$$

Deriving the base change formula (Proposition \ref{prop base change formula}) we have:

\begin{thm}[Derived base change formula]\label{thm der base change formula}
Consider a cartesian square
$$
\xymatrix{X' \ar[r]^{f'} \ar[d]^{g'} & Y' \ar[d]^g \\
X \ar[r]^f &Y}
$$
 in $\tbA$. Suppose that  $f:X\to Y$ satisfies (A3). Then there is an isomorphism in $\der ^+(A_{Y'})$, functorial in  $F \in \der ^+(A_{X})$:
$$\imin g \circ Rf_{\dex }F \simeq Rf'_{\dex }\circ g'{}^{-1}F.$$
\end{thm}

\pf
Let $F \in \der ^+(A_{X})$ and let $F'$ be a complex of $f$-soft  sheaves quasi-isomorphic to $F$ (Lemma \ref{def fl f soft}).  
By Lemma \ref{lem f-acy to f'-acy}, $g'^{-1}$ sends $f$-soft  sheaves to $f'_{\dex } $-acyclic  sheaves.  
So $g'^{-1}F'$ is a complex of $f'_{\dex } $-acyclic  sheaves quasi-isomorphic to $g'^{-1}F$. Therefore, since the full subcategory of $f$-soft sheaves is $f_{\dex}$-injective (Proposition  \ref{prop f-soft are f!-inj}) and the full additive subcategory of $f'_{\dex } $-acyclic  sheaves is $f'_{\dex } $-injective, by Proposition \ref{prop base change formula}
we have
\begin{eqnarray*}
\imin g\circ Rf_{\dex }F  & \simeq & g^{-1}\circ f_{\dex }F'\\
& \simeq & f'_{\dex }\circ g'^{-1}F' \\
& \simeq & Rf'_{\dex }\circ g'^{-1}F.
\end{eqnarray*}
\qed \\

Combining the the derived projection and base change formulas we obtain:

\begin{thm}[K\"unneth formula]\label{thm kunneth formula}
Consider a cartesian square
$$
\xymatrix{X' \ar[r]^{f'} \ar[d]^{g'} \ar@{-->} [dr]^\delta  & Y' \ar[d]^g \\
X \ar[r]^f &Y}
$$
in $\tbA$ where $\delta=f\circ g'=g\circ f'$. Suppose that  $f:X\to Y$ satisfies (A3).  Suppose that ${\rm wgld}(A)<\infty .$  There is a natural isomorphism
$$R\delta_{\dex }(g'{}^{-1}F \otimes f'^{-1}G) \simeq Rf_{\dex }F \otimes Rg_{\dex }G$$
for $F \in \der ^+(A_{X})$ and $G \in \der ^+(A_{Y'})$.
\end{thm}

\pf
Using the derived projection formula and the derived base change formula
we deduce
$$Rf'_{\dex }(g'{}^{-1}F \otimes f'^{-1}G) \simeq (Rf'_{\dex }\circ g'^{-1}F )\otimes
G \simeq (\imin g \circ Rf_{\dex }F )\otimes G.$$
Using the derived projection formula once again we find
$$Rg_{\dex }\circ Rf'_{\dex }(g'{}^{-1}F \otimes f'^{-1}G) \simeq Rg_{\dex }((\imin g \circ Rf_{\dex }F) \otimes G) \simeq Rf_{\dex }F \otimes Rg_{\dex }G$$
and the result follows since $R\delta_{\dex } \simeq Rg_{\dex } \circ Rf'_{\dex }$.
\qed \\

\end{subsection}

\begin{subsection}{A bound for the cohomology of proper direct image}\label{subsection f!dim}
Here we find a bound for the cohomology of the proper direct image. \\

The {\it cohomological dimension of $f_{\dex }$} is the smallest $n$ such that $R^kf_{\dex }F=0$ for all $k>n$ and all sheaves $F$ in $\mod (A_X)$. \\




\begin{thm}\label{thm cohomo  dim f!}
Let $f:X\to Y$ be a morphism in $\tbA$. Then the cohomological dimension of $f_{\dex }$  is bounded by  $\dim X$.
\end{thm}

\pf
Let $F \in \mod(A_X)$ and let $\alpha \in Y$. Taking a $f$-soft resolution of $F$ (Proposition \ref{prop f-soft are f!-inj}) one checks easily that $(Rf_{\dex }F)_\alpha \simeq R\Gamma_c(f^{-1}(\alpha);F_{|f^{-1}(\alpha)})$ (the Fiber formula - Corollary \ref{cor fibers of f!}).



By the method of Remark \ref{nrmk met} in combination with  \cite[Theorem 3.12]{ep1}, we see that $H^k(Rf_{\dex }F)_\alpha \simeq (R^kf_{\dex }F)_\alpha \simeq R^k\Gamma_c(f^{-1}(\alpha);F_{|f^{-1}(\alpha)})=0$ if $k > \dim X$. Since $\alpha$ was arbitrary the result follows. \qed \\


\end{subsection}

\begin{subsection}{Universal coefficients and K\"unneth formulas}\label{subsection  uc and kunneth} Here we prove the universal coefficients formula and the K\"unneth formula.\\

\begin{thm}[Universal coefficients formula]\label{cor univ coef formula}
Let $X$ be an object of  $\tDf$ such that $c$ is a normal and constructible family of supports on $X.$ Let  $M$ be a flat $A$-module. 
Then there is an isomorphism
$$H_c^*(X; M_X)\simeq H_c^*(X;A_X)\otimes M.$$
\end{thm}

\pf
Let $F'$ be a complex of $c$-soft sheaves quasi-isomorphic to $A_X$ (\cite[Proposition 3.7]{ep1}). By Lemma \ref{lem qc and tensor ks1}, $F'\otimes M_X$ is a complex of $c$-soft sheaves quasi-isomorphic to $A_X\otimes M_X\simeq M_X$. Therefore, by Lemma \ref{lem qc and tensor ks1}, 
$$R\Gamma _c(X; A_X \otimes M_X) \simeq \Gamma _c(X; F' \otimes M_X) \simeq \Gamma _c(X;F' )\otimes M\simeq R\Gamma _c(X;A_X) \otimes  M.$$
Hence, by the purely homological algebra result in \cite[Exercise I.24]{ks1}, applied to the exact bifunctor $\bullet\,\,\otimes\,\,\bullet :\mod (A)\times \mod(A)\to \mod (A),$ we have 
\begin{eqnarray*}
H^k_c(X; M_{X})& = & H^k(R\Gamma _c(X; M_X))\\
& \simeq &  H^k(R\Gamma _c(X;A_X)) \otimes  M.\\
&=&H^k_c(X; A_{X}) \otimes M.
\end{eqnarray*}
\qed \\

\begin{thm}[K\"unneth formula]\label{cor kunn coho}
Consider the cartesian square
$$
\xymatrix{X \times Y \ar[r]^{p_{Y}} \ar[d]^{p_X} \ar@{-->} [dr]^{a_{X \times Y}} & Y \ar[d]^{a_{Y}} \\
X \ar[r]^{a_X} & {\rm pt}}
$$
 in $\tbA$. Suppose that $a_X:X\to {\rm pt }$ satisfies (A3) and $A$ is a field. 
Then for any $k \in \ZZ$ there is a natural isomorphism
$$H^k_c(X \times Y; A_{X \times Y}) \simeq \bigoplus_{p+q=k}(H^q_c(X; A_{X}) \otimes H^p_c(Y; A_{Y})).$$
\end{thm}

\pf 
We have $A_{X \times Y} \simeq p_X{}^{-1}A_X \simeq p_Y{}^{-1}A_Y$ and  by Theorem \ref{thm kunneth formula} we obtain 
$$Ra_{X \times Y\dex }A_{X \times Y} \simeq Ra_{X\dex }A_X \otimes Ra_{Y\dex }A_Y.$$ 
Therefore, by the purely homological algebra result in \cite[Exercise I.24]{ks1}, applied to the exact bifunctor $\bullet\,\,\otimes\,\,\bullet :\mod (A)\times \mod(A)\to \mod (A),$ we have 
\begin{eqnarray*}
H^k_c(X\times Y; A_{X\times Y})& = & H^k(Ra_{X \times Y\dex }A_{X \times Y}) \\
&\simeq & \bigoplus _{p+q=k}H^p(Ra_{X\dex }A_X) \otimes H^q(Ra_{Y\dex }A_Y)\\
&=&\bigoplus_{p+q=k}(H^q_c(X; A_{X}) \otimes H^p_c(Y; A_{Y})).
\end{eqnarray*} \qed \\

\end{subsection}

\end{section}

\begin{section}{Poincar\'e-Verdier Duality}\label{section pvd}
In this section we prove the local and the global Verdier duality, we introduce the $A$-orientation sheaf and prove the Poincar\'e and the Alexander duality.\\

The proofs here follow in a standard way from the results already obtained in the previous section. Compare with the topological case in \cite[Chapter III, Sections 3.1 and 3.3]{ks1}.\\

\begin{subsection}{Poincar\'e-Verdier duality}\label{subsection f^!}
Here we show that the derived proper direct image functor $Rf_{\dex }$  extends to $\der (A_X)$ and both $Rf_{\dex }$ and its extension have  a right adjoint $f^{\dex }$. We then deduce the basic properties of the right adjoint $f^{\dex }$ and the local and the global Verdier duality.\\

Let $f:X\to Y$ be a morphism in $\tbA$.  Let $\mathcal{J}$ be the full additive subcategory of $\mod (A_X)$ of $f$-soft sheaves. As a consequence of the results we proved for $f_{\dex }$ and for $f_{\dex } $-acyclic sheaves   we have the following properties: \footnote{As usual in category theory (\cite{ks3}), to avoid set theoretic issues, we assumed throughout the paper that we are working in a (very big) fixed universe, so below by small sum we mean a sum of a family indexed by a set in this universe.}
\begin{equation*}\label{Brownhyp}
  \begin{cases}
    \text{$\mathcal{J}$ is cogenerating}; \\
    \text{$f_{\dex }$ has finite cohomological dimension}; \\
    \text{$\mathcal{J}$ is $f_{\dex } $-injective};\\
    \text{$\mathcal{J}$ is stable under small $\oplus$}; \\
    \text{$f_{\dex }$ commutes with small $\oplus$}.
  \end{cases}
\end{equation*}
Therefore, by \cite[Proposition 14.3.4]{ks3}  the  functor $Rf_{\dex }:\der ^{+}(A_X)\to \der ^{+}(A_Y)$
extends to a functor
$$Rf_{\dex }:\der (A_{X}) \to \der (A_{Y})$$
such that:

\begin{itemize}
\item[(i)]
for every $F\in \der (A_X)$ we have
$$Rf_{\dex }F\simeq f_{\dex }F'$$
where $F'$ is a complex of $f_{\dex } $-acyclic sheaves quasi-isomorphic to $F;$

\item[(ii)]
$Rf_{\dex }$ commutes with small $\oplus .$\\
\end{itemize}

\begin{thm}\label{thm f! dual}
Let $f:X\to Y$ be a morphism in $\tbA$. Then the  functor $Rf_{\dex }:\der (A_{X}) \to \der (A_{Y})$
admits a right adjoint
$$f^{\dex }: \der (A_{Y}) \to \der (A_{X}).$$
The functor $f^{\dex }$ will thus satisfy an isomorphism
$$\Ho  _{\der (A_Y)}(Rf_{\dex }F;G)\simeq \Ho _{\der (A_X)}(F;f^{\dex }G)$$
functorial in  $F\in \der (A_X)$ and $G\in \der (A_Y).$ Moreover, the restriction
$$f^{\dex }: \der ^+(A_{Y}) \to \der ^+(A_{X})$$
is well defined and it is the  right adjoint  to the restriction  $Rf_{\dex }: \der ^+(A_X) \to \der ^+(A_Y).$
\end{thm}



\pf
The existence of right adjoint $f^{\dex }: \der (A_{Y}) \to \der (A_{X})$ is a consequence of the Brown representability theorem (see \cite[Corollary 14.3.7]{ks3} for details).

For the second part we have to show that if $G \in \der ^+(A_Y)$, then $f^{\dex }G \in \der ^+(A_X)$.
We may assume that $G \in \der ^{\geq 0}(A_Y)$. Let $N_0$ be the  dimension of $X$. Then the cohomological dimension of $f_{\dex }$ is bounded by $N_0$ (Theorem \ref{thm cohomo dim f!}). Set $a=-N_0-1$ for short. If $F \in \der ^{\leq a}(A_{X})$, then $Rf_{\dex }F \in \der ^{\leq -1}(A_{Y})$ and
$$
0=\Ho_{\der (A_Y)}(Rf_{\dex }F,G) \simeq
\Ho_{\der (A_X)}(F,f^{\dex }G).
$$
Hence for each $F \in \der ^{\leq a}(A_{X})$ we have $\Ho_{\der (A_{X})}(F,f^{\dex }G)=0.$  In particular, if
$F=\tau^{\leq a}f^{\dex }G$, then
$$
\Ho_{\der (A_{X})}(\tau^{\leq a}f^{\dex }G,f^{\dex }G) \simeq \Ho_{\der ^{\leq a}(A_{X})}(\tau^{\leq a}f^{\dex }G,\tau^{\leq a}f^{\dex }G)=0
$$
and so $\tau^{\leq a}f^{\dex }G=0$.  This implies, by definition,  that  $f^{\dex }G \in \der ^+(A_{X})$. \qed \\

The bound of the $f$-soft dimension of $\mod(A_X)$ implies the following result:\\

\begin{prop}\label{prop general duality}
Let $f:X\to Y$ be a morphism in $\tbA$.  Let $F \in \mod(A_X)$, $G \in \mod(A_Y)$. Let $N_0$ be the dimension of $X$. Then
$$
\Ho(F,H^{-N_0}f^{\dex}G) \simeq \Ho(R^{N_0}f_{\dex}F,G).
$$
\end{prop}

\pf
Let $I^\bullet$ be a complex of injective objects quasi-isomorphic to $f^{\dex}G$. As in the proof of Theorem \ref{thm f! dual} we have $f^{\dex}G \simeq \tau^{\geq -N_0}I^\bullet$ and hence we have the exact sequence
$$
\lexs{H^{-N_0}f^{\dex}G}{I^{-N_0}}{I^{-N_0+1}}
$$
that implies the isomorphism
$$
\Ho(F,H^{-N_0}f^{\dex}G) \simeq \Ho_{\der^+(A_X)}(F,f^{\dex}G[-N_0]).
$$
On the other hand the complex $Rf_{\dex}F[N_0]$ is concentrated in negative degree and hence
$$
\Ho_{\der^+(A_Y)}(Rf_{\dex}F[N_0],G) \simeq \Ho(R^0f_{\dex}F[N_0],G) \simeq \Ho(R^{N_0}f_{\dex}F,G).
$$
Then the result follows from Theorem \ref{thm f! dual}. \qed \\

If  a morphism $f:X\to Y$  in $\tbA$ is a homeomorphism onto a locally closed subset of $Y$, then $f_{\dex }$ is (isomorphic to) the extension by zero functor. What about $f^{\dex }$?

\begin{prop}\label{closemb}
Let $f:X\to Y$ be a morphism in $\tbA$. If $f:X\to Y$ is a homeomorphism onto a locally closed subset of $Y$, then $$f^{\dex }G \simeq f^{-1} \circ \rh(A_{f(X)},G)\simeq f^{-1}\circ R\Gamma _{f(X)}(G)$$
for every  $G \in \der ^+(A_{Y})$.
In particular:
\begin{itemize}
\item
if $f:X\to Y$ is a closed immersion, then $\id \simeq  f^{\dex }\circ Rf_{\dex }$;

\item
if $f:X\to Y$ is an open immersion, then $f^{\dex }\simeq f^{-1}.$
\end{itemize}
\end{prop}

\pf
Let $G \in \der ^+(A_{Y})$ and $F \in \der ^+(A_{X})$. Then
\begin{eqnarray*}
\Ho _{\der ^{+}(A_Y)}(f_{\dex }F, G)& \simeq & \Ho _{\der ^+(A_{Y})}(f_{\dex }F, R\ho (A_{f(X)}, G))\\
& \simeq & \Ho _{\der ^+(A_{X})}(f^{-1}\circ f_{\dex }F, f^{-1}\circ R\ho (A_{f(X)}, G))\\
& \simeq & \Ho _{\der ^+(A_{X})}(F, f^{-1}\circ R\ho (A_{f(X)}, G)).
\end{eqnarray*}
Hence $f^{\dex }G \simeq f^{-1} \circ \rh(A_{f(X)},G).$

Suppose now that $f:X\to Y$ is a closed immersion. Then  $f$ is proper,  $Rf_* \simeq Rf_{\dex }$ and $f^{-1}\circ  Rf_*  \simeq \id $. On the other hand, we have the isomorphisms
 \begin{eqnarray*}
 f^{\dex }\circ Rf_*F & \simeq & \imin f \circ \rh(A_{f(X)},Rf_*F) \\
 & \simeq & \imin f  \circ Rf_*\rh(A_{f(X)},F) \\
 & \simeq & \imin f  \circ Rf_*F.
\end{eqnarray*}
Hence, $\id \simeq  f^{\dex }\circ Rf_{\dex }$.

Suppose now that $f:X\to Y$ is an open immersion. Then  $f^{-1}\circ  Rf_*  \simeq \id $ and $R\Gamma _{f(X)}  \simeq  Rf_*\circ f^{-1}.$ Hence, $f^{\dex }\simeq f^{-1}.$ \qed \\

We now prove several useful properties of the dual $f^{\dex }$ of the derived proper direct image functor $Rf_{\dex }.$

\begin{prop}\label{prop f! dual functorial}
Let $f:X\to Y$ and $g:Y\to Z$ be morphisms in $\tbA$. Then
$(g\circ f)^{\dex }\simeq f^{\dex }\circ g^{\dex }.$
\end{prop}

\pf
This follows from $R(g\circ f)_{\dex }\simeq Rg_{\dex }\circ Rf_{\dex }$ and the adjunction in Theorem \ref{thm f! dual}.
\qed \\

\begin{prop}[Dual projection formula]\label{prop der proj formula dual}
Let $f:X\to Y$ be a morphism in $\tbA$.  Let $F\in \der ^b(A_Y)$ and $G\in \der ^+ (A_Y).$ Suppose that  ${\rm wgld}(A)<\infty.$ Then we have a natural isomorphism
$$f^{\dex }\circ R\ho (F,G)\simeq R\ho (f^{-1}F, f^{\dex }G).$$
\end{prop}

\pf
This follows from the derived projection formula (Theorem \ref{thm der proj formula}), the adjunction in Theorem  \ref{thm f! dual} and the adjunction
$$\Ho _{\der ^+(A_Z)}(F\otimes H, G)\simeq \Ho _{\der ^+(A_Z)}(H, R\ho (F,G)).$$
\qed\\

\begin{prop}[Dual base change formula]\label{prop der base change formula dual}
Consider a cartesian square
$$
\xymatrix{X' \ar[r]^{f'} \ar[d]^{g'} & Y' \ar[d]^g \\
X \ar[r]^f &Y}
$$
 in $\tDf$. Suppose that $f:X\to Y$ and $f':X'\to Y'$ are in $\tbA$ and that $f:X\to Y$ satisfies (A3). Then there is an isomorphism in $\der ^+(A_{X})$, functorial in  $F \in \der ^+(A_{Y'})$:
$$f^{\dex } \circ Rg_*F \simeq Rg'_*\circ f'^{\dex }F.$$
\end{prop}

\pf
This follows from the derived base change formula (Theorem \ref{thm der base change formula}), the adjunction in Theorem \ref{thm f! dual} and the adjunction
$$\Ho _{\der ^+(A_Z)}(H, Rh_*G)\simeq \Ho _{\der ^+(A_X)}(h^{-1}H, G)$$
for every morphism $h:X\to Z$ in $\tbA.$
\qed\\

\begin{thm}[Local and global Verdier duality]\label{thm hom f! dual}
Let $f:X\to Y$ be a morphism in $\tbA$. Then for  $F \in \der ^b(A_{X})$ and $G \in \der ^+(A_{Y})$, we have the local Verdier duality
$$Rf_*\circ \rh(F,f^{\dex }G) \simeq \rh(Rf_{\dex }F,G)$$
and the global Verdier duality
$$R\Ho (F,f^{\dex }G) \simeq \Ho (Rf_{\dex }F,G).$$
\end{thm}

\pf
We obtain the morphism
$$Rf_*\circ \rh(F,f^{\dex }G) \to \rh(Rf_{\dex }F,G)$$
by composing the canonical morphism
$$Rf_*\circ \rh (F, f^{\dex }G)\to \rh (Rf_{\dex }F, Rf_{\dex }\circ f^{\dex }G)$$
 with the morphism $Rf_{\dex }\circ f^{\dex }G\to G$ obtained by adjunction (Theorem \ref{thm f! dual}).

Let $V\in \op (Y).$ Then we have
\begin{eqnarray*}
H^j(R\Gamma (V; Rf_*\circ \rh (F, f^{\dex }G)))& \simeq & \Ho _{\der ^+(A_{f^{-1}(V)})}(F_{|f^{-1}(V)}, f^{\dex }G[j]_{|f^{-1}(V)})\\
& \simeq & \Ho _{\der ^+(A_{V})}(Rf_{\dex }F_{|V}, G[j]_{|V})\\
& \simeq & H^j(R\Gamma (V;\rh (Rf_{\dex }F,G)))
\end{eqnarray*}
completing the proof of the first isomorphism. The second isomorphism is obtained from the first one by applying the functor $R\Gamma (Y;\bullet ).$
\qed \\

\end{subsection}

\begin{subsection}{Orientation and duality}\label{subsection orientation}
Here we introduce the $A$-orientation sheaf and prove the Poincar\'e and the Alexander duality theorems.\\

As before, as the reader can easily verify, in all of our previous results for the proper direct image $f_{\dex }$ of a morphism $f:X\to Y$ in $\tbA$  the assumptions (A0), (A1) and (A2) were used only to show:
\begin{itemize}
\item[(i)]
if $\alpha \in Y$ is closed, then  $c$ is a normal and constructible family of supports on $f^{-1}(\alpha );$
\item[(ii)]
fiber formula;
\item[(iii)]
the theory of $f$-soft sheaves.
\end{itemize}
For the morphism $a_X:X\to {\rm pt}$ we have that (i) holds if we assume that $c$ is a normal and constructible family of supports on $X,$ (ii) is Remark \ref{nrmk ax!} and (iii) is the theory of $c$-soft sheaves developed already in \cite[Section 3]{ep1}.\\






Let $X$ be an object of  $\tDf$ such that $c$ is a normal and constructible family of supports on $X$.  Then the functor  $Ra_{X\dex }: \der ^+(A_X) \to \der ^+(\mod (A))$ admits a right adjoint
$$a_X^{\dex }: \der (\mod (A)) \to \der (A_{X})$$
(Theorem \ref{thm f! dual}) and we  have the dual projection formula (Proposition \ref{prop der proj formula dual}), the  
local and the global Verdier duality (Theorem \ref{thm hom f! dual}) for $a_{X\dex }.$ In particular, we obtain the form of the global Verdier duality proved already in \cite{ep1} (assuming $A$ is a field), where $a_X^{\dex }A_X$ is the {\it dualizing complex}:\\

\begin{thm}[Absolute Poincar\'e duality] \label{thm global vd}
Let $X$ be an object of  $\tDf$ such that $c$ is a normal and constructible family of supports on $X$.  Then we have a natural isomorphism
$$R\Ho (F, a_X^{\dex }A)\simeq R\Ho (R\Gamma _c(X; F),A)$$
as $F$ varies through $\der ^{b}(A_X).$
\end{thm}


Let $X$ be an object of  $\tDf$ such that $c$ is a normal and constructible family of supports on $X$.  We want to define the notion of orientation on $X.$ We follow the definition for topological manifolds of dimension $n$ (\cite[page 194]{i}), however we have to impose the following condition which in the topological case is true. We say that $X$ {\it has an $A$-orientation sheaf} if
for every $U\in \opc (X)$  there exists an admissible (finite) cover $\{U_1, \ldots , U_\ell \}$ of $U$ such that for each $i$ we have
\begin{equation} \label{orientation}
H_c^p(U_i; A_X)=
\begin{cases}
A \qquad \textmd{if} \qquad p={\rm dim}X\\
\\
\,\,\,\,\,\,\,\,\,\,\,\,\,\,\,\,\,\,\,\,\,\,\,\,\,\,\,\,\,\,\,
\\
0\qquad \textmd{if} \qquad p\neq {\rm dim}X.
\end{cases}
\end{equation}
$\,$\\


\begin{thm}\label{thm poinc dual}
Let $X$ be an object of  $\tDf$ of dimension $n$ such that $c$ is a normal and constructible family of supports on $X$. Suppose that $X$ has an $A$-orientation sheaf.  Then the $A$-pre-sheaf ${\mathcal Or}_X$ on $X$ with sections
$$\Gamma(U;{\mathcal Or}_X)\simeq \Ho(H_c^{n}(U;A_X),A)$$
is an $A$-sheaf, called  the {\it $A$-orientation sheaf on $X$}, such that ${\mathcal  Or}_X\simeq H^{-n}a_X^{\dex}A.$
Moreover ${\mathcal Or}_X$ is locally constant and there is a quasi-isomorphism
$${\mathcal Or}_X[n] \simeq a_X^{\dex}A.$$
\end{thm}

\pf
Setting $f=a_X,$ $F=A_U$, $G=A$ in Proposition \ref{prop general duality} we obtain
$$\Gamma(U;H^{-n}a_X^{\dex}A)\simeq \Ho(H_c^{n}(X;A_U),A) \simeq \Ho(H_c^{n}(U;A_X),A),$$
where the second isomorphism follows from  \cite[Corollary 3.9]{ep1}.

On the other hand, \eqref{orientation},  implies $R\Gamma(U_i;a_X^{\dex}A) \simeq \Rh(R\Gamma_c(U_i;A_X),A) \simeq A[-n]$, i.e. $a_X^{\dex}A$ is concentrated in degree $-n$ and the sheaf $H^{-n}a_X^{\dex}A$ is locally isomorphic to $A_X$.

\qed \\

If $A=\ZZ$ we call ${\mathcal  Or}_X$ the {\it orientation sheaf on $X$}.\\

In particular we recover \cite[Theorem 4.11]{ep1}:

\begin{nrmk}[Poincar\'e duality in cohomology]
{\em
When $A$ is a field, setting $(\bullet)^\vee=\Ho(\bullet,A)$ we obtain:
$$H^p(X;{\mathcal Or}_X) \simeq H_c^{n-p}(X;A_X)^\vee.$$
}
\end{nrmk}

Using the pure homological algebra result \cite[Proposition VI.4.6]{i}  we also have:

\begin{cor}\label{cor pd z1}
Let $X$ be an object of  $\tDf$ of dimension $n$ such that $c$ is a normal and constructible family of supports on $X$. Suppose that $X$ has an orientation sheaf ${\mathcal Or}_X$. Then there is a short exact sequence of abelian groups:
$$
\exs{\operatorname{Ext}^1(H^{k+1}_c(X;\ZZ_X),\ZZ)}{H^{n-k}(X;{\mathcal Or}_X)}{\Ho(H^k_c(X,\ZZ_X),\ZZ)}.
$$
In particular $H^{n-k}(X;{\mathcal Or}_X) \simeq \Ho(H^k_c(X,\ZZ_X),\ZZ)$ when $H^{k+1}_c(X,\ZZ_X)$ has no torsion.\\
\end{cor}

\pf
By the  pure homological algebra result \cite[Proposition VI.4.6]{i}  we  have
$$
\exs{\operatorname{Ext}^1(H^{k+1}C^{\bullet },\ZZ)}{H^{-k}\Rh(C^{\bullet },\ZZ)}{\Ho(H^kC^{\bullet },\ZZ)}
$$
for any bounded complex $C^\bullet$ of abelian groups. 
Applying this to $C^{\bullet }=R\Gamma _c(X; \ZZ _X)$ and using $R\Ho (R\Gamma _c(X; \ZZ _X), \ZZ)\simeq R\Ho (\ZZ _X, a_X^{\dex }\ZZ)\simeq R\Gamma (X; {\mathcal Or}_X[n])$  (by Theorems \ref{thm global vd} and \ref{thm poinc dual}) the result follows.
\qed \\

\begin{cor}\label{cor pd z2}
Let $X$ be an object of  $\tDf$ of dimension $n$ such that $c$ is a normal and constructible family of supports on $X$. Suppose that $X$ has an orientation sheaf ${\mathcal Or}_X$.  Then there exists  an isomorphism
$$H^n(X;{\mathcal Or}_X)\simeq  \Ho (H^{0}_c(X;\ZZ_X), \ZZ)\simeq \ZZ ^l$$
where $l$ is the number of complete connected components of $X$.
\end{cor}

\pf
By Corollary \ref{cor pd z1} (with $k=0$) and since $H_c^0(X; \ZZ _X)= \ZZ ^l$ where $l$ is the number of complete connected components of $X$, the result follows once we show that $H_c^1(X; \ZZ _X)$ is torsion free. But  this is  \cite[Chapter I, Exercise 11 and Chapter II, Exercise 28]{b}. 
\qed \\

By an {\it $A$-orientation} we understand an isomorphism $A_X\simeq {\mathcal Or}_X$. We shall say that $X$ is {\it $A$-orientable}  if an $A$-orientation exists and {\it $A$-unorientable} in the opposite case. For $A=\ZZ$ we simply say {\it orientation}, {\it orientable} or {\it unorintable}.\\

From \cite[Theorem 3.12]{ep1} ($\dim X$ is a bound on the cohomological $c$-dimension of $X$) and Corollary \ref{cor pd z1}, arguing as in \cite[Proposition 4.13]{ep1} we have:

\begin{prop}\label{prop orient}
Let $X$ be an object of  $\tDf$ of dimension $n$ such that $c$ is a normal and constructible family of supports on $X$. Suppose that $X$ has an orientation sheaf ${\mathcal Or}_X$. Then
\begin{enumerate}
\item $H^{n}_c(X;\ZZ _X)\simeq \ZZ$ if $X$ is orientable.
\item $H^{n}_c(X;\ZZ _X)\simeq 0$ if $X$ is unorientable. \\
\end{enumerate}
\end{prop}

If is $Z$ a closed constructible subset of $X$, then setting $F=A_Z$ in Theorem \ref{thm global vd} we obtain:\\

\begin{thm}[Alexander duality]\label{thm alex dual}
Let $X$ be an object of  $\tDf$ of dimension $n$ such that $c$ is a normal and constructible family of supports on $X$. Suppose that $X$ is $A$-orientable. If $Z$ a closed constructible subset  of $X,$ then there exists  a quasi-isomorphism
$$R\Gamma_Z(X;A_X)\simeq R\Ho(R\Gamma_c(Z;A_X),A)[n].$$
\end{thm}

In particular we recover \cite[Theorem 4.14]{ep1}:

\begin{nrmk}[Alexander duality in cohomology]
{\em
When $A$ is a field, setting $(\bullet)^\vee=\Ho(\bullet,A)$ we obtain:
$$H_Z^p(X;A_X) \simeq H_c^{n-p}(Z;A_X)^\vee.$$
}
\end{nrmk}

Using the pure homological algebra result \cite[Proposition VI.4.6]{i}  we also have:

\begin{cor}\label{cor alex dual1}
Let $X$ be an object of  $\tDf$ of dimension $n$ such that $c$ is a normal and constructible family of supports on $X$. Suppose that $X$ is orientatable. Then there is a short exact sequence of abelian groups:
$$
\exs{\operatorname{Ext}^1(H^{k+1}_c(Z;\ZZ_X),\ZZ)}{H^{n-k}_Z(X;\ZZ_X)}{\Ho(H^k_c(Z,\ZZ_X),\ZZ)}.
$$
In particular $H^{n-k}_Z(X;\ZZ_X) \simeq \Ho(H^k_c(Z;\ZZ_X),\ZZ)$ when $H^{k+1}_c(Z;\ZZ_X)$ has no torsion.\\
\end{cor}

\pf
By the  pure homological algebra result \cite[Proposition VI.4.6]{i}  we  have
$$
\exs{\operatorname{Ext}^1(H^{k+1}C^{\bullet },\ZZ)}{H^{-k}\Rh(C^{\bullet },\ZZ)}{\Ho(H^kC^{\bullet },\ZZ)}
$$
for any bounded complex $C^\bullet$ of abelian groups.
 Applying this to $C^{\bullet }=R\Gamma _c(Z; \ZZ _X)$ and using $R\Ho (R\Gamma _c(Z; \ZZ _X), \ZZ)[n]\simeq  R\Gamma _Z(X; \ZZ _X)$  (by Theorem \ref{thm alex dual}) the result follows.
\qed \\

\begin{cor}\label{cor alex dual2}
Let $X$ be an object of  $\tDf$ of dimension $n$ such that $c$ is a normal and constructible family of supports on $X$. Suppose that $X$ is orientatable.  Then there exists  an isomorphism
$$H^n_Z(X;\ZZ_X)\simeq  \Ho (H^{0}_c(Z;\ZZ_X), \ZZ)\simeq \ZZ ^l$$
induced by the given orientation, where $l$ is the number of complete connected components of $Z$.
\end{cor}

\pf
By Corollary  \ref{cor alex dual1} (with $k=0$) and since $H_c^0(Z; \ZZ _X)=  \ZZ ^l$ where $l$ is the number of complete connected components of $Z$, the result follows once we show that $H_c^1(Z; \ZZ _X)$ is torsion free. But  this is  \cite[Chapter I, Exercise 11 and Chapter II, Exercise 28]{b}. 
\qed \\

\end{subsection}

\end{section}

\end{document}